\PassOptionsToPackage{dvipsnames,table}{xcolor}
\documentclass[sort&compress]{elsarticle}
\usepackage[utf8]{inputenc}
\usepackage{xcolor} %
\usepackage{amsmath,amsthm,amssymb}
\usepackage{mathtools}
\usepackage{graphicx}
\graphicspath{{./img/}}
\usepackage{enumerate}
\usepackage{subcaption}
\usepackage{hyperref}
\usepackage[percent]{overpic}

\usepackage{booktabs} %
\usepackage{paralist}%
\usepackage{siunitx}%
\usepackage{bm}
\usepackage{float}
\usepackage[ruled,vlined]{algorithm2e}
\usepackage{cleveref} %

\usepackage[acronym,nogroupskip,nonumberlist,automake]{glossaries-extra}
\setabbreviationstyle[acronym]{long-short}
\newacronym{DA}{DA}{Data Assimilation}
\newacronym{DMD}{DMD}{Dynamic Mode Decomposition}
\newacronym{POD}{POD}{Proper Orthogonal Decomposition}
\newacronym{SVD}{SVD}{Singular Value Decomposition}
\newacronym{RMSE}{RMSE}{Root Mean Squared Error}
\newacronym{PF}{PF}{Particle Filter}
\newacronym{ProjPF}{Proj-PF}{Projected Particle Filter}

\newacronym{OP-PF}{OP-PF}{Optimal Proposal Particle Filter}
\newacronym{ProjOPPF}{Proj-OP-PF}{Projected Optimal Proposal Particle Filter}
\newacronym{AUS}{AUS}{Assimilation in the Unstable Subspace}
\newacronym{LV}{LV}{Lyapunov Vectors}
\newacronym{SWE}{SWE}{Shallow Water Equations}
\newacronym{ESS}{ESS}{Effective Sample Size}
\newacronym{RESAMP}{RES\%}{Resampling Percentage}
\newacronym{XC}{XC}{Pattern Correlation Coefficient}
\newacronym{L96}{L96}{Lorenz'96 model}
\newacronym{EnKF}{EnKF}{Ensemble Kalman Filter}
\newacronym{KF}{KF}{Kalman Filter}
\newacronym{EKF}{EKF}{Extended Kalman Filter}
\newacronym{PDF}{PDF}{Probability Density Function}
\newacronym{ROM}{ROM}{Reduced Order Models}
\newacronym{4D-Var}{4D-Var}{Four-dimensional variational data assimilation}
\newglossaryentry{PDE}
{
  type=\acronymtype,
  name={PDE},
  description={partial differential equation},
  first={\glsentrydesc{PDE} (\glsentrytext{PDE})},
  plural={PDEs},
  descriptionplural={partial differential equations},
  firstplural={\glsentrydescplural{PDE} (\glsentryplural{PDE})}
}
\newglossaryentry{ODE}
{
  type=\acronymtype,
  name={ODE},
  description={ordinary differential equation},
  first={\glsentrydesc{ODE} (\glsentrytext{ODE})},
  plural={ODEs},
  descriptionplural={ordinary differential equations},
  firstplural={\glsentrydescplural{ODE} (\glsentryplural{ODE})}
}

\makeglossaries

\newcommand{\mat}[1]{\ensuremath{\mathbf{{#1}}}}
\renewcommand{\vec}[1]{\ensuremath{\bm{{#1}}}}

\newcommand{\C}{\mathbb{C}}
\newcommand{\R}{\mathbb{R}}

\newcommand{\ess}{\ensuremath{\mathrm{ESS}}}

\newcommand{\modelupdate}{\mat{F}} %
\newcommand{\state}{\vec{u}} %
\newcommand{\statetruth}{\state^{\textrm{truth}}}
\newcommand{\stateens}{\state^{\textrm{ens}}}
\newcommand{\modelerror}{\vec{\omega}} %
\newcommand{\modelerrorcovariance}{\mat{Q}}
\newcommand{\modeldimension}{M} %

\newcommand{\data}{\vec{y}} %
\newcommand{\dataerror}{\vec{\eta}} %
\newcommand{\dataerrorcovariance}{\mat{R}}
\newcommand{\dataoperator}{\mat{H}}
\newcommand{\datadimension}{D} %

\newcommand{\statereduction}{\mat{V}} %
\newcommand{\reducedmodeldimension}{M^{q}} %

\newcommand{\reducedmodelerror}{\modelerror^{q}} %
\newcommand{\reducedmodelupdate}{\mat{F}^{q}} %
\newcommand{\reducedstate}{\vec{\mathrm{v}}} %
\newcommand{\reconstructedstate}{\reducedstate^{r}} %
\newcommand{\reducedmodelerrorcovariance}{\modelerrorcovariance^q}

\newcommand{\datareduction}{\mat{U}} %
\newcommand{\reduceddata}{\vec{z}} %
\newcommand{\reduceddatadimension}{\datadimension^{q}} %
\newcommand{\reduceddataoperator}{\dataoperator^{q}}
\newcommand{\reduceddataerror}{\dataerror^{q}} %
\newcommand{\reduceddataerrorcovariance}{\dataerrorcovariance^{q}}

\newcommand{\projection}{\mat{\Pi}}

\newcommand{\timeindex}{t} %
\newcommand{\timestepobs}{\tau_{D}}
\newcommand{\timestepmod}{\tau_{M}}

\newcommand{\finaltime}{T}
\newcommand{\particleindex}{\ell} %
\newcommand{\totalparticles}{L}
\newcommand{\weight}{w}

\newcommand{\innovation}{\bm{\mathcal{I}}} %

\newcommand{\leftsingularvector}{\vec{\phi}}
\newcommand{\rightsingularvector}{\vec{\psi}}
\newcommand{\leftsingularvectors}{\mat{\Phi}}
\newcommand{\rightsingularvectors}{\mat{\Psi}}

\newcommand{\singularvalue}{\sigma}
\newcommand{\singularvalues}{\mat{\Sigma}}

\newcommand{\dmdmode}{\vec{\upsilon}}
\newcommand{\dmdmodes}{\mat{\Upsilon}}
\newcommand{\dmdoperator}{\mat{A}}
\newcommand{\dmdproblemsize}{R}
\newcommand{\dmdoperatorreduced}{\mat{A}_{\dmdproblemsize}}
\newcommand{\dmdfrequency}{\omega}
\newcommand{\dmdcoefficient}{b}
\newcommand{\dmdcoefficients}{\vec{b}}

\newcommand{\snapshotmatrix}{\mat{X}}
\newcommand{\probability}{P}

\newcommand{\nodespct}{p}

\newcommand{\pinv}{\dagger} %
\newcommand{\identity}{\mat{I}} %

\newcommand{\T}[0]{\mat{T}}

\DeclarePairedDelimiter{\norm}{\Vert}{\Vert}
\DeclarePairedDelimiter{\abs}{\vert}{\vert}
\DeclareMathOperator{\rank}{rank}
\DeclareMathOperator{\lspan}{span}

\DeclareMathOperator*{\argmin}{arg\,min}

\DeclareMathOperator{\re}{Re}
\DeclareMathOperator{\im}{Im}

\DeclareMathOperator{\rmse}{RMSE}
\DeclareMathOperator{\reducedrmse}{\rmse^{q}}

\theoremstyle{definition}

\begin{document} %
\title{Model and Data Reduction for Data Assimilation: Particle Filters Employing Projected Forecasts and Data with Application to a Shallow Water Model\tnoteref{label1}}
\tnotetext[label1]{This research was supported in part by NSF
grants DMS-1714195 and DMS-1722578 and originated as part of the AIM 2019 MCRN summer school and academic year engagement program.}

\author[focal1]{Aishah Albarakati\fnref{fn1}}
\ead{albaraaa@clarkson.edu}
\address[focal1]{Department of Mathematics, Clarkson University, Potsdam, NY 13699.}

\author[focal1]{Marko Budi\v{s}i\'{c}\fnref{fn2}}
\ead{mbudisic@clarkson.edu}

\author[focal3]{Rose Crocker\fnref{fn3}}
\ead{rose.crocker@adelaide.edu.au}
\address[focal3]{School of Mathematical Sciences, University of Adelaide, Adelaide SA 5005, Australia.}

\author[focal4]{Juniper Glass-Klaiber\fnref{fn4}}
\ead{glass22j@mtholyoke.edu}
\address[focal4]{Department of Physics, Mount Holyoke College, South Hadley, MA 01075.}

\author[focal5]{Sarah Iams\fnref{fn5}}
\ead{siams@seas.harvard.edu}
\address[focal5]{Paulson School of Engineering and Applied Sciences, Harvard University, Cambridge, MA 02138.}

\author[focal3]{John Maclean\fnref{fn5b}}
\ead{John.Maclean@adelaide.edu.au}

\author[focal6]{Noah Marshall\fnref{fn6}}
\ead{noah.marshall2@mail.mcgill.ca}
\address[focal6]{Department of Mathematics and Statistics, McGill University, Montreal, Quebec H3A 0B9, Canada.}

\author[focal7]{Colin Roberts\fnref{fn7}}
\ead{colin.roberts@rams.colostate.edu}
\address[focal7]{Department of Mathematics, Colorado State University, Ft. Collins, CO 80523.}

\author[focal8]{Erik S. Van Vleck\fnref{fn8}}
\ead{erikvv@ku.edu}
\address[focal8]{Department of Mathematics, University of Kansas, Lawrence, KS 66045.}

\begin{abstract}
  The understanding of nonlinear, high dimensional flows, e.g, atmospheric and ocean flows, is critical to address the impacts of global climate change.
\gls{DA} techniques combine physical models and observational data, often in a Bayesian framework, to predict the future state of the model and the uncertainty in this prediction.
Inherent in these systems are noise (Gaussian and non-Gaussian), nonlinearity, and high dimensionality that pose challenges to making accurate predictions.
To address these issues we investigate the use of both model and data dimension reduction based on techniques including \gls{AUS}, \gls{POD}, and \gls{DMD}.
Algorithms to take advantage of projected physical and data models may be combined with \gls{DA} techniques such as \gls{EnKF} and \gls{PF} variants.
The projected \gls{DA} techniques are developed for the optimal proposal particle filter and applied to the \gls{L96} and \gls{SWE} to test the efficacy of our techniques
  in high dimensional, nonlinear systems.
\end{abstract}


\maketitle

\glsreset{SWE}

\tableofcontents

\section{Introduction}

Several important challenges impede the development of \glsreset{DA}\gls{DA} techniques: nonlinear physical models, high dimensional models and data, and non-Gaussian posterior distributions.
Particle filters and their variants are well suited to handle nonlinearities and non-Gaussian distributions, but due to so-called filter degeneracy still struggle to make accurate predictions for high dimensional problems, especially \glspl{PDE}.
Variants of particle filters have been developed to ameliorate this difficulty, including implicit particle filters, proposal density methods, the optimal proposal,  \cite{Chorin10,Leeuwen10,Snyder2011,MTAC12}.
 A related set of works analyze the performance of particle filters: \cite{SnyderEtAl08,Leeuwen15} show that particle filter degeneracy is induced by a `curse of dimensionality' associated with the data and/or model dimension.

Our contribution in this paper is to develop an approach to particle filtering based on reduced-dimension physical and data models, employing projections of the data and model spaces. The crucial benefit of this approach is that it directly targets the filter degeneracy induced by, for example, simulating some high-dimensional \glspl{PDE}, while maintaining the Bayesian framework of the particle filter suitable for nonlinear, non-Gaussian \gls{DA}.
We build on the substantial development of Assimilation in the Unstable Subspace (\gls{AUS}) methods (see, e.g., \cite{Trevisan2004, Uboldi2005, Trevisan11, PaCaTr13, PT15}), and recent work on projected data models for particle filtering~\cite{MVV20}.
The AUS methods are based on state space projections to subspaces spanned by Lyapunov vectors corresponding to the dominant Lyapunov exponents of the system dynamics. The AUS state space projections can greatly improve \gls{KF} and variational \gls{DA} schemes: {for example, an AUS implementation of 3D-Var works efficiently in high dimensions \cite{Carrassi08}.} However for ensemble-based \gls{DA} schemes, including the particle filter, the AUS approach has limited promise as the ensemble forecast already tends to align with the dominant Lyapunov exponents~\cite{Bocquet2017}.
Here we extend the AUS approach by developing projections based on other common model reduction techniques, such as \glsreset{POD}\gls{POD} and \glsreset{DMD}\gls{DMD}, and by considering combinations of projected physical and data models.
Although our focus here is on particle filters (both the standard and optimal proposal forms), the projected models developed here are also applicable to \gls{KF} and variational techniques.
We combine dimension reduction in both the physical model and the data model and compare projections based on \gls{AUS}, \gls{POD}, and \gls{DMD}.
Using these projected models, we develop projected particle filter algorithms and apply them to the \glsreset{L96}\gls{L96} and the \glsreset{SWE}\gls{SWE}.

While the present paper is motivated in large part by techniques developed using \gls{AUS} projections, see, e.g.,
\cite{PaCaTr13,ProjShadDA,MVV20}, we are also motivated by recently developed techniques for localized particle filters.
Recent work has often focused on the issue of localization, such as~\cite{Farchi18}, and currently two localized particle filtering algorithms~\cite{Poterjoy16, Potthast19} have been applied in an operational geophysical framework.
In the localized particle filter of~\cite{Potthast19}
observations are projected onto the subspace spanned by the ensemble of model forecasts to reduce the dimension of the observations.
In~\cite{Poterjoy2016}, the authors apply a localized particle filter which reduces the number of particles needed for effective assimilation.
In this scheme, particle weights are updated locally near observations, but are preserved away from observations to mimic the covariance localization in \glsreset{EnKF}\gls{EnKF}.
Other techniques such as the \emph{Dynamically Orthogonal} formulation in~\cite{SL09,Sapsis} can be interpreted in terms of reduced order physical and data models (see~\cite{SL13,MajdaQiSapsis14,QiMajda15}).

Although the original contribution of this paper is in combining data-driven \gls{ROM} with particle filters, both \gls{POD} and \gls{DMD} have been used in concert with other data assimilation techniques.
Kalman filter assimilation with a \gls{DMD}-\gls{ROM} has been used to predict wind turbine wakes in~\cite{Iungo2015}, while in~\cite{Mehta2018} the \gls{DMD} was used to enhance a Bayesian-optimized Kalman filter to predict events in the upper atmosphere.

For medium- to high-dimensional models, POD can be used to determine the dominant energy modes and, by reducing the model to the corresponding subspace, exploit a possible low dimensional structure of the model space for use in the nonlinear filtering problem \cite{wang_proper_2020}, with \gls{EnKF}~\cite{Popov2021}, and with the \gls{4D-Var} assimilation scheme \cite{cao_reduced-order_2007,dimitriu_comparative_nodate, Du2013, Fang2009}. Likewise, POD, tensorial POD, and discrete empirical interpolation have been used in the \gls{4D-Var} scheme to reduce the state space in application to the shallow water equations \cite{stefanescu_poddeim_2015}.
In~\cite{Panda2021}, efficacy of assimilation based on merging \gls{DMD}, neural networks, and \gls{4D-Var} was evaluated on chaotic dynamical systems.

While we primarily focus on using \gls{DMD} and \gls{POD} to enhance data assimilation algorithms, it is interesting to point out that \cite{nonomura_dynamic_2018,nonomura_extended-kalman-filter-based_2019}  used data assimilation techniques, specifically the Kalman filter and its variants, to enhance the computation of \gls{DMD}, in order to reduce the contribution of system noise and lead to constructing better estimates on the eigenvalues corresponding to DMD modes.

In addition, our use of \gls{POD} and \gls{DMD} as dimension reduction techniques provides a bridge between techniques developed to assimilate coherent structures and the model reduction literature.
Methods such as those developed in~\cite{MSJ17} and~\cite{MALO18} assimilate coherent structures extracted from data, but without an explicit form for the likelihood of the coherent structures. Instead these works use likelihood-free sequential Monte Carlo methods, or an ad hoc \emph{perturbed observations} approach.
In this paper we derive an explicit likelihood that corresponds to coherent structures extracted by \gls{POD} or \gls{DMD}.

In \cref{background_sec} we present background on data assimilation techniques including ensemble Kalman filters and particle filters.
In \cref{sec:proj-pf} we formulate, using abstract projections, the projected physical and data models and their use in the context of standard and optimal proposal particle filters.
We outline in \cref{sec:tech-model-reduction} the basics of the \gls{POD}, \gls{DMD}, and \gls{AUS} model reduction techniques.
Section \ref{results_sec} contains numerical results obtained using the algorithms developed to take advantage of these projected physical and data models.
These methods are applied to two nonlinear models: \gls{L96} and \gls{SWE}.
Numerical results show the efficacy of the techniques that have been developed.
In \cref{discuss_sec} we summarize and analyze our results and outline directions for future research.

\section{Background on Data Assimilation}
\label{background_sec}

\subsection{Data Assimilation: Nonlinearity and Non-Gaussian Posteriors}\label{nonlinDA}

Data assimilation is a suite of methods commonly used for obtaining accurate estimates of states and/or parameters associated with large-scale geophysical systems such as the climate and atmosphere.
DA schemes seek to optimally combine the information contained in an \emph{observation} $\data_{\timeindex }$, informed by collected data, with that in a \emph{forecast} $\state_{\timeindex }$, given by a mathematical model of the system.
The observation and forecast naturally contain associated error, due to factors such as instrumentation error in the observations and model errors or noise in the forecast.
The key aim of any \gls{DA} scheme is to propitiously balance these different sources of error and the defining characteristic of a particular \gls{DA} scheme is how it does so~\cite{Kalnay03}.

We formulate the data assimilation problem in terms of a physical model and a data model.
Consider the discrete time stochastic model with additive noise
\begin{equation}\label{Model}
\state_{\timeindex +1}= \modelupdate_{\timeindex }\left(\state_{\timeindex }\right)+\modelerror_{\timeindex },
\end{equation}
where $\state_\timeindex \in \R^{\modeldimension}$ is the model state at time index $\timeindex$, $\modelerror_{\timeindex } \in \R^{\modeldimension}$ are Gaussian errors in the model, $\modelerror_{\timeindex } \sim \mathcal{N}(0,\modelerrorcovariance_{\timeindex })$, and the function $\modelupdate_{\timeindex }:\R^\modeldimension\rightarrow \R^\modeldimension$ is the deterministic component of the model.
The data model represents the measurement of the physical state as observations $\data_{\timeindex }$
given by
\begin{equation}\label{Data}
\data_{\timeindex } = \dataoperator\state_{\timeindex } + \dataerror_{\timeindex },
\end{equation}
where $\dataerror_{\timeindex } \in \R^{\datadimension}$ are Gaussian errors in the observation, $\dataerror_{\timeindex } \sim  \mathcal{N}\left(0,\dataerrorcovariance_{\timeindex }\right)$, while  $\dataoperator :\R^\modeldimension\rightarrow \R^\datadimension, \datadimension \leq \modeldimension$ is the matrix\footnote{In this work we focus only on linear observation operators $\dataoperator$, in order to obtain closed equations for the OP-PF algorithm in \cref{sec:OP-PF} and the reduced data models in \cref{sec:dimension-reduction}. In general the observation operator may be nonlinear, and we discuss extensions of our work to this case in \cref{discuss_sec}.} representing the observation operator.

\Cref{Data} implies, first, that the data $\data_{\timeindex }$ is generated from the true (unknown) state by $\data_{\timeindex } = \dataoperator\state^\text{truth}_{\timeindex } + \dataerror_{\timeindex }$, and second, that there is a framework to convert estimates of the state $\state_{\timeindex }$ into `data space' via  $\dataoperator\state_{\timeindex }$.
The task for the \gls{DA} algorithm is to provide the best estimate of the state $\state_{\timeindex }$ using the combined information from \cref{Model,Data}.
Complete observation ($\dataoperator=\identity$) is useful in analysis, as it reduces the \gls{DA} problem to removing the influence of observation noise.

We are interested in \gls{DA} problems in which the model \cref{Model} is sufficiently nonlinear to pose challenges for \gls{DA} algorithms that rely implicitly on assumptions of linearity, chiefly the \gls{EnKF}. In order to expand on this criterion, we now formalize the \gls{DA} problem from a Bayesian perspective.

Let us formally write our estimates of the state and data as \glspl{PDF}.
The distribution of the state evolves from time $\timeindex-1$ to time $\timeindex$, including information from the data at time $\timeindex$, according to Bayes' theorem
\begin{equation}\label{Bayes}
  \probability(\state_{\timeindex }|\state_{\timeindex-1 },\,\data_{\timeindex }) =
  \frac{\probability(\data_{\timeindex }|\state_{\timeindex })\,\probability(\state_{\timeindex }|\state_{\timeindex-1})}{\probability(\data_{\timeindex })},
  \quad\text{or}\quad
  \probability(\state_{\timeindex }|\state_{\timeindex-1 },\,\data_{\timeindex }) \propto
  \probability(\data_{\timeindex }|\state_{\timeindex })\,\probability(\state_{\timeindex }|\state_{\timeindex-1 }).
\end{equation}
The left-hand-side of this equation is the distribution we want to approximate with a \gls{DA} algorithm. It is the \emph{posterior} distribution, of the state at time $\timeindex$ given (or conditioned on) the past history of the state and also given the data collected at time $\timeindex$. The right-hand-side of \cref{Bayes} is a blueprint for the action of a \gls{DA} scheme. The \gls{PDF} $\probability(\state_\timeindex|\state_{\timeindex-1 })$, the \emph{prior} distribution, contains our understanding of the state evolution from time $\timeindex-1$ to $\timeindex$.
The second term in the numerator of \cref{Bayes} is the \emph{likelihood} $\probability(\data_{\timeindex }|\state_{\timeindex })$ and can be explicitly evaluated from \cref{Data}. The normalizing constant $\probability(\data_\timeindex)$ in the denominator of \cref{Bayes} is typically impossible to calculate directly; however the \gls{PF} we develop below implicitly calculates $\probability(\data_\timeindex)$ in a discrete setting, by normalising weights.

The formulation of \gls{DA} via \cref{Bayes} may appear incomplete (for example, we have employed without proof a formulation in which data from previous observation times is not used again at time $\timeindex$), but can be rigorously justified with some weak assumptions: \cite[\S3]{Wikle2007} provides a clear exposition for the standard \gls{DA} formulation, and \cite[e.g.]{Bocquet2014} covers formulations of \gls{DA} in which the same observations may be assimilated at multiple times.

Having set up the \gls{DA} problem, let us clarify the goal of this paper. We are specifically interested in \gls{DA} problems in which the nonlinearity of the model \cref{Model}, together with factors such as sparse observations, $\datadimension \ll \modeldimension$, results in a non-Gaussian posterior distribution in \cref{Bayes}. Standard methods such as the \gls{EnKF} are inappropriate for strongly non-Gaussian posteriors, but the \gls{PF} described in \cref{sec:PF} is suitable. The key difficulty, motivated earlier and developed explicitly for \glspl{PF} in \cref{sec:PF,sec:OP-PF}, is applying the \gls{PF} to high-dimensional models and data---and the remainder of the paper develops a suite of methods to address this difficulty.

\glsreset{PF}
\subsection{The Standard \gls{PF}}
\label{sec:PF}
Particle Filters, in their many forms, are one of the most common \gls{DA} methods for nonlinear systems, as they can be proven to reproduce the true target posterior in the limit of large particle populations~\cite{KodyBook}.
Although this property makes \gls{PF}s very useful, they also suffer from several issues, including particle degeneracy and the curse of dimensionality, which will be discussed in the context of the standard particle filter.

The simplest form of particle filter is the \emph{Sequential Importance Resampling},
 or the \emph{Standard Particle Filter}~\cite[e.g.]{Vetra2018, Leeuwen2015}.
In this method, the probability density function for the state is approximated by a weighted set of guesses at the state variable.
 The probability density function for the state at time $\timeindex-1$ is
\begin{equation}
    \label{pfDist}
    \probability(\state_{\timeindex-1 }) = \sum_{\particleindex=1}^\totalparticles w_{\timeindex-1}^l \delta\left(\state_{\timeindex-1 }-\state_{\timeindex-1 }^\particleindex\right),
\end{equation}

where each particle is a state vector $\state_{\timeindex-1 }^\particleindex \in \R^\modeldimension$, each weight $w_{\timeindex-1}^l\ge0$, the sum of all weights $\sum_{l=1}^L w_{\timeindex-1}^l = 1$ and $\delta(\cdot)$ is the Dirac-delta distribution. This density function is quite simple to interpret: it says that we have $L$ guesses for the state, and that each guess is judged to be good or bad by its weight.

The Standard PF consists of the following two steps:

\begin{itemize}
\item \textbf{ Particle Update:}

    Each particle $\state_{\timeindex-1 }^\particleindex$, $\particleindex$ from $1$ to $\totalparticles$, is updated to time $\timeindex$ by running \cref{Model}. This completes the step, but let us pause to interpret what has happened. In the Bayesian formulation of \cref{Bayes}, consider the prior distribution $\probability(\state_{\timeindex }|\state_{\timeindex-1 })$: then \cref{Model} implies that for any particular $\state_{\timeindex-1 }^\particleindex$
    $$\probability(\state_{\timeindex }^\particleindex|\state_{\timeindex-1 }^\particleindex) \sim \mathcal{N}\left(\modelupdate_{\timeindex-1 }\left(\state_{\timeindex-1 }^\particleindex\right),\, \modelerrorcovariance \right)\;.$$
    Combining this equation with \cref{pfDist}, we see the exact prior distribution is a Gaussian mixture, and updating each particle according to \cref{Model} approximates the exact prior distribution by sampling once from each Gaussian distribution in the mixture.

\item \textbf{ Weight Update:}
    $$
    w_\timeindex^\particleindex = \frac{1}{c} \exp\left[
        -\frac{1}{2}\left(\data_{\timeindex} - \dataoperator\state_{\timeindex}^{\particleindex} \right)^\top
            \dataerrorcovariance^{-1}
        \left(\data_{\timeindex} - \dataoperator\state_{\timeindex}^{\particleindex} \right)
    \right]
    w_{\timeindex-1}^\particleindex \;,
    $$
where $c$ is a normalization constant such that $\sum_{l=1}^L w_{\timeindex}^l = 1$.

As before this brief description concludes the step, but we pause to connect the step to \cref{Bayes}. The exponential in the weight update is (up to normalization) precisely the term $\probability(\data_{\timeindex }|\state_{\timeindex }^\particleindex)$ (compare \cref{Data}). The normalization constant $c$ implicitly calculates the remaining factor from \cref{Bayes}, the denominator $\probability(\data_\timeindex)$. That is, the weight update multiplies each particle by the precise terms needed to complete \cref{Bayes}.
\end{itemize}

Despite the standard PF being suitable for \gls{DA} in nonlinear systems, the algorithm suffers from several related issues.
The method performs badly if one of the particle weights approaches 1.
When this occurs, all others approach zero and so the posterior distribution is effectively being approximated by a single particle. %
This issue is addressed in the Standard PF by monitoring the \gls{ESS}, $\sum_{\particleindex=1}^\totalparticles 1/\left( \weight_\timeindex^\particleindex \right)^2$, and resampling (for which there are standard algorithms, see for example~\cite[\S3.3.2]{HandbookDA}) when the \gls{ESS} drops below some threshold to refresh the particle ensemble.

The PF's tendency to exhibit particle degeneracy is related to what is often referred to as \emph{the curse of dimensionality}, a phenomenon suffered by all importance sampling algorithms in which efficiency decreases rapidly with increasing dimension of the state space~\cite{Surace2019}.
In high-dimensional spaces, resampling is insufficient to prevent degeneracy and $\totalparticles$ must be very large to give an appropriate estimate of the posterior, decreasing the computational efficiency of the algorithm.
In fact, it has been shown that the theoretically required sample size to avoid particle degeneracy scales exponentially with the state dimension~\cite{Snyder2008}.
Several methods have been developed to reduce the sample size needed to counter particle degeneracy in high dimensions, such as the \gls{OP-PF}.

\subsection{Optimal Proposal Particle Filter (OP-PF)}
\label{sec:OP-PF}

The \gls{OP-PF} ameliorates the issue of degeneracy in \gls{PF}s by attempting to ensure all posterior particles have similar weights. The idea now is to use the data in both, the particle and weight updates---first to nudge all particles towards the observations, and then to weight them.
The `proposal' in \gls{OP-PF} refers to the distribution used to update particles from one time step to the next.
In the standard particle filter the proposal density is $\probability(\state_{\timeindex }^\particleindex|\state_{\timeindex -1}^\particleindex)\sim \mathcal{N}(\modelupdate_{\timeindex -1}(\state_{\timeindex -1}^\particleindex),\modelerrorcovariance)$, as the particles are updated using the model.
Comparatively, in the OP-PF, the proposal distribution is conditioned on $\data_{\timeindex }$, as in $\probability(\state_{\timeindex }^\particleindex | \state_{\timeindex -1}^\particleindex, \data_{\timeindex })$.
Given the additive noise of the model~\eqref{Model}, the optimal proposal update in each particle is Gaussian with $\probability(\state_{\timeindex }^\particleindex|\state_{\timeindex -1}^\particleindex,\data_{\timeindex }) \sim \mathcal{N}(\vec{m}_{\timeindex }^\particleindex,\modelerrorcovariance_{p})$, and the particle update is

\begin{align}
  \state_{\timeindex }^\particleindex  =&  \vec{m}_{\timeindex}^\particleindex + \bm{\varphi},\quad \bm{\varphi} \sim \mathcal{N}(0,\modelerrorcovariance_p)\label{OPPFpred}
\shortintertext{where}
    \modelerrorcovariance^{-1}_{p} &= \modelerrorcovariance^{-1}+\dataoperator^{\top}\dataerrorcovariance^{-1}  \dataoperator\label{opS}\\
    \vec{m}_{\timeindex }^\particleindex &= \modelupdate_{\timeindex -1}(\state_{\timeindex-1}^\particleindex) + \modelerrorcovariance_{p}\dataoperator^{\top} \dataerrorcovariance^{-1}\left(\data_{\timeindex } - \dataoperator \modelupdate_{\timeindex -1}(\state_{\timeindex -1}^\particleindex)\right)\;.\label{opM}
\end{align}

Having employed the data in the particle update, the weight update must now be adjusted so that the overall scheme obeys Bayes' law. By rearranging \cref{Bayes} it can be shown \cite[e.g.]{Snyder2011} that the $\particleindex$th particle weight drawn from the proposal distribution satisfies $\weight_{\timeindex }^\particleindex \propto \probability(\data_{\timeindex }|\state_{\timeindex -1}^\particleindex)\weight_{\timeindex -1}^\particleindex$ and is also Gaussian, so that

\begin{equation}
\label{ref15}
    w_{\timeindex }^\particleindex \propto \exp \left[-\frac{1}{2}(\innovation^\particleindex_{\timeindex })^{\top}(\dataoperator\modelerrorcovariance\dataoperator^{\top} +\dataerrorcovariance)^{-1}(\innovation_{\timeindex }^\particleindex)\right]w_{\timeindex -1}^\particleindex
\end{equation}
where $\innovation_{\timeindex }^\particleindex := \data_{\timeindex } - \dataoperator \modelupdate_{\timeindex -1}(\state_{\timeindex -1}^\particleindex)$ is the \emph{innovation vector}.
Degeneracy as characterized by a single particle or a few particles with large weight still affects the OP-PF, but less than in many other particle filters.
In \cite{Snyder15} it is shown that among all \gls{PF} techniqes that obtain $\state_{\timeindex }^\particleindex$ using $\state_{\timeindex -1}^\particleindex$ and $\data_{\timeindex }$, the `optimal proposal' has the minimum variance in the weights and suffers the least from weight degeneracy.
This was extended in \cite{VanLeeuwen18} to any \gls{PF} scheme that obtains $\state_{\timeindex }^\particleindex$ using the ensemble at the previous time step $\state_{\timeindex -1}^{1:L}$ and $\data_{\timeindex }$.
The distributions used in the \gls{OP-PF} are not always available, as the additive model error and the form of the linear observation operator are required to obtain closed forms for particle and weight update schemes.

However, in \cite{Snyder2011} it is shown that the optimal proposal requires an ensemble size $\totalparticles$ satisfying $\log(\totalparticles) \propto \modeldimension \times \datadimension$ for a linear model, or will suffer from filter degeneracy.
Thus, degeneracy is deeply connected to the dimension of model and observations, posing a fundamental stumbling block to the use of \gls{PF}s in high dimensional problems.

\section{Particle Filters with Dimension Reduction in State and Observation}\label{sec:proj-pf}

\subsection{Dimension reduction of state and observation}\label{sec:dimension-reduction}

Consider the physical model~\eqref{Model} and the data model~\eqref{Data}.
There are two routes to dimension reduction using these models, namely \emph{physical model projection} and \emph{data model projection}.
 Recall that the physical state of the system is given by $\state_{\timeindex } \in \R^\modeldimension$ and our observation data is given by $\data_{\timeindex } \in \R^\datadimension$.
As previously discussed, the issue with geophysical models is that the dimensions of the physical and data space, $\modeldimension$ and $\datadimension$ respectively, can be extremely large.
This poses a problem with data assimilation methods, for example, particle filters, due to the constant need to re-draw and re-weight the particles.
In that case, the benefit to assimilation is lost.
We cannot move forward in time more than a few steps without re-sampling the set of particles, which overwrites previously gained information about the system.
Lowering either the physical model dimension or the data model dimension helps in mitigating this problem.
Here, we focus exclusively on linear, projection-based reduction of order, which amount to techniques for selecting the dimension and coordinates for the target subspace on which the models will be projected.

Starting with dimension reduction of the state, consider a matrix \(\statereduction_{\timeindex } \in \R^{\modeldimension \times \reducedmodeldimension}\) whose columns form a time-dependent orthonormal basis (\(\statereduction_{\timeindex }^{\top}\statereduction_{\timeindex }\equiv\identity\)) for the \(\reducedmodeldimension\)-dimensional subspace of the state space.
\emph{Dimension reduction} or, simply, \emph{reduction} of the state vector \(\state_{\timeindex }\) is given by inner products with columns of \(\statereduction_{\timeindex }\), which can be interpreted as the matrix multiplication \(\statereduction_{\timeindex }^{\top} : \R^{\modeldimension} \to \R^{\reducedmodeldimension}\):
\begin{equation}
  \label{eq:projection}
  \reducedstate_{\timeindex } = \statereduction_{\timeindex }^{\top} \state_{\timeindex }, \quad \reducedstate_{\timeindex } \in \R^{\reducedmodeldimension}.
\end{equation}
Since typically \(\reducedmodeldimension \ll \modeldimension\), this operation is not invertible.

The \emph{reconstruction} \(\statereduction_{\timeindex } : \R^{\reducedmodeldimension} \to \R^{\modeldimension}\) generates the state \(\reconstructedstate_{\timeindex }\)  corresponding to \(\reducedstate_{\timeindex }\)
\begin{equation}
\reconstructedstate_{\timeindex } := \statereduction_{\timeindex } \reducedstate_{\timeindex }
=\statereduction_{\timeindex }\statereduction_{\timeindex }^{\top} \state_{\timeindex }
\label{eq:reconstruction}
\end{equation}
\(\reconstructedstate_{\timeindex }\) is an element of the full state space \(\R^{\modeldimension}\), restricted to the subspace spanned by columns of \(\statereduction_{\timeindex }\), or the $\lspan$ of \( \statereduction_{\timeindex }\).
Due to orthogonality (\(\statereduction_{\timeindex }^{\top}\statereduction_{\timeindex }=\identity\)), reducing the reconstruction recovers the reduced state
$ \statereduction_{\timeindex }^{\top}(\statereduction_{\timeindex } \reducedstate_{\timeindex }) = \reducedstate_{\timeindex },$
however computing the reconstruction after a reduction does not recover the state \(\state_{\timeindex }\) itself, since
\(\state_{\timeindex } \not = \statereduction_{\timeindex }\statereduction_{\timeindex }^{\top} \state_{\timeindex }\),
unless the \(\state_{\timeindex }\in\lspan\statereduction_{\timeindex }\) to begin with.
In general,~\eqref{eq:reconstruction} computes the element of \(\lspan\statereduction_{\timeindex }\)  nearest to  \(\state_{\timeindex }\)
\begin{equation}
\reconstructedstate_{\timeindex } = \argmin_{\vec{x}\in\lspan \statereduction_{\timeindex }} \norm{ \vec{x}- \state_{\timeindex }}_{2}.\label{eq:l2-minimizer}
\end{equation}
The transformation \(\projection_{\timeindex }^p : \R^{\modeldimension} \to \R^{\modeldimension}\) whose matrix is given by
\begin{equation}
  \label{eq:projection-matrix}
  \projection_{\timeindex }^p = \statereduction_{\timeindex } \statereduction_{\timeindex }^{\top},
\end{equation}
is the orthogonal projection onto the \(\lspan \statereduction_{\timeindex }\).
In certain applications, the projection matrix may be given first, in which case an associated reduction matrix \(\statereduction_{\timeindex }\) can be computed via \gls{SVD} of \(\projection_{\timeindex }^{p}\).

To evolve reduced states \(\reducedstate_{\timeindex }\) using the physical model, we first reconstruct the state to form \(\statereduction_{\timeindex }\reducedstate_{\timeindex }\), apply~\eqref{Model} to it, and then reduce the output using \(\statereduction_{\timeindex }^{\top}\):
\begin{equation}\label{projphys}
  \reducedstate_{\timeindex +1}
  = \statereduction_{\timeindex +1}^{\top}\left(  \modelupdate_{\timeindex }(\statereduction_{\timeindex } \reducedstate_{\timeindex }) +  \modelerror_{\timeindex }\right)
  = \statereduction_{\timeindex +1}^{\top} \modelupdate_{\timeindex }(\statereduction_{\timeindex } \reducedstate_{\timeindex }) + \statereduction_{\timeindex +1}^{\top} \modelerror_{\timeindex }
\equiv \reducedmodelupdate_{\timeindex } (\reducedstate_{\timeindex }) + \reducedmodelerror_{\timeindex }.
\end{equation}
Since the orthogonal map of Gaussian random variables results in Gaussian outputs, \(\reducedmodelerror_{\timeindex } \sim \mathcal{N}(0, \reducedmodelerrorcovariance_{\timeindex })\), where
$\reducedmodelerrorcovariance_{\timeindex } = \statereduction_{\timeindex +1}^{\top} \modelerrorcovariance_{\timeindex } \statereduction_{\timeindex +1}$.
To this evolution, in the reduced-dimensional space \(\R^{\reducedmodeldimension}\), corresponds the evolution in the full space \(\R^{\modeldimension}\):
\begin{equation}\label{Model1}
\projection^{p}_{\timeindex +1}\state_{\timeindex +1}= \projection^{p}_{\timeindex +1} \modelupdate_{\timeindex }\left(\projection^{p}_{\timeindex }\state_{\timeindex }\right)+\projection^{p}_{\timeindex +1}\modelerror_{\timeindex }.
\end{equation}

The observation space can be similarly reduced using another set of vectors \(\datareduction_{\timeindex }\).
Here, we follow~\cite{MVV20} in assuming that  \(\datareduction_{\timeindex }\) is a \(\modeldimension \times \reduceddatadimension\) matrix, that is the reduction of the dimension still acts on the state space, although we will use it to reduce the observation.
This allows for the comparison of model-based order reduction methods and data-driven order reduction methods on an equal footing.
In this way, it is possible to use the same procedure to derive both \(\datareduction_{\timeindex }\) and  \(\statereduction_{\timeindex }\), although this is not necessary.

We start by defining the reduction of the observation space as
\begin{equation}
  \label{eq:reduced-data-model}
  \reduceddata_{\timeindex } := \datareduction_{\timeindex }^{\top}(\dataoperator^{\pinv} \data_{\timeindex }),
\end{equation}
where $\dagger$ denotes the Moore-Penrose pseudoinverse. Working with \(\dataoperator^{\pinv}\data_{\timeindex }\) instead of simply \(\data_{\timeindex }\) allows for the use of \(\datareduction_{\timeindex }^{\top}\) whose input space is the state space, as explained above.
Since we assume that \(\dataoperator\) has a full row-rank, the pseudoinverse \(\dataoperator^{\pinv}\) is an injection, so no information is lost in the process.

Applying this reduction to observations of states constrained to the model-reduced subspace \(\data_{\timeindex } = \dataoperator (\statereduction_{\timeindex } \reducedstate_{\timeindex } ) + \dataerror_{\timeindex } \) results in the observation equation
\begin{equation}\label{projdata}
  \begin{aligned}
  \reduceddata_{\timeindex } &=
  \datareduction_{\timeindex }^{\top}\dataoperator^{\pinv}
  [\dataoperator (\statereduction_{\timeindex } \reducedstate_{\timeindex } ) + \dataerror_{\timeindex }]
  &= \underbrace{\datareduction_{\timeindex }^{\top}\dataoperator^{\pinv}\dataoperator \statereduction_{\timeindex }}_{\reduceddataoperator_{\timeindex }} \reducedstate_{\timeindex }  +
\underbrace{\datareduction_{\timeindex }^{\top}\dataoperator^{\pinv}\dataerror_{\timeindex }}_{\reduceddataerror_\timeindex}
\end{aligned}
\end{equation}

The composition \(\dataoperator^{\pinv}\dataoperator = \projection_{\dataoperator}\) is the projection on the row (input) space of the data operator, that is the observable subspace of the state space.
The reduced data operator \(\reduceddataoperator_{\timeindex } : \R^{\reducedmodeldimension} \to \R^{\reduceddatadimension} \),
\begin{equation}
\reduceddataoperator_{\timeindex } \coloneqq \datareduction_{\timeindex }^{\top}\projection_{\dataoperator} \statereduction_{\timeindex }\label{eq:reduced-data-operator}
\end{equation}
is therefore the composition of the model reconstruction, the projection to the observable subspace space, and the data reduction.
The reduced noise is again Gaussian \(\reduceddataerror_{\timeindex } \sim \mathcal{N}(0,\reduceddataerrorcovariance_{\timeindex })\), with \(\reduceddataerrorcovariance_{\timeindex } = \datareduction_{\timeindex }^{\top}\dataoperator^{\pinv}\dataerrorcovariance(\dataoperator^{\pinv})^{\top} \datareduction_{\timeindex }\).

The alternative to the state-space based reduction of observations is to directly employ a reduction of the observation space, via the matrix \(\datareduction_{\timeindex } \in \R^{\datadimension \times \reduceddatadimension}\), and reduction \(\datareduction_{\timeindex }^{\top}:\R^{\datadimension} \to \R^{\reduceddatadimension}\).
This avoids the need to first pull-back observations into the state space by \(\dataoperator^{\pinv}\) as in~\eqref{eq:reduced-data-model}, allowing for the reduction analogous to~\eqref{eq:projection}:
\begin{equation}
  \label{eq:data-model-direct}
  \reduceddata_{\timeindex } := \datareduction_{\timeindex }^{\top} \data_{\timeindex },
\end{equation}
leading directly to the reduced model
\begin{equation}
  \label{eq:projdata-alt}
  \begin{aligned}
  \reduceddata_{\timeindex } = \datareduction_{\timeindex }^{\top} [\dataoperator (\statereduction_{\timeindex } \reducedstate_{\timeindex } ) + \dataerror_{\timeindex }]
  = \underbrace{\datareduction_{\timeindex }^{\top} \dataoperator \statereduction_{\timeindex }}_{\reduceddataoperator_{\timeindex }} \reducedstate_{\timeindex } + \underbrace{\datareduction_{\timeindex }^{\top}\dataerror_{\timeindex }}_{\reduceddataerror_{\timeindex }},
\end{aligned}
\end{equation}
with the Gaussian reduced noise
\(\reduceddataerror_{\timeindex } \sim \mathcal{N}(0,\reduceddataerrorcovariance)\), with \(\reduceddataerrorcovariance_{\timeindex } = \datareduction_{\timeindex }^{\top}\dataerrorcovariance \datareduction_{\timeindex }\).

In summary, the original state and observation equations \cref{Model} and \cref{Data} are replaced by reduced order equations through the use of orthogonal system matrices \(\statereduction_{\timeindex }\) and \(\datareduction_{\timeindex }\),
\begin{equation}
  \reducedstate_{\timeindex +1} =
            \reducedmodelupdate_{\timeindex } (\reducedstate_{\timeindex }) + \reducedmodelerror_{\timeindex },\quad
            \reducedmodelupdate_{\timeindex }( \reducedstate ) = \statereduction_{\timeindex +1}^{\top} \modelupdate_{\timeindex }(\statereduction_{\timeindex } \reducedstate),
            \quad
            \reducedmodelerror_{\timeindex } \sim  \mathcal{N}\left(0,\reducedmodelerrorcovariance_{\timeindex }\right)
            \quad
            \reducedmodelerrorcovariance_t = \statereduction_{\timeindex +1}^{\top} \modelerrorcovariance_{\timeindex } \statereduction_{\timeindex +1}
  \label{eq:reduced-model-summary},
\end{equation}
\begin{equation}
  \reduceddata_{\timeindex } =\reduceddataoperator_{\timeindex } \reducedstate_{\timeindex }  +
 \reduceddataerror_{\timeindex },
 \quad
 \reduceddataerror_{\timeindex } \sim \mathcal{N}\left(0,\reduceddataerrorcovariance_{\timeindex }\right),
 \label{eq:reduced-data-summary}
\end{equation}
where the two options for the reduced data model are
\begin{subequations}
\begin{align}
  \reduceddataoperator_{\timeindex } &=
    \datareduction_{\timeindex }^{\top}\dataoperator^{\pinv}\dataoperator \statereduction_{\timeindex }
  ,\quad
  \reduceddataerrorcovariance_t  =
  \datareduction_{\timeindex }^{\top}\dataoperator^{\pinv}\dataerrorcovariance{(\dataoperator^{\pinv})}^{\top} \datareduction_{\timeindex } , & \text{ when } \datareduction_{\timeindex }^{\top}:\R^{\modeldimension} \to \R^{\reduceddatadimension}, \label{eq:model-based-reduction}\\
  \shortintertext{or}
  \reduceddataoperator_{\timeindex } &=
    \datareduction_{\timeindex }^{\top} \dataoperator \statereduction_{\timeindex }, \quad
  \reduceddataerrorcovariance_t  =
    \datareduction_{\timeindex }^{\top}\dataerrorcovariance \datareduction_{\timeindex }, & \text{ when } \datareduction_{\timeindex }^{\top}:\R^{\datadimension} \to \R^{\reduceddatadimension}. \label{eq:data-based-reduction}
\end{align}\label{eq:reduction-summary}
\end{subequations}
The choice between two options is decided by whether \(\datareduction_{\timeindex }^{\top}\) is developed as a model-based reduction, or an observation-based reduction.
Additionally, the projected optimal proposal particle filter is modified to take the reduced state variable $\statereduction_{\timeindex }$ as the input by lifting it to the full space and then applying the observation model to it
\begin{equation}\label{projphysdatamodel}
\data_{\timeindex } = \dataoperator \statereduction_{\timeindex } \reducedstate_{\timeindex } + \dataerror_{\timeindex }.
\end{equation}

The orthonormal bases $\statereduction_{\timeindex }$ and $\datareduction_{\timeindex }$ can be obtained from different dimension reduction techniques, in particular, \gls{POD}, \gls{DMD}, and \gls{AUS}, described in more detail in~\cref{sec:tech-model-reduction}.
Since \gls{AUS} is exclusively a model-based reduction, in this work we employ only the first alternative~\eqref{eq:model-based-reduction}, to allow for a direct comparison between AUS and the others.
Regardless of how they are computed, the equations~\eqref{eq:reduction-summary} are used to formulate projected versions of particle filters.

\glsreset{ProjPF}
\subsection{\gls{ProjPF}}

Using the formulated projected models in the previous section, we can now formulate projected versions of the standard particle filter and the optimal proposal particle filter.
The basic idea is to use either the original physical model~\eqref{Model} or the physical model-based projection only, \eqref{projphys} together with~\eqref{projphysdatamodel}, for the particle update and the full projected models~\eqref{projphys} and~\eqref{projdata} for the weight update.
We also will employ a projected resampling technique based on both the physical model and data model projections.
Let
$\reducedstate_{\timeindex }^\particleindex \coloneqq \statereduction_{\timeindex }^{\top}\state_{\timeindex }^\particleindex$ for $\particleindex=1,\dots,\totalparticles$ denote the $\particleindex$th projected particle at time $\timeindex $.

The following formulations detail the alterations to the particle update and weight update routines of the particle filter utilized in our projected particle filter:
\noindent
\begin{itemize}
\item \textbf{ Particle Update:} Use~\eqref{projphys} to form
  \begin{equation}
  \label{projPFpartupdate}
  \reducedstate_{\timeindex }^\particleindex = \reducedmodelupdate_{\timeindex -1}(\reducedstate_{\timeindex -1}^\particleindex) + \reducedmodelerror_{\timeindex },\quad \particleindex=1,\dots,\totalparticles.
  \end{equation}

  \noindent
\item \textbf{ Weight Update:} Using the projected data model~\eqref{projdata} or~\eqref{eq:reduced-data-summary},~\eqref{eq:model-based-reduction}
  \begin{equation}\label{projPFweightupdate}
  w_{\timeindex }^\particleindex \propto \exp(-\frac{1}{2}(\bm{\innovation}_{\timeindex }^\particleindex)^{\top} (\reduceddataerrorcovariance_{\timeindex })^{-1} (\bm{\innovation}_{\timeindex }^\particleindex))w_{\timeindex -1}^\particleindex,\quad \particleindex=1,\dots,\totalparticles.
  \end{equation}
  where $\bm{\innovation}_{\timeindex }^\particleindex \coloneqq \reduceddata_{\timeindex } - \reduceddataoperator_{\timeindex } \reducedstate_{\timeindex }^\particleindex$.
\end{itemize}

\glsreset{ProjOPPF}
\subsection{\gls{ProjOPPF}}
In addition to a projected particle filter, we also employ a projected OP-PF.
Accordingly, the following formulations detail the alterations to the particle update and weight update routines of the OP-PF utilized in our projected OP-PF:
\noindent
\begin{itemize}
\item \textbf{ Particle Update:} Use the optimal proposal particle update~\eqref{OPPFpred},~\eqref{opS},~\eqref{opM} applied to the projected physical model~\eqref{projphys} together with the corresponding data model~\eqref{projphysdatamodel} to form
  \begin{equation}\label{projOPPFpred}
    \reducedstate_{\timeindex }^\particleindex  =  \vec{m}_{\timeindex }^\particleindex + \bm{\varphi},\quad \bm{\varphi} \sim \mathcal{N}(0,\modelerrorcovariance_p)
  \end{equation}
  where
  \begin{align}
    \label{profopS} \modelerrorcovariance^{-1}_{p} =& (\reducedmodelerrorcovariance_{\timeindex })^{-1} +  (\dataoperator \statereduction_{\timeindex })^{\top}  \dataerrorcovariance^{-1}  (\dataoperator \statereduction_{\timeindex })  \;, \\
    \label{projopM} \vec{m}_{\timeindex }^\particleindex =& \reducedmodelupdate_{\timeindex -1}(\reducedstate_{\timeindex -1}^\particleindex) + \modelerrorcovariance_{p}(\dataoperator \statereduction_{\timeindex })^{\top} \dataerrorcovariance^{-1}\left(\data_{\timeindex } - \dataoperator \statereduction_{\timeindex } \reducedmodelupdate_{\timeindex -1}(\reducedstate_{\timeindex -1}^\particleindex)\right)\;.
  \end{align}
  Alternatively, the particle update corresponding to the unprojected physical model (\(\statereduction_{\timeindex }\equiv\identity\))  could be employed.

  \medskip
  \noindent

  \glsreset{ESS}

\item \textbf{ Weight Update:} We employ the projected physical model~\eqref{projphys} and either the state space based projected data model~\eqref{projdata} or the observation space based projected data model~\eqref{eq:reduced-data-summary},~\eqref{eq:model-based-reduction}, their corresponding covariance matrices $\reducedmodelerrorcovariance_{\timeindex } = \statereduction_{\timeindex }^{\top} \modelerrorcovariance_{\timeindex } \statereduction_{\timeindex }$ and $\reduceddataerrorcovariance_{\timeindex } = \datareduction_{\timeindex }^{\top}\dataoperator^{\pinv}\dataerrorcovariance_{\timeindex }(\dataoperator^{\pinv})^{\top}\datareduction_{\timeindex }$ or $\reduceddataerrorcovariance_{\timeindex } = \datareduction_{\timeindex }^{\top} \dataerrorcovariance_{\timeindex } \datareduction_{\timeindex }$ respectively, and the projected observation operator,
  $\reduceddataoperator_{\timeindex } = \statereduction_{\timeindex }^{\top}\dataoperator^{\pinv} \dataoperator \statereduction_{\timeindex }$ or $\reduceddataoperator_{\timeindex } = \datareduction_{\timeindex }^{\top}\dataoperator \statereduction_{\timeindex }$, respectively.
  Form the matrix
  \begin{equation}\label{profopSproj}
  \mat{Z}^q_{\timeindex }:= (\reduceddataoperator_{\timeindex })\reducedmodelerrorcovariance_{\timeindex }(\reduceddataoperator_{\timeindex })^{\top} +\reduceddataerrorcovariance_{\timeindex }
  \end{equation}
  and then update the weights as
  \begin{equation}\label{projOPPFweightupdate}
  w_{\timeindex }^\particleindex \propto \exp[-\frac{1}{2}(\bm{\innovation}_{\timeindex }^\ell)^{\top} (\mat{Z}^q_{\timeindex })^{-1}(\bm{\innovation}_{\timeindex }^\ell)]w_{\timeindex -1}^\particleindex,\quad \particleindex=1,\dots,\totalparticles,
  \end{equation}
  where $\bm{\innovation}_{\timeindex }^\particleindex \coloneqq \reduceddata_{\timeindex } - \reduceddataoperator_{\timeindex } \reducedstate_{\timeindex }^\particleindex$.
\end{itemize}
We employ an extension of the projected resampling scheme proposed in \cite{MVV20}.
When the \gls{ESS}, given by
\begin{equation}
\ess = \frac{\left(\sum_{\particleindex=1}^\totalparticles w^\particleindex\right)^2}{\sum_{\particleindex=1}^\totalparticles \left(w^\particleindex\right)^2}\label{eq:ess},
\end{equation}
falls below a threshold (e.g., $\ess < \frac{1}{2}L$), then we resample.
For a given $\alpha\in [0,1],$ noise of the following form is added to resampled particles
\begin{equation}\label{resampling}
\statereduction_{\timeindex }^{\top}[\alpha \datareduction_{\timeindex } \datareduction_{\timeindex }^{\top}+(1-\alpha)\identity]\dataerror
\end{equation}
with $\dataerror \sim \mathcal{N}(\mathbf{0},\omega \identity), $ where $\omega\geq 0$ is a tuneable parameter.
The pseudocode summary of the algorithm is given in \cref{alg:ProjOPPFalg}.

\begin{algorithm}[H]
  \glsreset{ProjOPPF}
  \caption{\gls{ProjOPPF}}
    \SetAlgoLined
    $\alpha \gets \textrm{user input}$\;
    $\omega \gets \textrm{user input}$\;
    \For{$\timeindex = 1,\dots, \finaltime$}{
        \For{$\particleindex =1,\dots,\totalparticles$}{
          $\mat{m}_{\timeindex }^\particleindex = \reducedmodelupdate_{\timeindex -1} (\reducedstate_{\timeindex -1}^\ell)+\modelerrorcovariance_p\left\{\dataoperator \statereduction_{\timeindex }^\top \dataerrorcovariance^{-1}[\data_{\timeindex }-\dataoperator \statereduction_{\timeindex }\modelupdate_{\timeindex -1} (\reducedstate_{\timeindex -1}^\ell)]\right\}$\; %
            $\reducedstate_{\timeindex }^\particleindex \gets \mat{m}_{\timeindex }^\ell + \bm{\varphi}$,\qquad  $\bm{\varphi} \sim \mathcal{N}(0,\modelerrorcovariance_p)$\;
            $w_{{\timeindex }}^\particleindex \propto \exp \left[ -\frac{1}{2} (\innovation_{\timeindex}^\particleindex)^\top (\mat{Z}_{\timeindex}^q)^{-1}(\innovation_{\timeindex}^\particleindex)\right] w_{\timeindex -1}^\particleindex$;%
        }
        \If{$\ess < \frac{1}{2} \totalparticles$ \tcp{Resample if ESS below threshold}}{
        \(\statereduction_{\timeindex}^{\top}[\alpha \datareduction_{\timeindex} \datareduction_{\timeindex}^{\top}+(1-\alpha)\identity]\dataerror\),
        \qquad
        \(\dataerror \sim \mathcal{N}(\mathbf{0},\omega \identity)\)\;
        }
    }
  \label{alg:ProjOPPFalg}
\end{algorithm}

Significant reductions in computational complexity may be achieved via the reduced state space and the observation space dimensions employed in the
\gls{ProjOPPF} algorithm developed here.
If the number of particles $\totalparticles$ is fixed, then in unprojected \gls{OP-PF} the innovation vectors $\innovation^\particleindex_\timeindex$ are multiplied by $\modelerrorcovariance_p \dataoperator^{\top}\dataerrorcovariance^{-1}$ for the particle updates and by $(\dataoperator \modelerrorcovariance \dataerrorcovariance^{\top} + \dataerrorcovariance)^{-1}$ in the weight update.
If these matrices are independent of time, then the cost becomes the cost of multiplying the $\totalparticles$ innovation vectors by these matrices, operations of order $\mathcal{O}(\modeldimension \datadimension \totalparticles)$ and $\mathcal{O}(\modeldimension^2 \totalparticles)$, respectively.
Similarly, for \gls{ProjOPPF} if the projections do not depend on time, then the computational cost of the operations is of order $\mathcal{O}(\reduceddatadimension \reducedmodeldimension \totalparticles)$ and $\mathcal{O}((\reduceddatadimension)^2 \totalparticles)$, respectively.
If the projections depend on time, as with the \gls{AUS} projection, then there is an additional cost of forming or factoring matrices, e.g., using LU or Cholesky, in addition to multiplication by the innovation vectors.

\section{Techniques for Model Reduction}\label{sec:tech-model-reduction}

Development of particle filters on a subspace of the state or observation, detailed in \cref{sec:proj-pf}, does not depend on any particular technique for computing the dimension reduction matrices \(\statereduction_{\timeindex }\) and \(\datareduction_{\timeindex }\).
In this section we outline three techniques for computing these matrices.
\gls{POD} and \gls{DMD} are data-driven (model-free) techniques that only require a set of simulation snapshots to calculate the reduction subspace, while the \gls{LV} computation requires access to derivatives of the deterministic part of the model update equation \(\modelupdate\) (see~\eqref{Model}).

\subsection{Proper Orthogonal Decomposition (POD)}\label{sec:pod}

\glsreset{POD}\gls{POD} is the model-reduction technique based on computation of a parsimonious orthogonal basis for the state space subspace occupied by a given evolution a dynamical system.
It is ubiquitous in applied mathematics, and in other contexts is known as principal component analysis, Karhunen--Lo\'eve decomposition, and empirical orthogonal function decomposition.
Here, we review only the necessary about \gls{POD}; an excellent short review of the main features with a wealth of references can be found in \cite[\S 22.4]{Tropea2007}.
In the context of fluid dynamics~\cite{Berkooz1993}, the orthogonal basis is calculated using eigenvectors of the cross-correlation matrix of the simulated data.

Given a recording of evolution of state vectors\footnote{If the data-based model reduction is needed, as in~\eqref{eq:data-based-reduction}, then evolution of observations should be used instead.} (called \emph{snapshots}) \(\state_{\timeindex } \in \R^\modeldimension\), over time \(\timeindex =1,\dots,\finaltime\), stored as a \emph{snapshot matrix}
\begin{equation}
  \snapshotmatrix :=
\begin{bmatrix}
  \state_{1} & \state_{2} & \dots & \state_{\finaltime}
\end{bmatrix},\label{eq:snapshot-matrix}
\end{equation}
\gls{POD} amounts to a separation-of-variables ansatz
\begin{equation}
 \state_{\timeindex } \approx \sum_{m=1}^{\modeldimension} \leftsingularvector_{m} \singularvalue_{m} \rightsingularvector_{\timeindex ,m}.\label{eq:POD-ansatz}
\end{equation}
Here vectors \(\leftsingularvector_{m}\), $m=1, \dots, \modeldimension$, are the normalized ``spatial'' profiles of the state \(\state_{\timeindex } \in \R^\modeldimension\), vectors \(\rightsingularvector_{\timeindex ,m} \coloneqq  \begin{bmatrix} \psi_{1,m} & \psi_{2,m} & \cdots & \psi_{\finaltime,m} \end{bmatrix}^\top \) are the normalized time evolutions, while \(\sigma_{m}\) are the linear combination coefficients, i.e., magnitudes.
While there are many possible separation-of-variable decompositions, \gls{POD} is specified by the requirement that \(\{\leftsingularvector_{m}\}_{m=1}^\modeldimension\) and \(\{\rightsingularvector_{m}\}_{m=1}^\modeldimension\) should be orthogonal sets.

In matrix notation, and over a fixed period of time, this ansatz corresponds to \glsreset{SVD}\gls{SVD} of the snapshot matrix \(\snapshotmatrix\)
\begin{equation}
  \label{eq:svd-of-u}
  \snapshotmatrix =
  \begin{bmatrix}
    \leftsingularvector_{1} &  \leftsingularvector_{2} & \dots
  \end{bmatrix}
  \begin{bsmallmatrix}
    \singularvalue_{1} &  & \\
    & \singularvalue_{2} & & \\
    & & \ddots
  \end{bsmallmatrix}
  \begin{bmatrix}
    \rightsingularvector_{1} &  \rightsingularvector_{2} & \dots
  \end{bmatrix}^{\top}
  ,
\end{equation}
The rank of the snapshot matrix \(\snapshotmatrix\) is equal to the number of nonzero singular values, \(\sigma_{m}\).
It is common to order the singular values in decreasing values, and refer to those \(\sigma_{m}\) and vectors \(\leftsingularvector_{m}\) as \emph{dominant} if they have a low index.
Furthermore, singular values that are equal to zero, and associated singular vectors, are sometimes omitted to form the ``economy'' version of \gls{SVD}.

To reduce the dimension of \(\snapshotmatrix\), while preserving the character of dynamics, the reduction matrix \(\statereduction\) (see~\eqref{eq:projection}) is formed,
\begin{equation}
  \label{eq:POD-projection}
  \mat{\statereduction}^{(r)} =
  \begin{bmatrix}
    \leftsingularvector_{1} & \cdots & \leftsingularvector_{r}
  \end{bmatrix},
\end{equation}
containing the first \(r < \modeldimension\) left singular vectors  \(\leftsingularvector_{m}\) (spatial profiles).
By the Eckart--Young theorem~\cite{Eckart1936}, the projected snapshot matrix
\(\statereduction^{(r)}{\statereduction^{(r)}}^{\top} \snapshotmatrix \)
is the best approximation of \(\snapshotmatrix\) among all matrices of rank $r$ as measured by the Frobenius norm, i.e., element-wise \(\ell^{2}\) norm~\cite[\S 2.4]{Golub2013}.

In general, the choice $r$ that obtains a parsimonious, yet usable, reduced-order approximation can be problem dependent, although there are prescriptions of optimal rank in absence of problem-dependent guidances~\cite{Gavish2014}.
Ideally, a gap or a jump in a singular value plot is an indication that there is a sharp change in the approximation error as the number of dimensions retained is changed across the gap.
In other cases, no such gap may be seen, which can be in certain cases traced to model-agnostic application of the technique \cite[\S 19.4]{Kutz}.

\begin{algorithm}[hbt!]
    \caption{Data projection using \glsreset{POD}\gls{POD}}
    \label{alg:pod}
    \SetAlgoLined\DontPrintSemicolon
    $\snapshotmatrix \gets \begin{bmatrix} \state_1 & \state_2 & \cdots & \state_T \end{bmatrix}$ \tcp{Form snapshot matrix}
    $ \snapshotmatrix = \leftsingularvectors \singularvalues \rightsingularvectors^\top$ \tcp{Singular Value Decomposition}
    $r \gets \textrm{user input}$\;
    $\statereduction_{\textrm{POD}} \gets \begin{bmatrix} \leftsingularvector_{ 1} & \leftsingularvector_{ 2} & \cdots & \leftsingularvector_{ r}\end{bmatrix}$ \tcp{Dimension reduction matrix}
  \end{algorithm}
  \subsection{Dynamic Mode Decomposition (DMD)}\label{sec:dmd}
\glsreset{DMD}\gls{DMD}~\cite{Schmid2010,Rowley2009} provides a route to order reduction by approximating the evolution of snapshots \(\state_{\timeindex}\) by a separation-of-variables ansatz in the form
\begin{equation}
  \label{eq:DMD-approximation}
  \state_{\timeindex} \approx \sum_{m=1}^{M} \dmdmode_{m} e^{\timeindex \dmdfrequency_{m} \timestepmod} \dmdcoefficient_{m},
\end{equation}
where \(\timestepmod\) is the timestep separating adjacent snapshots, \(\dmdmode_{m}\) are DMD modes, corresponding to a spatial profile of the component of dynamics, \(\dmdfrequency_{m} \in \C\) are complex-valued DMD frequencies governing growth, decay, and oscillation of time evolution, while \(b_{m} \in \C\) are linear combination coefficients.
For real- valued input vectors \(\state_{\timeindex}\), complex-valued DMD frequencies and modes come in conjugate pairs.

The primary assumption is that there exists a high-dimensional time-invariant matrix \(\dmdoperator\) that relates pairs of snapshots
\begin{equation}
  \label{eq:DMD-operator}
  \state_{\timeindex+1} = \dmdoperator \state_{\timeindex}.
\end{equation}
While this assumption does not hold exactly, the Koopman operator theory~\cite{Budisic2012Chaos, Tu2014a, Korda2018, Kutz_DMD, Mezic2013} asserts that for time-invariant systems, including~\eqref{Model}, the state  \(\state_{\timeindex}\) can be embedded in an infinite-dimensional space on which the matrix \(\dmdoperator\) is the linear infinite-dimensional Koopman operator and the equation analogous to~\eqref{eq:DMD-approximation} does hold for expected values of these state variables.

If~\eqref{eq:DMD-operator} holds and we have access to \(\dmdoperator\)  then modes \(\dmdmode\) can be computed as eigenvectors  of \(\dmdoperator\)
\begin{equation}
  \label{eq:DMD-modes}
  \dmdoperator \dmdmode_{m} = \lambda_{m} \dmdmode_{m},
\end{equation}
while frequencies \(\dmdfrequency_{m}\)  are derived from eigenvalues \(\lambda_{m}\) by the formula
\begin{equation}
  \label{eq:DMD-frequencies}
  \lambda_{m} = \exp( \dmdfrequency_{m} \timestepmod ),
\end{equation}
where \(\timestepmod\) again is the timestep that separates adjacent snapshots \(\state_{m}\).
As is well-known from linear systems analysis, \(\re \dmdfrequency_{m}\) determines the rate of exponential decay (\(\re (\dmdfrequency_{m}) < 0\)) or growth (\(\re( \dmdfrequency_{m}) > 0\)) of DMD modes, while \(\im (\dmdfrequency_{m})\) corresponds to the angular frequency of oscillations.
Modes with frequencies at the origin \(\dmdfrequency_{m} = 0 \) correspond to constant components in the evolution.

Typically, though, we do not have access to \(\dmdoperator\) directly and we need to approximate its eigenvalues and eigenvectors.
\gls{DMD} is a family of numerical procedures approximates \(\dmdmode_{m}\) and \(\dmdfrequency_{m}\) while avoiding both the explicit embedding of~\eqref{eq:DMD-operator} in the high-dimensional space, and the computation and storage of the full matrix \(\dmdoperator\), which in practice is prohibitively large, and theoretically infinite.
Here we present the so-called \emph{exact} \gls{DMD} algorithm~\cite{Tu2014a} which is commonly a starting point for the DMD analysis.

Regression of the discrete dynamics equation~\eqref{eq:DMD-operator} onto the snapshots can be solved approximately for all adjacent time steps as an \(\ell^{2}\) minimization problem
\begin{equation}
 \min_{\dmdoperator} \sum_{\timeindex = 0}^{T-1}\norm{ \state_{\timeindex+1} - \dmdoperator \state_{\timeindex}  }_{2}^{2} =  \min_{\dmdoperator} \norm{ \snapshotmatrix_{2} - \dmdoperator \snapshotmatrix_{1}}.
\label{eq:variational-DMD}
\end{equation}
The second formula is the equivalent matrix notation using involving two submatrices of the snapshot matrix, \(\snapshotmatrix_{1}\) and \(\snapshotmatrix_{2}\), formed by respectively erasing the first and the last column of \(\snapshotmatrix\).
In principle, this problem could be solved by a Moore--Penrose inverse as
\(\dmdoperator_{\text{DMD}} := \snapshotmatrix_{2} \snapshotmatrix_{1}^{\pinv}, \)
but this is generally avoided as the size of \(\dmdoperator_{\text{DMD}}\) is quadratic in the dimension of snapshot vectors, and in practice may be prohibitively large store. Furthermore, any numerical errors arising in the process could fatally affect the well-posedness of the computation.

As the ultimate goal is not the calculation of \(\dmdoperator\) but rather a (small) subset of its eigenvectors and eigenvalues~\eqref{eq:DMD-modes}, most variants of \gls{DMD} employ an order reduction step to improve numerical robustness and reduce the size of  numerical linear algebra calculations.
Here, we use a \gls{POD} based order reduction, essentially using \cref{sec:pod} as the substep of the \gls{DMD} analysis.
Compute \gls{SVD} of the ``left'' snapshot matrix
\begin{equation}
\snapshotmatrix_{1} = \leftsingularvectors\singularvalues\rightsingularvectors^{\top}
\end{equation}
and truncate the involved matrices to \(\dmdproblemsize\) dominant vectors, forming   \(\leftsingularvectors_{\dmdproblemsize}\), \(\singularvalues_{\dmdproblemsize}\), \(\rightsingularvectors_{\dmdproblemsize}\).
Since the goal of this step is not to perform the full order reduction, but merely to numerically stabilize the problem and reduce the size of numerical problems solved below, \(\dmdproblemsize\) could be fairly large, even \(\dmdproblemsize = 0.9 \modeldimension\).

Using the matrices just computed, form
\begin{equation}
  \label{eq:DMD-reduced}
  \dmdoperator_{\dmdproblemsize} = \leftsingularvectors_{\dmdproblemsize} \snapshotmatrix_{2} \rightsingularvectors_{\dmdproblemsize} \singularvalues_{\dmdproblemsize},
\end{equation}
and compute its eigenvalues and eigenvectors \(\dmdoperator_{\dmdproblemsize} \hat{\dmdmode}_{m} = \lambda_{m} \hat{\dmdmode}_{m}\).
Matrix \(\dmdoperator_{\dmdproblemsize}\) is a \(\dmdproblemsize \times \dmdproblemsize\) matrix whose eigenvalues are the same as eigenvalues of the \(\modeldimension \times\modeldimension\) matrix \(\dmdoperator\), and whose eigenvectors \(\hat{\dmdmode}_{m}\) can be used to reconstruct DMD modes by
\begin{equation}
  \label{eq:reconstructing-DMD-modes}
  \dmdmode_m = \lambda_m^{-1} \snapshotmatrix_2 \rightsingularvectors_{\dmdproblemsize} \singularvalues_{\dmdproblemsize}^{-1} \hat{\dmdmode}_m.
\end{equation}
Depending on the algorithm to compute eigenvectors \(\hat{\dmdmode}\), \(\dmdmode_{m}\) may need to be normalized to unit \(\ell^{2}\) norm.

Combination coefficients \(\dmdcoefficient_{m}\)  can be used to rank the DMD modes by importance.
They can be solved by regression, that is solving~\eqref{eq:DMD-approximation} at the initial snapshot, or even many steps
\begin{equation}
  \label{eq:DMD-coefficients}
  \dmdcoefficients \gets \dmdmodes^{\pinv} \state_{0}, \quad \text{ or } \quad \dmdcoefficients \gets \argmin_{\vec{x} \in \C^{\modeldimension}} \sum_{\timeindex = 0}^{\finaltime}\norm*{\state_{\timeindex} - \sum_{m=1}^{M} \dmdmode_{m} e^{\timeindex \dmdfrequency_{m} \timestepmod} x_{m} }.
\end{equation}
Solving~\eqref{eq:DMD-coefficients} at only the initial state can result in a relative approximation error that is unevenly spread between growing and decaying modes; conversely, using all steps in the regression balances the error, but is numerically expensive.
As a compromise, implementation used below uses five steps equally spaced across all available snapshots.

One advantage of using \gls{DMD} is in additional information that can be used to choose what subset of modes to use for dimension reduction.
\gls{POD} modes are ranked solely by their \(L^{2}\) norms (singular values), therefore most approaches simply choose some number of dominant modes.
A similar effect can be achieved by ranking \gls{DMD} modes by absolute values \(\abs{\dmdcoefficient_{m}}\), which represent contributions of modes to the initial condition.
Alternatively, \gls{DMD} modes can be ordered by \(L^{2}\) norms of time evolution for each \gls{DMD} mode
\begin{equation}
  \label{eq:mean-dmd-coefficient}
  \bar{\dmdcoefficient}_{m}^{2} = \int_{0}^{\finaltime} \abs*{e^{\dmdfrequency_{m} \timeindex} \dmdcoefficient_{m}}^{2} dt = \abs{\dmdcoefficient}^{2} \frac{\exp (2\re \dmdfrequency_{m} \finaltime) - 1}{2 \re \dmdfrequency_{m} \finaltime}, \quad \text{ and } \bar{\dmdcoefficient}_{m}^{2} \coloneqq \abs{\dmdcoefficient}_{m}^{2} \text{ if } \re \dmdfrequency_{m} = 0,
\end{equation}
which give the same weight to modes that grow and those that decay at rates, everything else being the same.
In the applications below, we rank DMD modes in the descending order of \(\bar{\dmdcoefficient}_{m}\), and refer to those with large \(\bar{\dmdcoefficient}_{m}\) as \emph{dominant}.

Alternatively, \gls{DMD} modes can be chosen based on the real or imaginary parts of \(\dmdfrequency_{m}\), e.g., if there is a reason to choose only modes corresponding to a certain frequency band.
We did not pursue this direction further in this paper as we did not suspect that evolutions of  either \gls{L96} or \gls{SWE} were concentrated in a particular frequency band.

Choosing the \(\reducedmodeldimension\) dominant \gls{DMD} modes \(\dmdmode_{m}\), the final step is to form the orthogonal projection matrix \(\projection_{\textrm{DMD}}\).
The truncated \(\dmdmodes\) is not an orthogonal projection because DMD modes, unlike \gls{POD} modes, do not form an orthonormal system.
Additionally, \(\dmdmode_{m}\) are complex-valued, so conjugate pairs of columns \(\dmdmode_{m}, \dmdmode_{m+1} = \dmdmode^{\ast}\) should be replaced by their real-valued cartesian components \(\re \dmdmode, \im \dmdmode\) in \(\dmdmodes\) matrix.
Care should be taken to always include either both element of a pair, or neither.
After these steps that prepare a truncated real-valued basis for a subspace of DMD modes \(\hat{\dmdmodes}\), we compute orthogonal dimension reduction matrix \(\statereduction_{\text{DMD}}\) as left singular vectors of \(\hat{\dmdmodes}\), or state projection matrix \(\projection_{\text{DMD}}\) using QR decomposition of \(\hat{\dmdmodes}\).

\begin{algorithm}[hbt!]
    \caption{Dynamic Mode Decomposition}
    \label{alg:dmd}
    \SetAlgoLined
    $\snapshotmatrix_{1} \gets \begin{bmatrix} \state_0 & \state_1 & \cdots & \state_{T-1} \end{bmatrix}, \ \snapshotmatrix_{2} \gets \begin{bmatrix} \state_1 & \state_2 & \cdots & \state_{\finaltime} \end{bmatrix},$ \tcp{Left/Right Snapshot Matrix}
    $\snapshotmatrix_{1} = \leftsingularvectors \singularvalues \rightsingularvectors^\top$ \tcp{\gls{SVD}}
    $\dmdproblemsize \gets \textrm{user input}$ \tcp{Choice of the DMD problem size \(\dmdproblemsize < \rank \snapshotmatrix \) }
    $\leftsingularvectors_{\dmdproblemsize}, ~\singularvalues_{\dmdproblemsize}, ~\rightsingularvectors_{\dmdproblemsize}$ \tcp{Truncation of SVD matrices}
    $\dmdoperatorreduced \gets \leftsingularvectors_{\dmdproblemsize} \snapshotmatrix_2 \rightsingularvectors_{\dmdproblemsize} \singularvalues_{\dmdproblemsize}$ \tcp{Compressed DMD matrix}
    $\dmdoperatorreduced\hat{\dmdmode}_m = \lambda_m \hat{\dmdmode}_m$ \tcp{Eigendecomposition of DMD matrix}
    $\dmdmode_m \gets \lambda_m^{-1} \snapshotmatrix_2 \rightsingularvectors_{\dmdproblemsize} \singularvalues_{\dmdproblemsize}^{-1} \hat{\dmdmode}_m$ \tcp{Computing and normalizing DMD modes}
    $\dmdcoefficients = \argmin_{\vec{x} \in \C^{\modeldimension}} \sum_{\timeindex = 0}^{\finaltime}\norm*{\state_{\timeindex} - \sum_{m=1}^{M} \dmdmode_{m} e^{\timeindex \dmdfrequency_{m} \timestepmod} x_{m}}$ \tcp{DMD coefficients via \(\ell^{2}\) regression}
    \(\bar{\dmdcoefficient}_{m}^{2} = \abs{\dmdcoefficient}^{2} [\exp (2\re \dmdfrequency_{m} T) - 1]/(2 \re \dmdfrequency_{m} T)\) \tcp{Rank DMD modes in descending order of \(\bar{\dmdcoefficient}\)}
    $\reducedmodeldimension \gets \textrm{user input}$ \tcp{Choice of size of dimension reduction subspace \(\reducedmodeldimension < \dmdproblemsize\) }
    \(\hat{\dmdmodes} \gets \begin{bmatrix}\dmdmodes_{1}&\dmdmodes_{2}&\dots&\dmdmodes_{\reducedmodeldimension}\end{bmatrix}\) \tcp{Truncate DMD modes and convert to real vectors}
    \(\statereduction_{\text{DMD}} \gets \operatorname{SVD}(\hat{\dmdmodes})\) \tcp{Left singular vectors form the dimension reduction matrix}
  \end{algorithm}

\subsection{Assimilation in the Unstable Subspace (AUS)}

The last of the model reduction methods employed here, which is restricted to projecting in physical model space, is based on computational techniques for Lyapunov exponents and Finite Time Lyapunov exponents.
Again, calculated modes are used to select which dimensions are most dynamically significant and should be retained in the reduced basis.
In \gls{AUS}, these modes are determined by employing the discrete QR algorithm \cite{DVV07,DiVV15}.
For the discrete time model $\state_{\timeindex +1} = \modelupdate_{\timeindex }(\state_{\timeindex }) + \modelerror_{\timeindex }$ with $\state_{\timeindex }\in\R^\modeldimension$, let $\mathbf{U}_0\in\R^{\modeldimension\times\reducedmodeldimension }$
$(\reducedmodeldimension\leq \modeldimension)$ denote a random matrix such that $\mathbf{U}_0^{\top} \mathbf{U}_0 = \identity$ and
\begin{align}\label{QR}
  \mathbf{U}_{\timeindex +1} \T_{\timeindex} =& \modelupdate_{\timeindex }'(\state_{\timeindex })\mathbf{U}_{\timeindex } \approx \frac{1}{\epsilon}[\modelupdate_{\timeindex }(\state_{\timeindex } + \epsilon \mathbf{U}_{\timeindex }) - \modelupdate_{\timeindex }(\state_{\timeindex })],\quad \timeindex =0,1,\dots
\end{align}
where $\mathbf{U}_{\timeindex +1}^{\top} \mathbf{U}_{\timeindex +1} = \identity$ and $\T_{\timeindex }$ is upper triangular with positive diagonal elements.
With a finite difference approximation the cost is that of an ensemble of size $p$ plus a reduced QR via modified Gram--Schmidt to re-orthogonalize.
Lyapunov exponents or finite time approximations can be formed and monitored by taking time averages of the natural logarithm of the diagonal elements of the upper triangular $p \times p$ matrices $\T_{\timeindex }$ (see, e.g., \cite{DVV07,DiVV15}).
Time dependent orthogonal projections to decompose state space are given by $\projection_{\timeindex } = \mathbf{U}_{\timeindex } \mathbf{U}_{\timeindex }^{\top}$ and can be employed to form both physical model and data model projections.
While \gls{AUS} traditionally refers to projection/restriction of the physical model onto the neutral and unstable subspace, here we use \gls{AUS} more generally to refer to model-space-based projections using an approximate Lyapunov basis: accordingly, our implementation of \gls{AUS} may include only some unstable modes, or may include the entire neutral-unstable subspace and some stable modes.

\section{Numerical results}\label{results_sec}
\glsreset{L96}
\glsreset{SWE}
To evaluate the performance of \gls{ProjOPPF}, we apply the presented techniques to two commonly used models: \gls{L96} and a configuration of the \gls{SWE} corresponding to a barotropic instability.

The experiments were chosen to evaluate how a particular choice of the order reduction technique, and the dimensions of reduced model and data dimensions, resp. \(\reducedmodeldimension\) and \(\reduceddatadimension\), influence the accuracy of assimilation, as well as protect the particle filter from weight collapse.
In all cases, we use the model-based reduction of the observation space (see \cref{sec:dimension-reduction}) with a full row-rank linear observation operator, which allows for a fair comparison between order reduction techniques.

\subsection{Common experimental setup}\label{sec:exp-setup}

For both \gls{L96} and \gls{SWE} we use a common setup to evaluate the effectiveness of the data assimilation scheme.
We run one simulation of the model that serves as the ``truth'', $\statetruth_{\timeindex}$.
Noisy observations of the ``truth'' are used as inputs into assimilation, and the success of the assimilation is measured by how accurately the state of the estimator matches the state of the ``truth'' simulation.
In all cases we use \gls{ProjOPPF} as the data assimilation scheme, summarized in Algorithm~\ref{alg:ProjOPPFalg}.
The scheme uses a set of $\totalparticles$ particles about the initial condition of the truth $\statetruth$ with added noise from the same Gaussian distribution employed to simulate noise in the physical model.
All particles are initialized with equal weights $(1/L)$ and are propagated forward in time using the chosen model.
\glsreset{ESS} The \gls{ESS} is then calculated~\eqref{eq:ess} and projected resampling is performed with the spread of particles governed by~\eqref{resampling}, where the proportion of resampling variance inside the projection subspace is always taken to be $\alpha = 0.99$,  and total resampling variances $\omega = 10^{-2}$,  for \gls{L96}, and $\omega = 10^{-4}$, for \gls{SWE}
when $\ess<\frac{1}{2}L$.

\glsreset{RMSE}\glsreset{RESAMP}
To evaluate \gls{ProjOPPF} we report on two quantities:
\begin{itemize}
\item \gls{RMSE} between estimate of the state and the true state,
  \begin{equation}
  \rmse(\statetruth,\stateens) \coloneqq \norm*{\statetruth - \stateens}_2/\sqrt{\modeldimension},\label{eq:rmse}
\end{equation}
 where $\statetruth$ denotes the truth and $\stateens$ denotes the particle ensemble mean, and
\item \gls{RESAMP}, which measures the proportion of observation times in which the particle population needed to be resampled.
\end{itemize}
The lower each of the quantities are, the better the assimilation scheme is performing.
We also report in some experiments on the projected \gls{RMSE},
\begin{equation}
\reducedrmse(\statetruth,\stateens) \coloneqq \rmse(\statereduction\statereduction^\top\statetruth,\statereduction\reducedstate^{\textrm{ens}}) = \norm*{\statereduction\statereduction^\top\statetruth - \statereduction\reducedstate^{\textrm{ens}}}_{2}/\sqrt{\reducedmodeldimension}.\label{eq:projected-RMSE}
\end{equation}
which measures the error in the subspace of the reduced model.
A low projected \gls{RMSE} together with a significantly higher ``full'' \gls{RMSE}~\eqref{eq:rmse} is an indication that the assimilation is being effective when restricted to the subspace of the reduced model, but the projected model does not sufficiently resolve the full model.
In several of the figures, we will compare with the results of the optimal proposal particle filter with no model or data reduction and we will use (NON) to denote these results.
All numerical results are obtained by averaging over 10 randomized trials.

\glsreset{L96}
\subsection{Lorenz '96 Equations}

\subsubsection{Model and parameters}\label{sec:l96-model}

We first consider \gls{ProjOPPF} applied to the extensively used medium-dimensional dynamical system \gls{L96}. This model, developed by Edward Lorenz \cite{Lorenz96}, represents a nonlinear chaotic system that captures some multiscale features of the global horizontal circulation of the atmosphere.
\gls{L96} has become one of the most commonly used test problems in data assimilation since its introduction.

The model is presented as a system of \glspl{ODE} in \(\state = (u_{i})_{i=1}^\modeldimension\) of an arbitrary dimension \(\modeldimension\),
\begin{align} \label{eq:l96}
    \frac{d u_{i}}{d\timeindex} =& \left(u_{i+1}-u_{i-2}\right)u_{i-1} -u_{i} +F, \quad i=1,\dots,\modeldimension,
\end{align}
where the value of a constant (typically positive) forcing term $F$ determines qualitatively whether the evolution will be regular or chaotic. In its original form, the model was introduced as a variable-order system of \glspl{ODE}, but it can be interpreted as a 2nd order finite-difference discretization of a viscous Burgers-type equation with periodic boundary conditions (see \cite{Blender2013} for derivation):
\begin{equation}
  \label{eq:l96-continuous}
  \partial_\timeindex u = -u\partial_x u - \frac{1}{3}(\partial_x u) ^{2} -\frac{1}{6} u \partial_{xx} u  - u + F.
\end{equation}
In this form, the model has a nonlinear convective term corresponding to $u\partial_x u$, a diffusive terms, and a dissipative term $(-u)$. Note, however, that the discretization does introduce nonlinear effects not present in the full model~\cite{vanKekem2018a}, so that the ODE--\gls{PDE} correspondence is not straightforward.

\begin{figure}[ht]
    \centering
    \begin{subfigure}[t]{0.33\linewidth}\centering
\includegraphics[width=\textwidth]{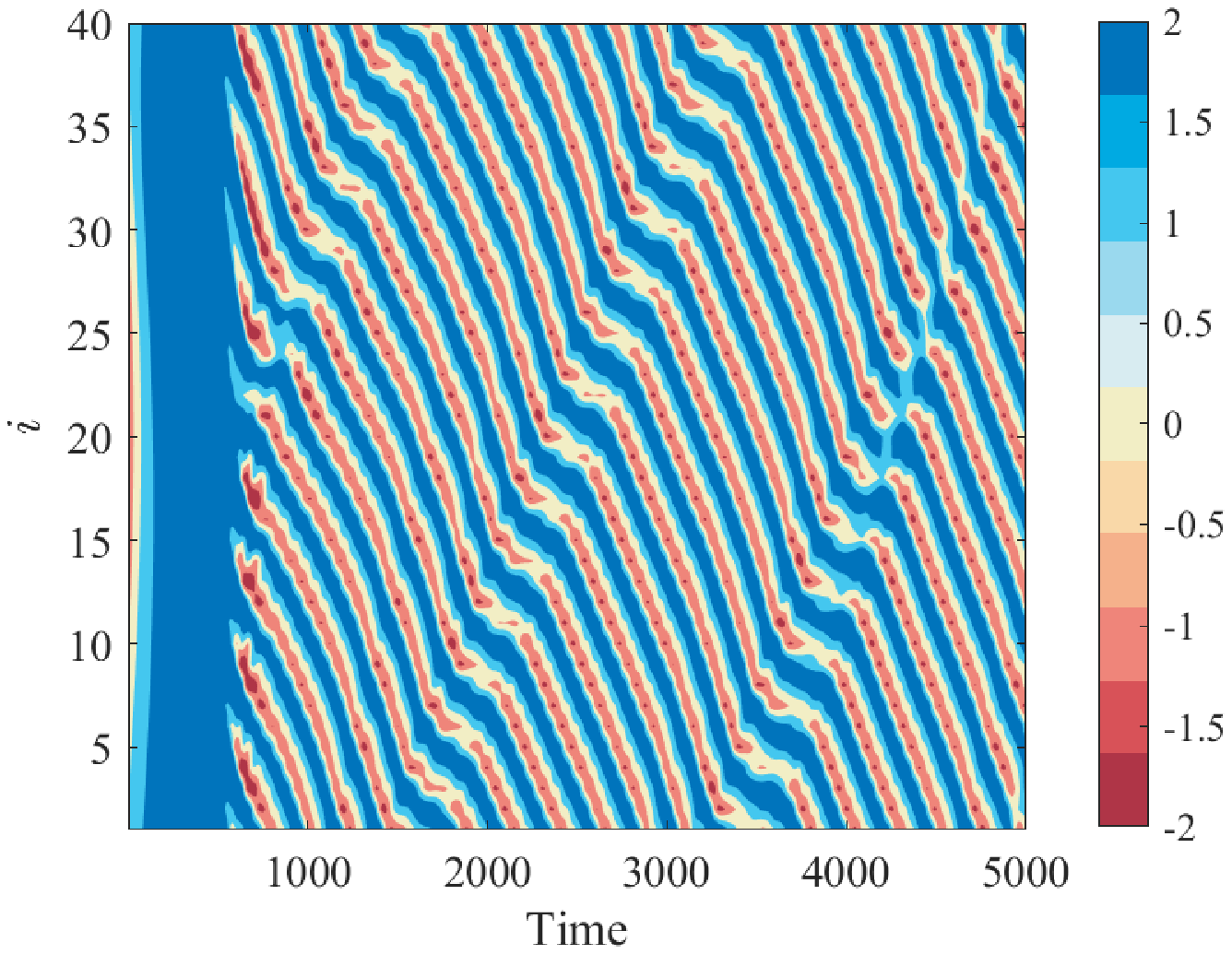}
  \caption{\(F=3\)}
\end{subfigure}
    \begin{subfigure}[t]{0.33\linewidth}\centering
\includegraphics[width=\textwidth]{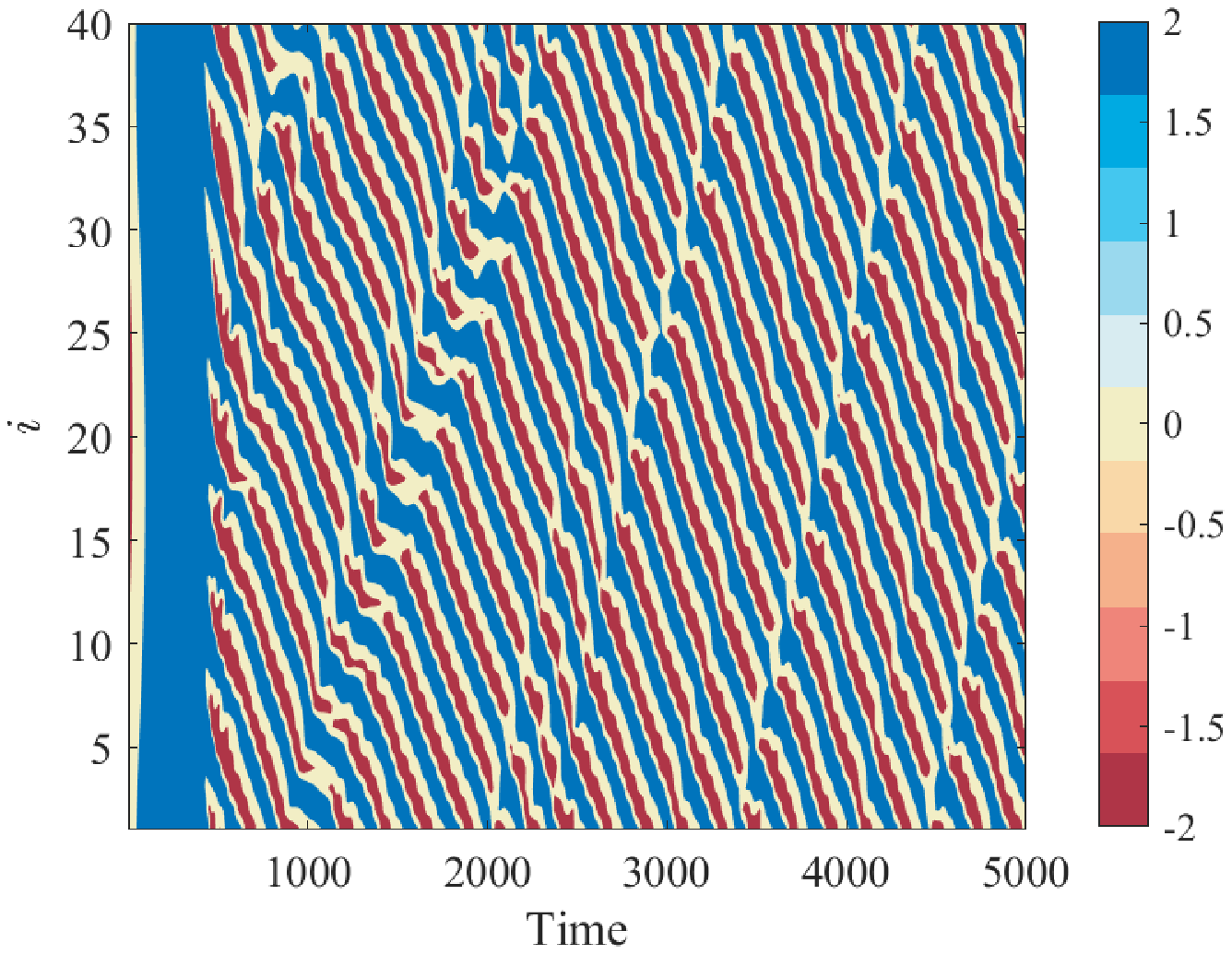}
  \caption{\(F=4\)}
\end{subfigure}

\begin{subfigure}[t]{0.33\linewidth}\centering
\includegraphics[width=\textwidth]{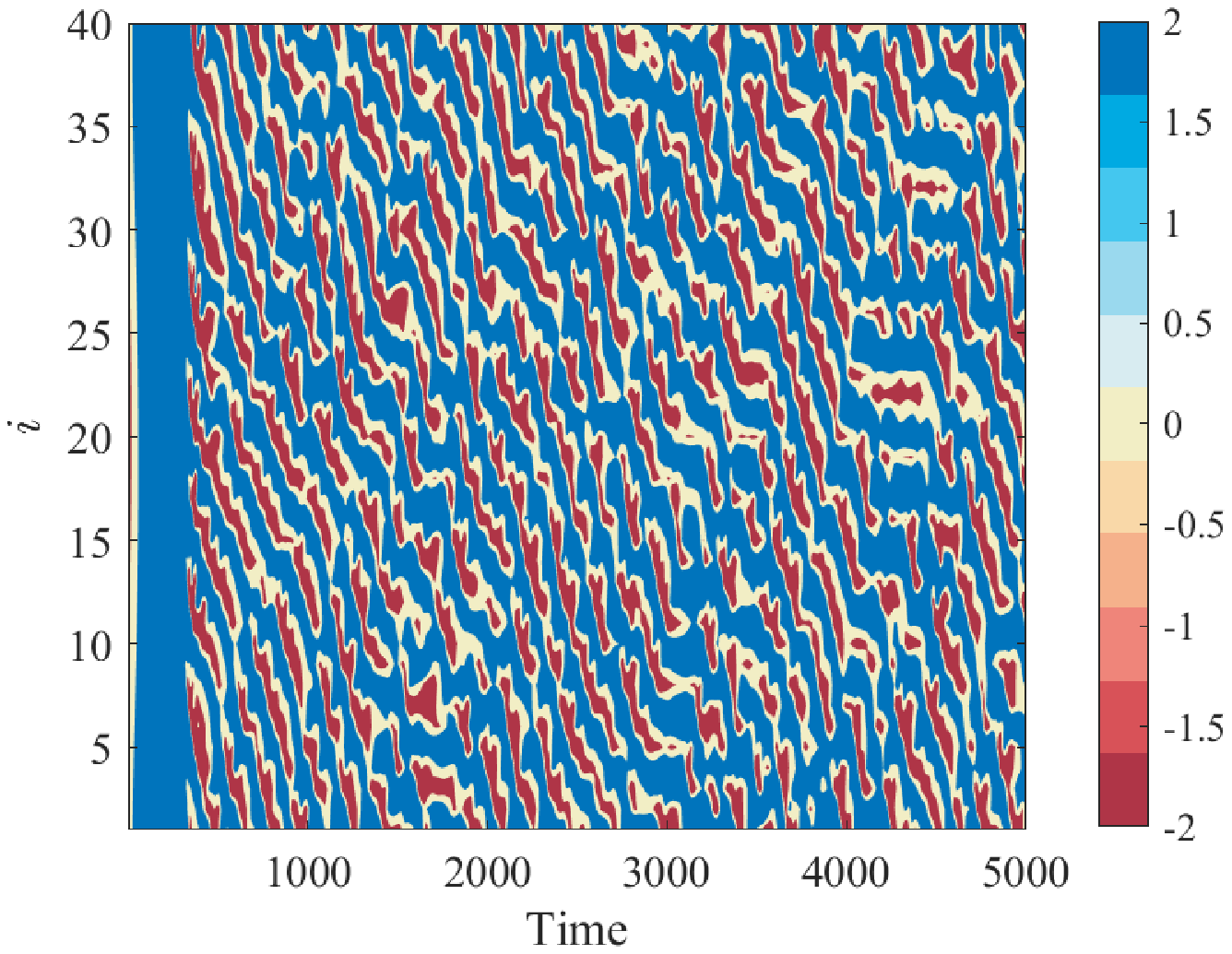}
\caption{\(F=6\)}
\end{subfigure}
\begin{subfigure}[t]{0.33\linewidth}\centering
\includegraphics[width=\textwidth]{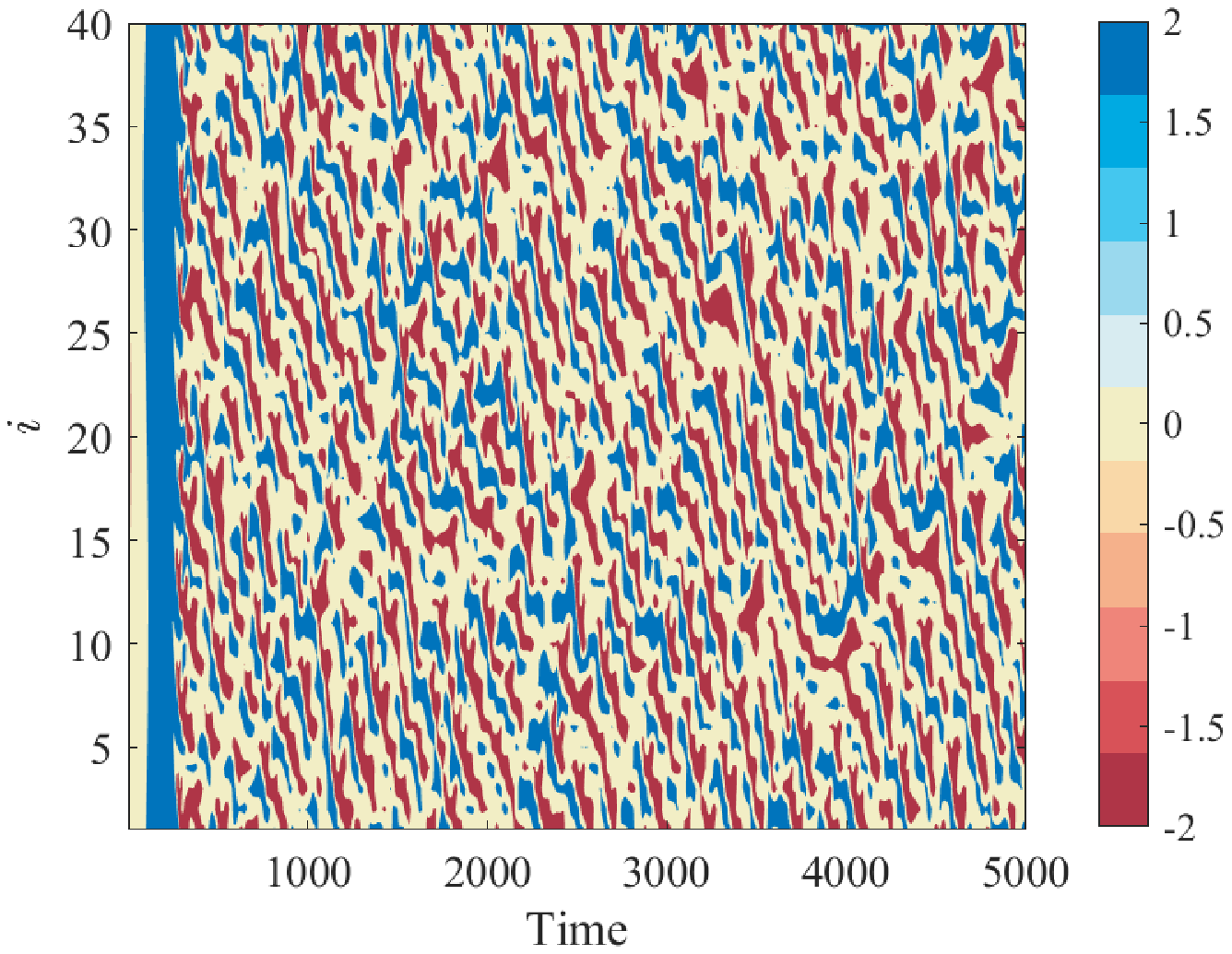}
\caption{\(F=8\)}
\end{subfigure}
\caption{Solutions \(u_{i}(\timeindex)\), $i=1,\dots,\modeldimension\equiv 40$ \; of \eqref{eq:l96} demonstrating regular and chaotic spatiotemporal patterns in \gls{L96} model for various values of forcing \(F\). Initial condition is a cosine bump in all cases. %
}
  \label{fig:l96-space-time}
\end{figure}

To produce the true evolution \(\statetruth_\timeindex\), the \gls{L96} model~\eqref{eq:l96}  is evolved in time using the Dormand--Prince
pair (MATLAB's \texttt{ode45}), with solution resampled at multiples of the
fixed time step $\timestepmod=0.01$.
Figure \ref{fig:l96-space-time} illustrates the space-time behavior of typical solutions of \gls{L96} ranging from $F=3$ (regular structure) to $F=8$ (chaotic structure). For $\modeldimension=40$ used in most calculations here, the onset of chaos is between $F=3$ and $F=4$~\cite{vanKekem2018a}; the similar behavior appears to hold for $\modeldimension=400$ as well.
In the numerical experiments, the initial condition is a random vector to reduce the burn-in time, and provide some variability between trials.

In the setup used here, all model variables are observed so that $\dataoperator = \identity$.
The noise is uncorrelated among any two variables of the system (within and between state and observation vectors).
Moreover, it is uniform for all states and for all observation variables, which means that the correlation matrices \(\modelerrorcovariance, \dataerrorcovariance\) are scaled identity matrices.
We compare several levels of such noise in the physical model, yielding
$\modelerrorcovariance=\alpha\cdot \identity$ with scalar variances
$\alpha = 0.1, 1.0$.
The observation error covariance is always fixed to \(\dataerrorcovariance = 0.01\,\identity\) resulting in the standard deviation of observation error of \(0.1\) which is included for comparison in the figures when reporting on \(\rmse\).

The number of particles is fixed with $\totalparticles=20$
The observations are computed every 5 steps, yielding the effective time step of assimilation \(\timestepobs = 0.05\).
For the purposes of assimilation, the 4th order Dormand--Prince integrator is used but with a fixed step size \(\timestepmod =0.01\).
The assimilation is performed over $10,000$ observation times, after 1000 time steps have elapsed.
The average \gls{RMSE} over time is calculated based on the \gls{RMSE} over the last $5,000$ observation times (see \Cref{fig:l96-RMSE-in-time} where the last $5,000$ observation times correspond to the second half of the assimilation window), to more accurately represent the asymptotic value of the \gls{RMSE}.
The resampling percentage \gls{RESAMP} is computed as the proportion of \emph{all} $10,000$ observation times in which resampling was performed.

For the classical values of the \gls{L96} system where $\modeldimension = 40$ and $F=8$, the system is chaotic with $13$ positive and $1$ neutral Lyapunov exponent.
The percentage of positive Lyapunov exponents approximately scales with the dimension $\modeldimension$ and there are generally a smaller percentage of positive Lyapunov exponents for smaller values of $F$.
For $\modeldimension=40$,  $\modeldimension=400$ and $F\in[3,8]$ we approximate the Lyapunov dimension by Kaplan--Yorke formula \(D_{L}\) defined for ordered Lyapunov exponents $\lambda_1 \geq \lambda_2 \geq \cdots \geq \lambda_M$ given by
\begin{equation}
D_L = k + \frac{\lambda_1 + \lambda_2 + \cdots + \lambda_k}{|\lambda_{k+1}|} \label{eq:kaplan-yorke}
\end{equation}

where $k$ is the maximum value of $i$ such that $\lambda_1 + \lambda_2 + \cdots + \lambda_i>0$.
Approximations are computed over a time interval of length $10,000$ with the integer part tabulated in \Cref{tab:l96-lyap-dim} and indicates that for smaller values \(F\) the dynamics can be thought of as inherently low dimensional as \(D_{L}/\modeldimension \ll 1\), but for larger values of \(F\) this is not the case.

\begin{table}[H]
  \centering
  \rowcolors{1}{}{white!90!black}\small
  \begin{tabular}{ >{\bfseries\normalsize}l | c | c | c | c |}
\(F\) & 3 & 4 & 6 & 8\\\hline
\(D_{L}\) (\(\modeldimension=40\)) &  1 & 3 & 22 & 28\\
    \(D_{L}\) (\(\modeldimension=400\)) & 5 & 12 & 224 & 270\\
  \end{tabular}
  \caption{Lyapunov dimensions for the \(\gls{L96}\) model for various forcing values and model dimensions \(\modeldimension\), estimated by the Kaplan--Yorke formula~\eqref{eq:kaplan-yorke}.}
    \label{tab:l96-lyap-dim}
\end{table}

\begin{figure}[H]
    \centering
\begin{subfigure}[t]{0.32\linewidth}\centering
\includegraphics[width=\textwidth]{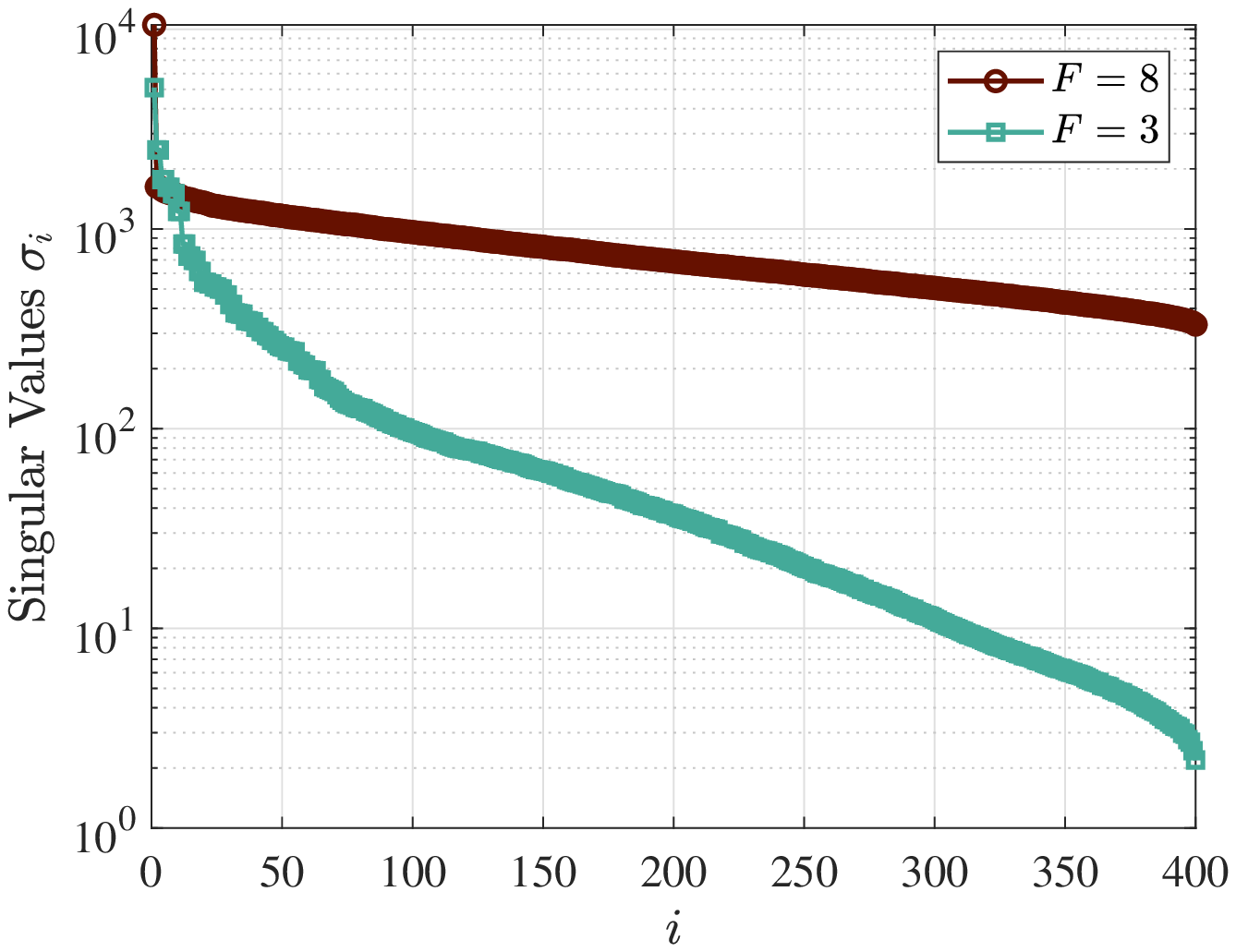}
  \caption{Singular values associated with \gls{POD} modes.}
  \label{fig:l96-svd}
\end{subfigure}
\begin{subfigure}[t]{0.32\linewidth}\centering
\includegraphics[width=\textwidth]{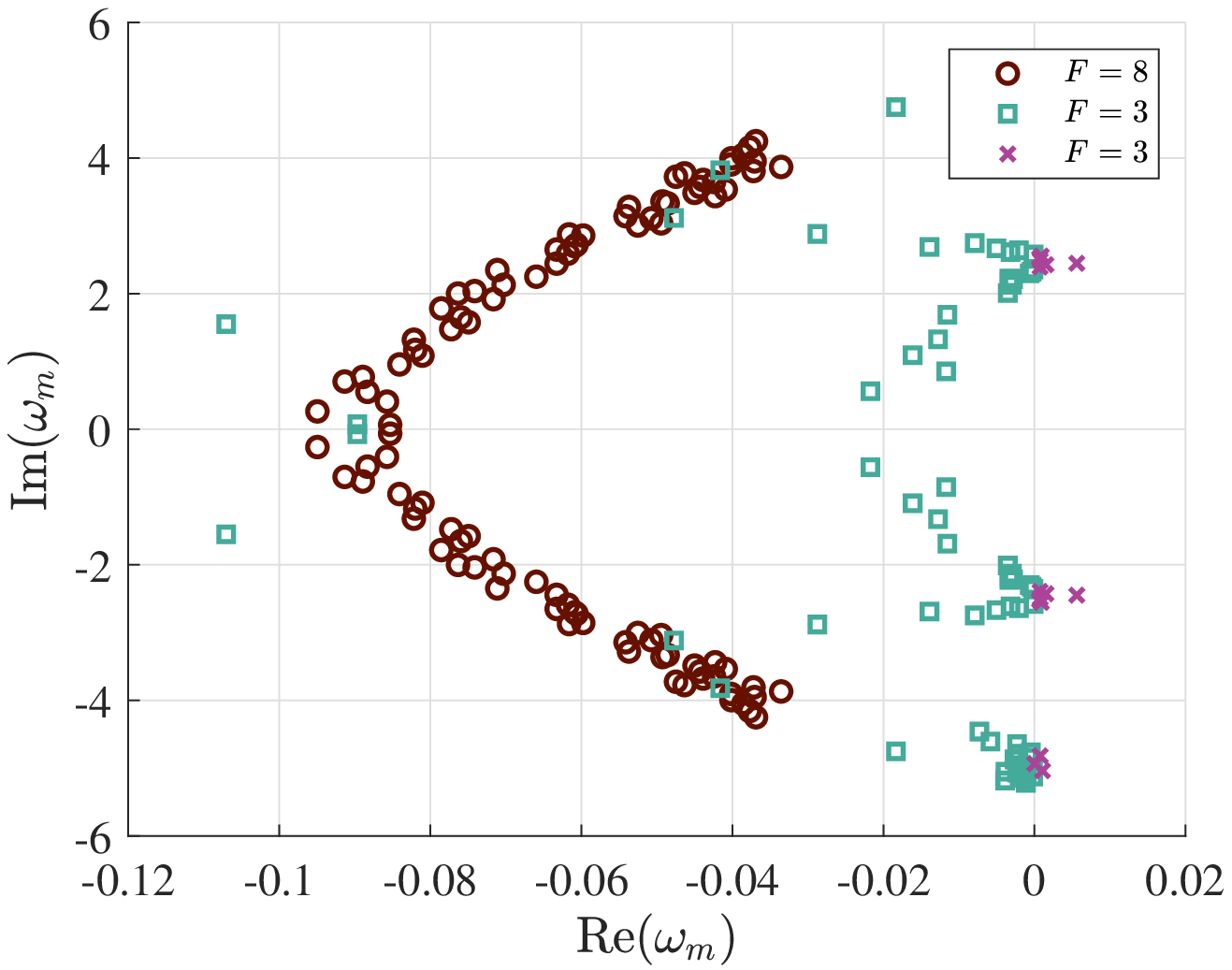}

\caption{\gls{DMD} eigenvalues in the complex plane.
\(\re (\dmdfrequency) < 0\) signify decaying modes, \(\im (\dmdfrequency)\) is proportional to the frequency of oscillations.}
\label{fig:l96-dmd}
\end{subfigure}
\begin{subfigure}[t]{0.32\linewidth}\centering
\includegraphics[width=\textwidth]{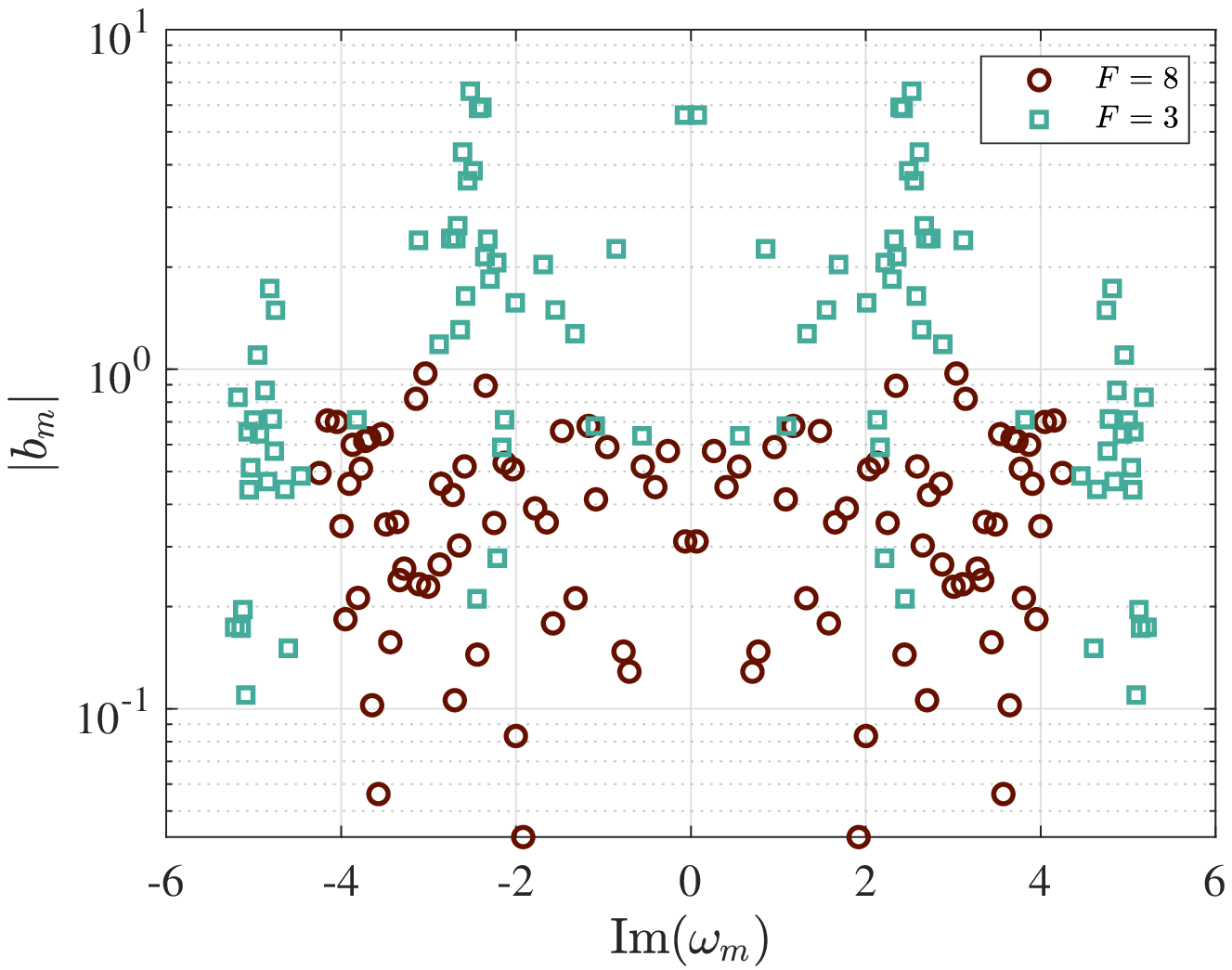}
  \caption{\gls{DMD} coefficients~\eqref{eq:DMD-coefficients}.}
\end{subfigure}
\caption{Singular values and DMD eigenvalues spectra for \gls{L96} model with \(\modeldimension=400\).
Lower values of \(F\), associated with more-regular behavior, demonstrate a faster decay of the singular values, and clustering of \gls{DMD} frequencies.
Higher values of \(F\), associated with disordered behavior, have a slower decay of the singular values and more evenly distributed spectrum.
}
  \label{fig:l96-singular-and-eigenvalues}
\end{figure}

To determine suitable \gls{POD} and \gls{DMD} modes for \gls{L96}, we  simulate the model with a random initial condition (different from those used to compute the ``truth''), over the entire time interval over which the assimilation is performed at the model time step \(\timestepmod = 0.01\).
Since the model evolution is largely in a steady state, both \gls{POD} and \gls{DMD} modes from such simulation can be used for the assimilation process.
This setting can be justified in cases where the model evolution is largely in a steady state, that is, not undergoing a regime change.
A more realistic setup could employ decompositions of assimilated data over some moving interval, which is something we will explore in future work.
For \gls{L96} the \gls{POD} and \gls{DMD} projections are computed separately for each of the 10 trial using snapshots obtained from a random initial condition.

Order reduction techniques are effective when the dynamics can be reproduced by a relatively small number of modes.
For POD, this is indicated by jumps (``gaps'') in the singular value spectrum.
Additionally, sharp changes in slope between segments of the singular values plot indicate the presence of mechanisms evolving at different scales.
\Cref{fig:l96-singular-and-eigenvalues} demonstrates that the regular $F=3$ evolution results in an ``elbow'' around subspace of dimension 70, while for \(F=8\) there is no discernible boundary between system scales (elbows in \Cref{fig:l96-svd}).

The duration of the time window over which \gls{POD} and \gls{DMD} are computed depends on the type of dynamics.
For regular steady-state behavior, the time window should be large enough for the trajectory to trace out the orbit at least once.
For irregular steady-state behavior, e.g., chaotic attractors, the time window should be long enough so that the points in the trajectory sample the attractor well.
All parameters $F$ considered for \gls{L96} lead to steady state behavior, although we leave the investigation of the effects of window duration for future work.

For transient behavior, e.g., trajectories shadowing heteroclinic orbits between two invariant sets, the duration of the window has to be chosen carefully.
In those cases, a common strategy is to employ a sliding-window computation of POD or DMD~\cite{Kutz2016a,Brunton2016}, leading to time-varying reduction/projection matrices.
The choice of the duration of the window leads to delicate effects, further explored in~\cite{Page2019}.
Theoretical connection between sliding window DMD and the stability theory of time-varying dynamical systems has been developed in~\cite{Macesic2018,Macesic2020,Crnjaric-Zic2020}.
A follow-up publication will investigate how time-varying projection matrices perform in particle filters, and we plan to explore these matters in more detail.

For DMD, regular steady-state dynamics are indicated by DMD frequencies \(\dmdfrequency_{m}\), see \eqref{eq:DMD-frequencies}, concentrated along the vertical axis, and by a ``peaked'' graph of combination coefficients \(\dmdcoefficient\) against the frequency.
\Cref{fig:l96-singular-and-eigenvalues} shows the singular values and DMD spectrum for the range of forcing parameters \(F\) distinguishing between DMD spectrum with positive and negative real parts for \(F=3\) in \Cref{fig:l96-dmd}.
It is evident that we can expect a better performance of order reduction techniques for lower values of \(F\), which are associated with regular behavior of the spatiotemporal evolution.

To evaluate the performance of \gls{ProjOPPF} we conduct three experiments with the intent of comparing how three reduction schemes performed over a range of parameters. Details of experiments are given in \Cref{tab:l96-exp}.
\begin{table}[H]
  \centering
  \rowcolors{1}{}{white!90!black}\small
  \begin{tabular}{ >{\bfseries\normalsize}c | c | c | c | c | c | c | c   }
    {Exp.} & {\(F\)} & {\(\modeldimension=\datadimension \)}  & {\(\modelerrorcovariance\) }   & Reduction & {\(\reducedmodeldimension\)}  & {\(\reduceddatadimension\)} \\ \hline
    1 & 3,8 & 40  & \(0.1\,\identity, 1.0\,\identity\) & AUS, POD, DMD & 5 -- 40 & 5  \\
    2 & 3,8 & 400 & \(0.1\,\identity\) & POD & 100 & 1 -- 100 \\
    3 & 3,4,6,8 & 400  & \(0.1\,\identity, 1.0\,\identity\) & POD, DMD & 10--400 &  5 \\
    \end{tabular}
    \caption{Parametrization of experiments used to test \gls{ProjOPPF} for \gls{L96}. In all cases, we set observation noise covariance to $\dataerrorcovariance=0.01\,\identity$, number of particles $\totalparticles=20$, and observe all model variables $\dataoperator = \identity$. }\label{tab:l96-exp}
  \end{table}
\subsubsection{Experiment 1 ($F=3, 8$, $\modeldimension=40$)}\label{sec:L96-ex1}

\begin{figure}[h!]
  \centering
\begin{subfigure}[t]{0.495\linewidth}\centering
\includegraphics[width=\textwidth,
  trim={0 0 0 28pt},clip
  ]{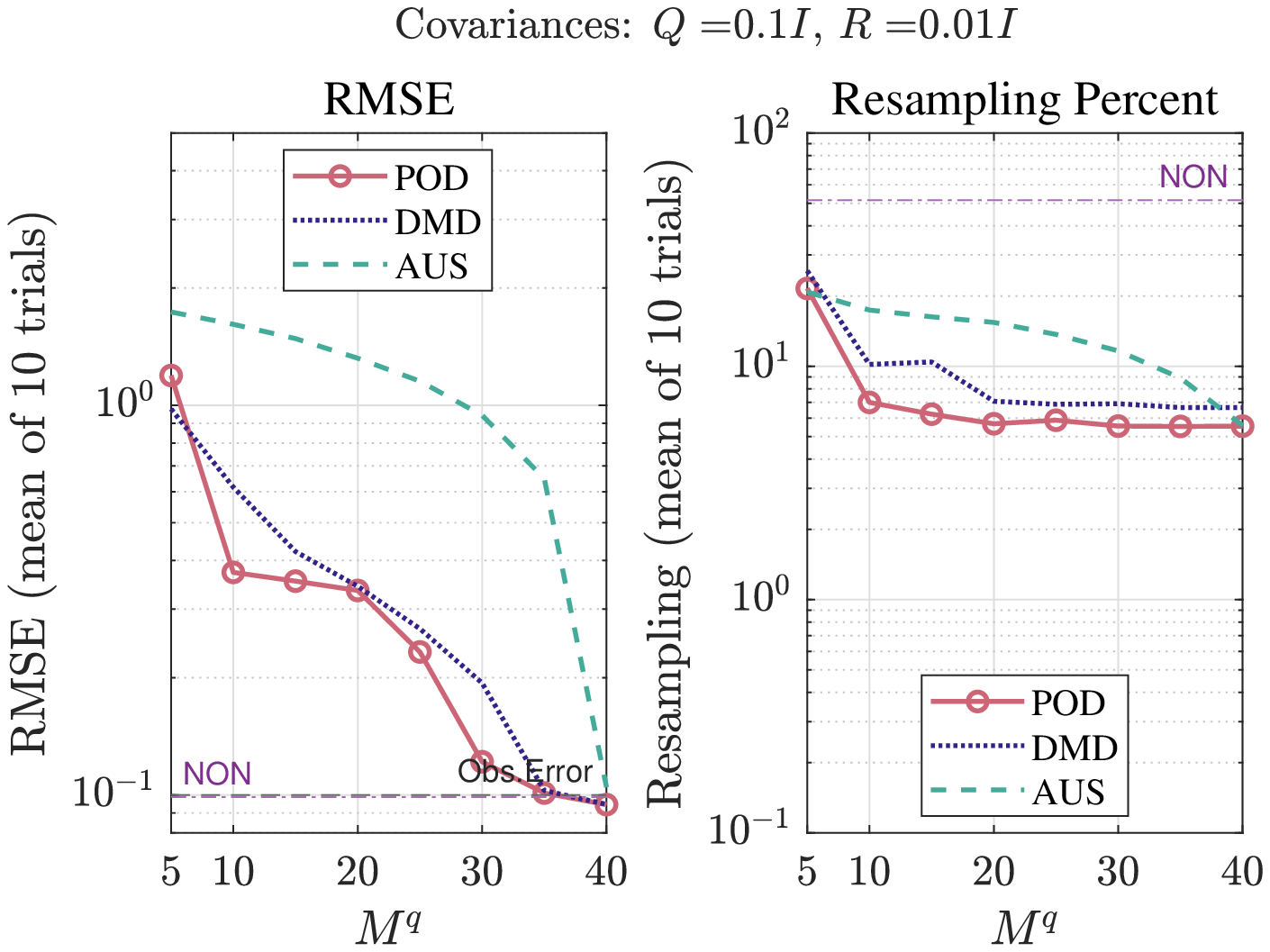}
  \caption{\(F=3\), \(\modelerrorcovariance = 0.1\,\identity\)}
  \end{subfigure}
  \begin{subfigure}[t]{0.495\linewidth}\centering \includegraphics[width=\textwidth,
  trim={0 0 0 28pt},clip
  ]{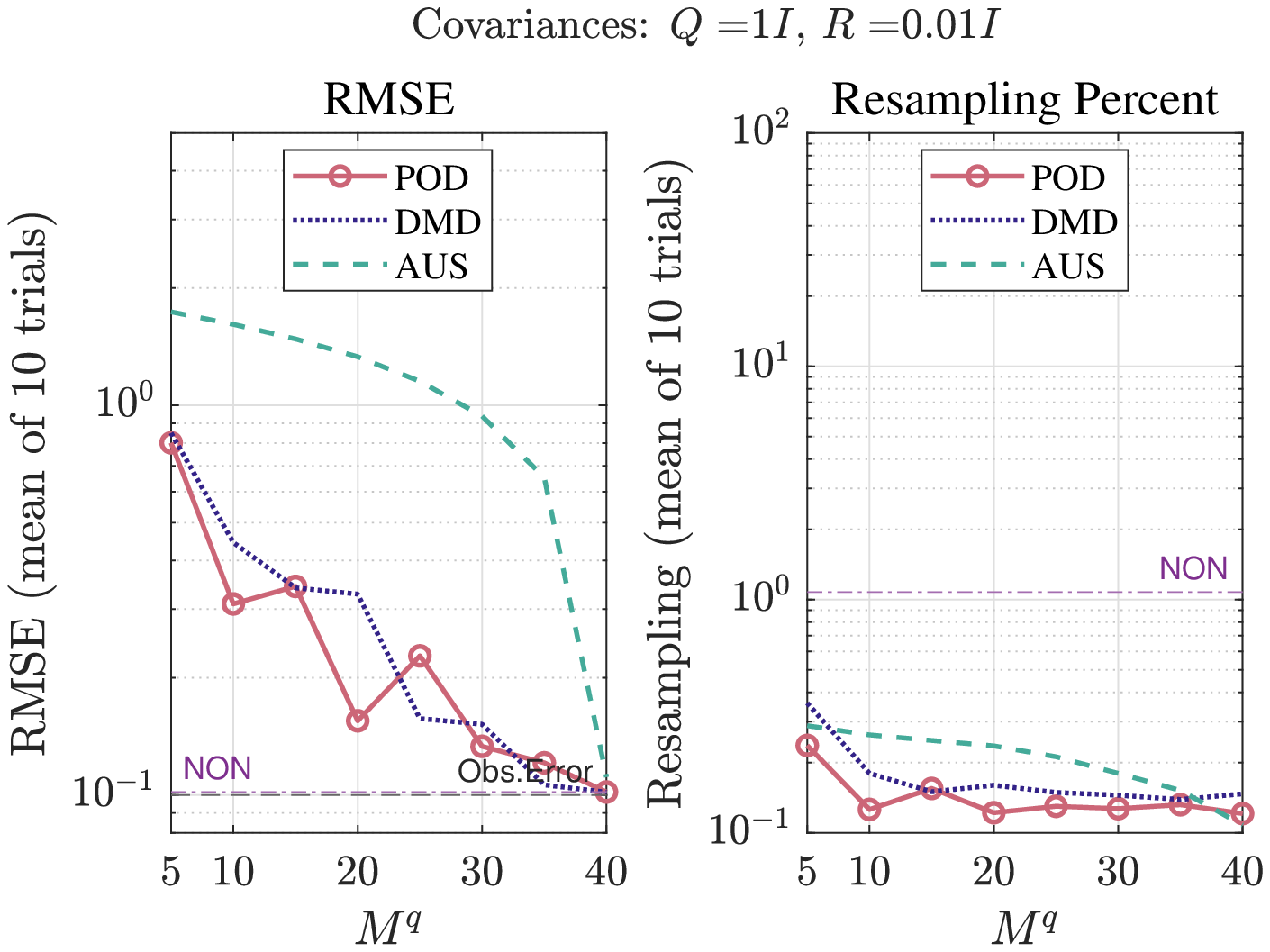}
  \caption{\(F=3\), \(\modelerrorcovariance = 1.0\,\identity\)}
  \end{subfigure}

\begin{subfigure}[t]{0.495\linewidth}\centering
\includegraphics[width=\textwidth,
  trim={0 0 0 28pt},clip
  ]{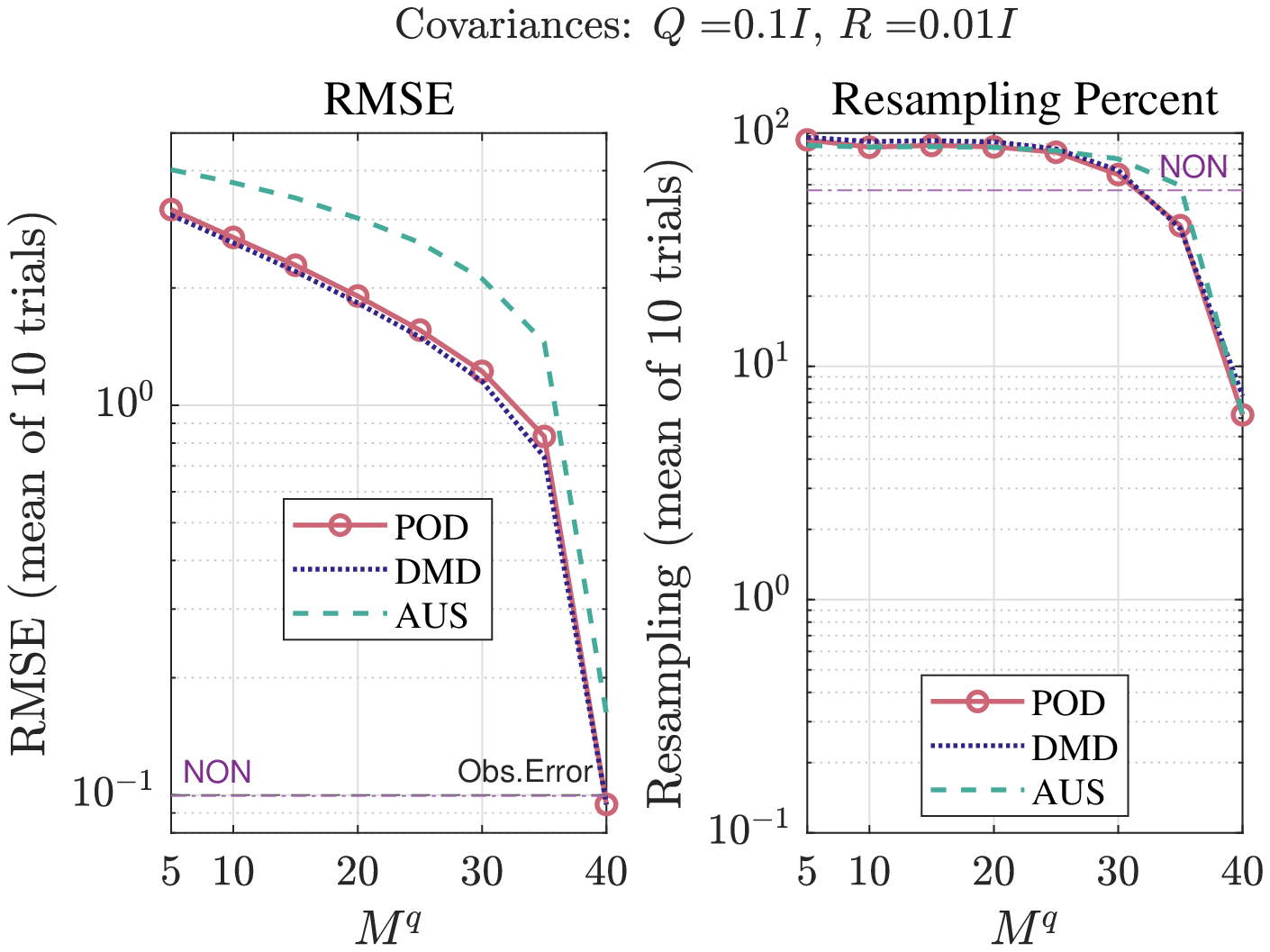}
  \caption{\(F=8\), \(\modelerrorcovariance = 0.1\,\identity\)}
  \end{subfigure}
  \begin{subfigure}[t]{0.495\linewidth}\centering \includegraphics[width=\textwidth,
  trim={0 0 0 28pt},clip
  ]{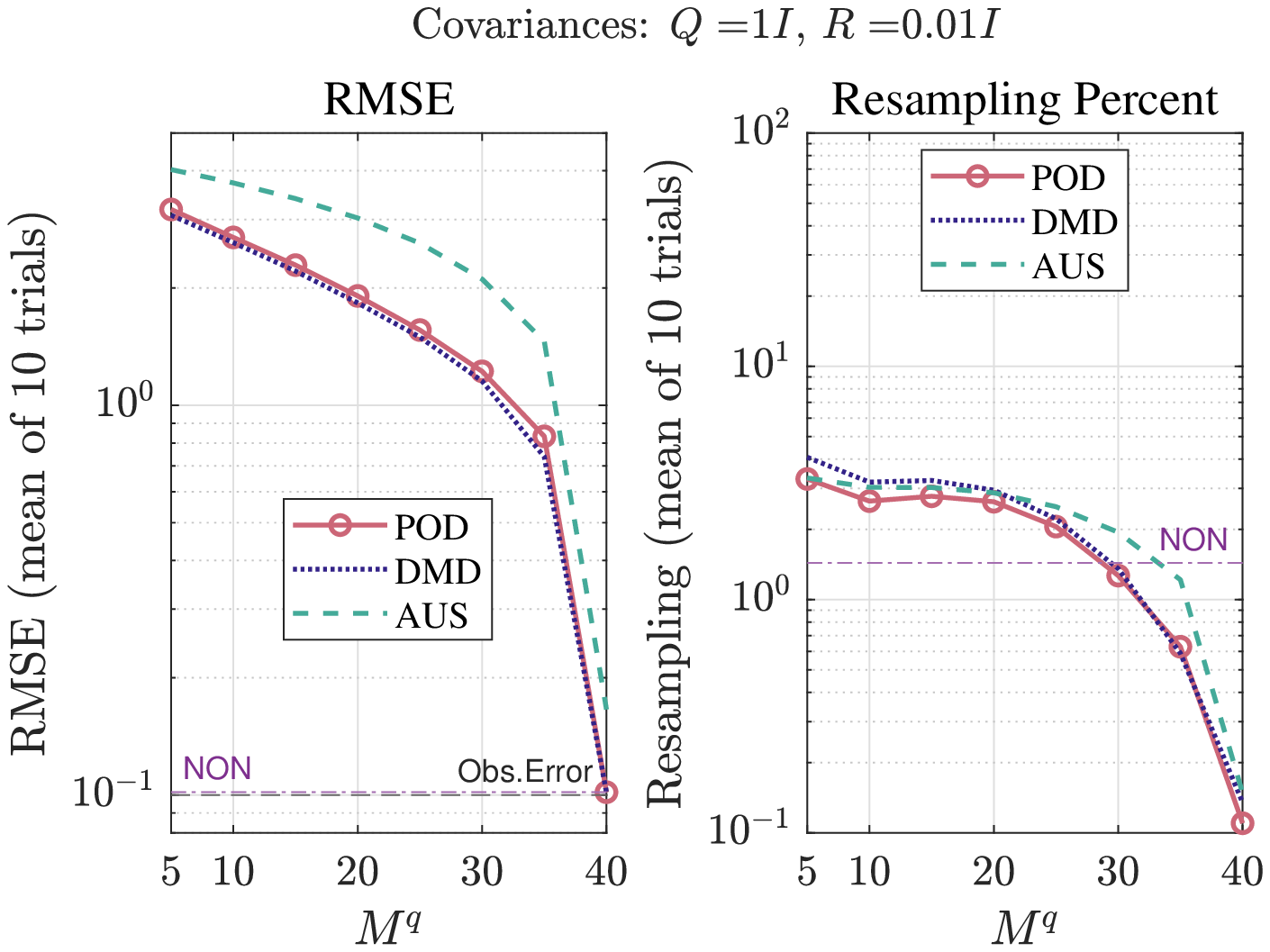}
  \caption{\(F=8\), \(\modelerrorcovariance = 1.0\,\identity\)}
  \end{subfigure}
  \caption{Experiment 1 for the \gls{L96} model with $\modeldimension=40$, \gls{ProjOPPF} using  \gls{DMD}, \gls{POD}, and \gls{AUS} for physical and data models with fixed data projection dimension $\reduceddatadimension=5$ and varying the rank of the physical model projection dimension ($\reducedmodeldimension = 5,10,\dots,\modeldimension = 40$).
Standard deviation of the observation noise (\(\dataerrorcovariance = 0.01\,\identity\)) is given by the dash-dot line in \gls{RMSE} panels. NON corresponds to \gls{ProjOPPF} with the identity projection $(\reducedmodeldimension=\modeldimension\quad {\rm and}\quad \reduceddatadimension = \datadimension)$ or equivalently no reduction of both the model and the data.
See \cref{sec:L96-ex1} for further details.}
  \label{l96_pod_dmd_aus1}
  \end{figure}

Using \gls{AUS} for reducing the order of the \emph{data} model  was investigated in \cite{MVV20, carrassi2020}; here we investigate the efficacy when \gls{AUS} is used for the reduction of both the physical and data models.
Additionally, we compare \gls{AUS} with order reductions derived using \gls{POD} and \gls{DMD}.
\Cref{l96_pod_dmd_aus1} shows the mean \gls{RMSE} and \gls{RESAMP} trends with increasing model projection rank for the three projection methods.
For all projection types the \gls{RMSE} does not approach the observation noise until the rank of the projections are at least $35$ although the \gls{RMSE}s for \gls{POD} and \gls{DMD} are slightly lower than for \gls{AUS}. With $\modelerrorcovariance=1.0\,\identity$ the \gls{RESAMP} are much smaller than for $\modelerrorcovariance = 0.1\,\identity$
since the optimal proposal particle filter with larger model error covariance leads to greater particle diversity and hence smaller \gls{RESAMP}.
We focus on dimension reduction using \gls{POD} and \gls{DMD} for the remaining experiments.

\subsubsection{Experiment 2 ($F=3, 8$, $\modeldimension=400$)}\label{sec:L96-ex2}

\begin{figure}[ht!]
  \centering
\begin{subfigure}[t]{0.495\linewidth}\centering
  \includegraphics[trim={0 0 0 28pt},clip,width=\textwidth]{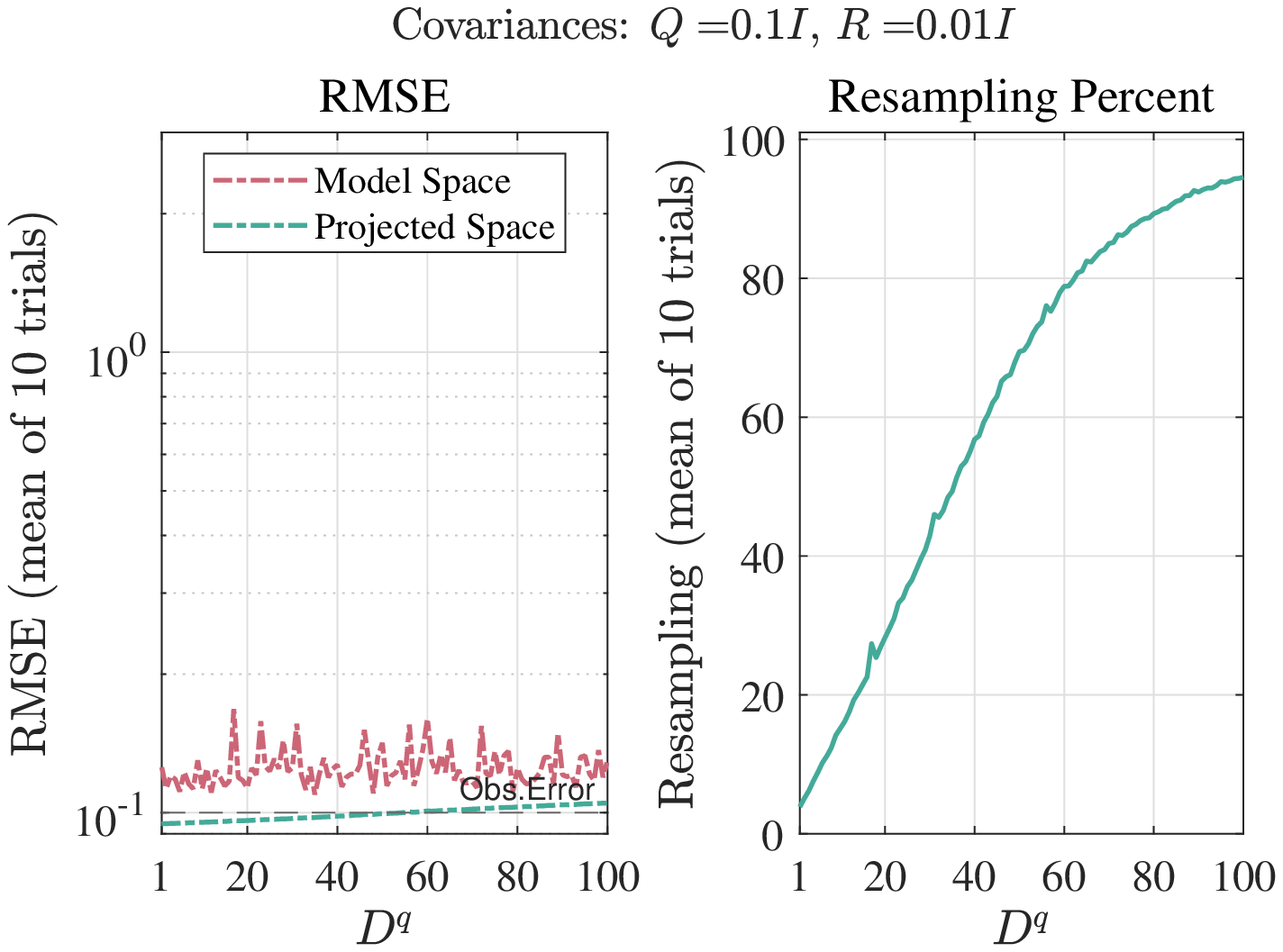}
  \caption{\(F=3\)}\end{subfigure}
\begin{subfigure}[t]{0.495\linewidth}\centering
  \includegraphics[trim={0 0 0 28pt},clip,width=\textwidth]{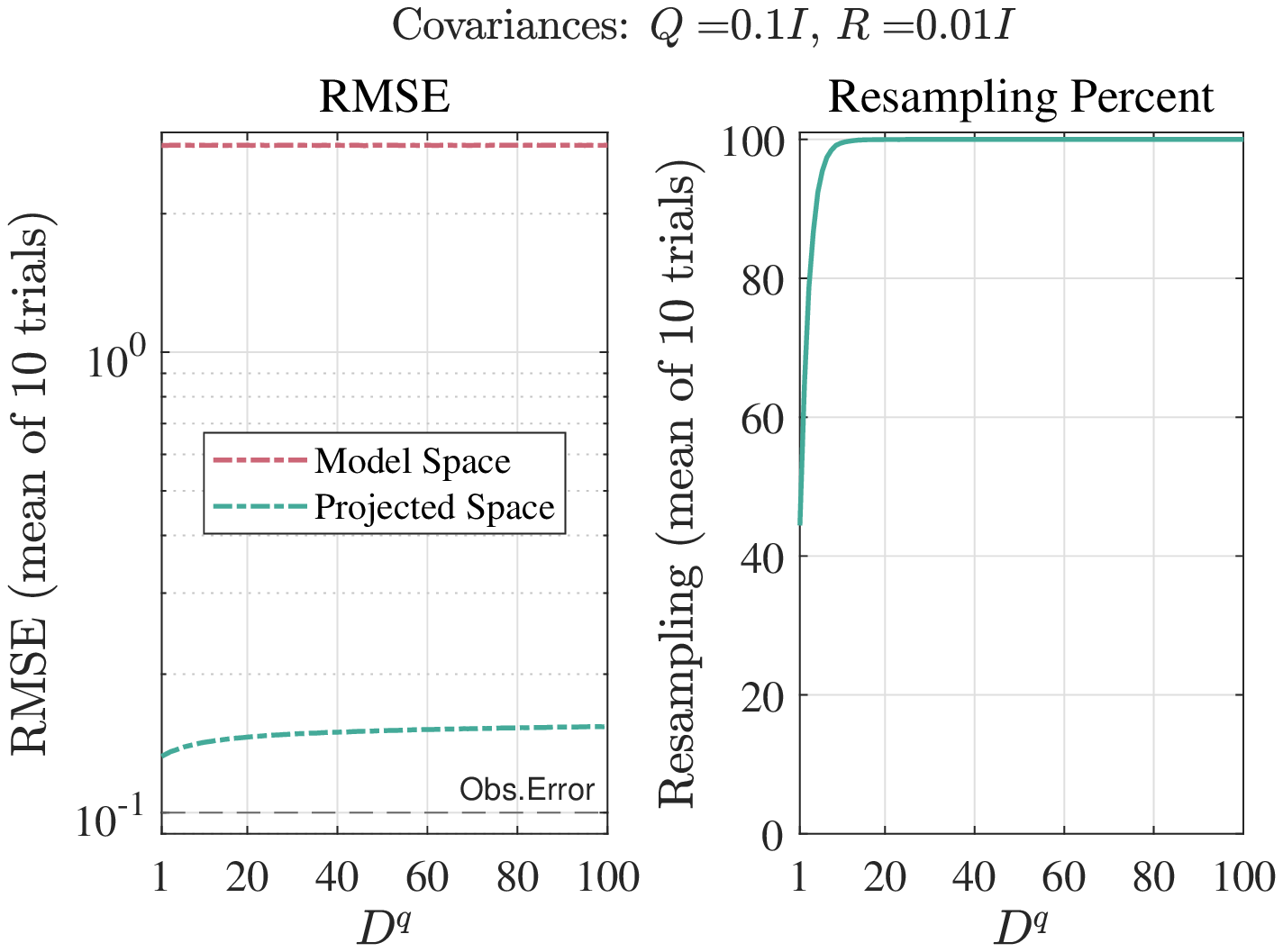}
  \caption{\(F=8\)}\end{subfigure}
\caption{Influence of the reduction of order of the observation space for assimilation of \gls{L96} model with \(F=3\) and \(F=8\) for the experiment 2 (see \cref{tab:l96-exp}). Assimilation was performed using \gls{ProjOPPF} with \gls{POD} projection of physical and data models \(\modeldimension=\datadimension=400\) with fixed physical projection dimension $\reducedmodeldimension=100$ and varying the data model projection with rank $\reduceddatadimension=1,2,\dots,\reducedmodeldimension=100$. See \cref{sec:L96-ex2} for further details. }
  \label{L96_data}
\end{figure}

We set the model dimension to \(\modeldimension = 400\) and fix \(\reducedmodeldimension = 100\), then vary \(\reduceddatadimension\) between \(1\) and \(100\). The results are shown in \Cref{L96_data}.
The \gls{RMSE} and projected \gls{RMSE} are both relatively constant over the range of data dimensions.
The projected \gls{RMSE} is at the level of standard deviation of observation noise, which indicates that the assimilation is effective in the states \(\state\) belonging to the subspace of the reduced model.
The higher value of the ``full'' \gls{RMSE} indicates that this is not sufficient to constrain the full state, meaning that the projected model does not sufficiently resolve the full model.

The proportion of time steps in which particle resampling was performed, \gls{RESAMP}, increases steadily with the dimension of the projection of the data space.
This indicates that the primary effect of \(\reduceddatadimension\) is to mitigate the weight collapse of the particle filter, without significantly affecting the accuracy of the assimilation.
This is a confirmation of the effectiveness of using data-driven order reduction techniques, \gls{POD} and \gls{DMD}, instead of the model-based \gls{AUS}, as detailed in~\cite{MVV20}.

\subsubsection{Experiment 3 ($F=3,4,6,8$, $\modeldimension=400$)}\label{sec:L96-ex3}

Next we consider \gls{L96} system with $\modeldimension = 400$ and varying values of the forcing term $F$ and focus on time-independent \gls{POD} and \gls{DMD} based projections. We revert to $\reduceddatadimension=5$. As before, we average the results over 10 trials with randomized initial conditions and noise realization.

\begin{figure}[h!]
  \centering
  \begin{subfigure}[t]{0.495\linewidth}\centering \includegraphics[width=\textwidth]{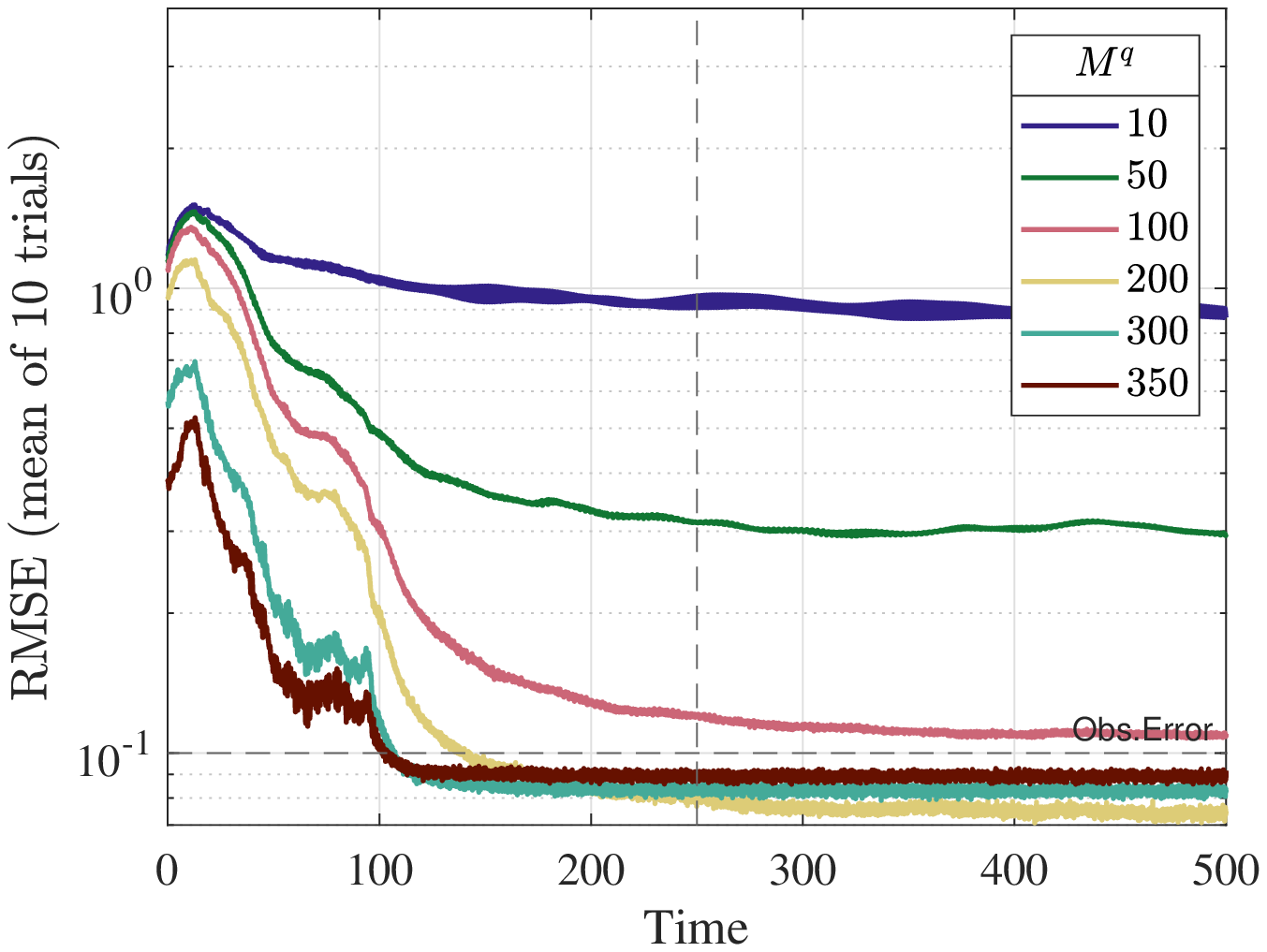}
  \caption{\(F=3\)}
\end{subfigure}
\hfill
\begin{subfigure}[t]{0.495\linewidth}\centering \includegraphics[width=\textwidth]{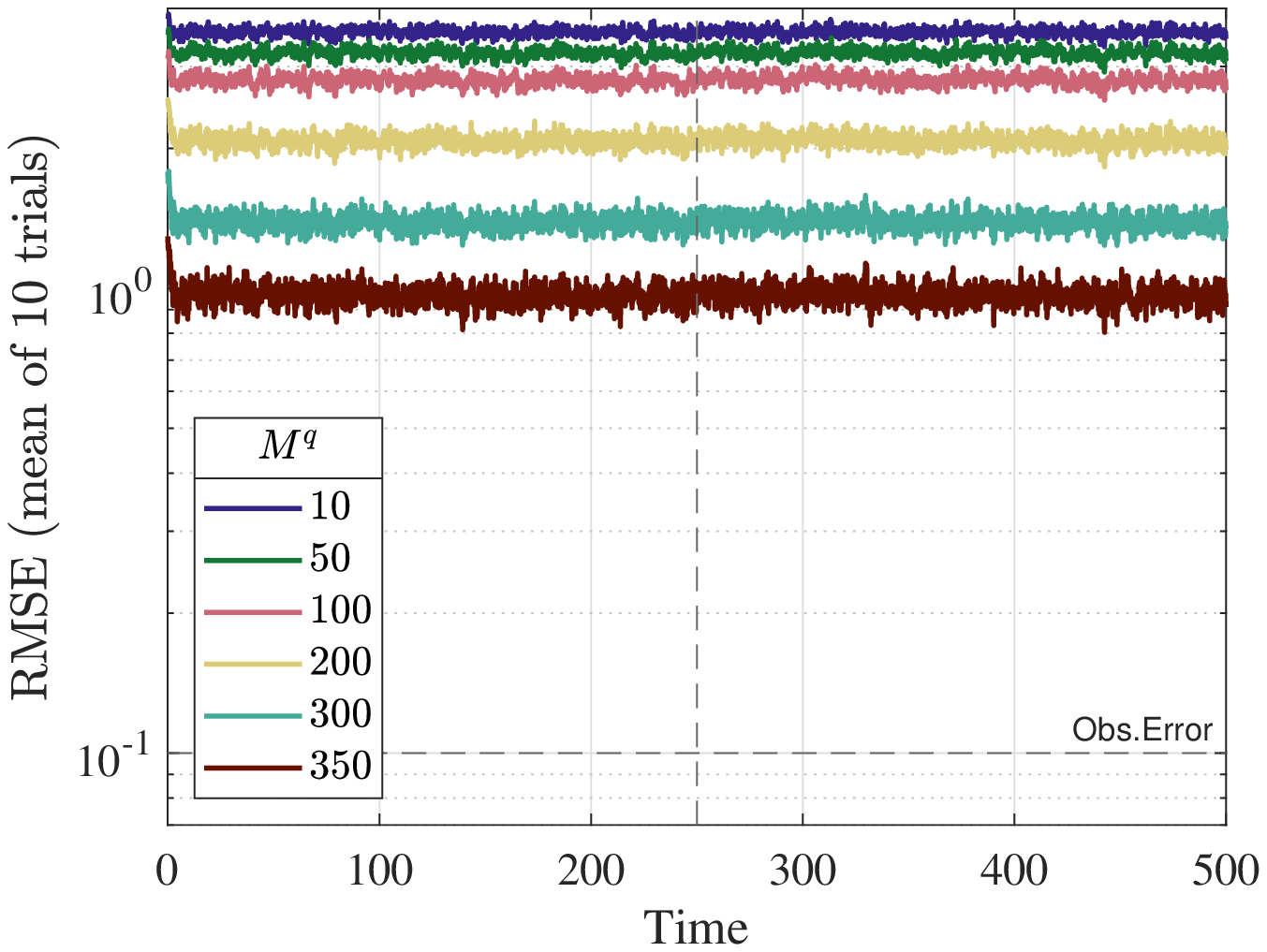}
  \caption{\(F=8\)}
  \end{subfigure}
  \caption{Time evolution of \gls{RMSE} for \gls{ProjOPPF} using \gls{POD} for the \gls{L96} model in regular \(F=3\) and chaotic \(F=8\) regime. Covariances \(\modelerrorcovariance = 0.1\,\identity, \dataerrorcovariance = 0.01\,\identity\), and the remaining parameters as in Experiment 3 (see \Cref{tab:l96-exp}). Horizontal line at \(0.1 = \sqrt{0.01}\) denotes the standard deviation of the observation noise, while vertical line at \(t=250\) indicates the start of the period used to compute mean RMSE in \Cref{L96_data} and \Cref{l96_pod_dmd} .}
  \label{fig:l96-RMSE-in-time}
\end{figure}

Time evolution of RMSE shown in \Cref{fig:l96-RMSE-in-time} for $\reduceddatadimension=5$ shows two representative cases of reduced-order assimilation in regular \(F=3\) and chaotic \(F=8\) regimes.
For \(F=3,\) the dynamics are regular, and evolve in lower-dimensional subspace of the state space, as indicated by the calculation of singular values and Lyapunov dimensions mentioned in \Cref{sec:l96-model}.
Consequently, fewer models, likely between \(\reducedmodeldimension=50\) and \(\reducedmodeldimension=100\) are needed for accurate data assimilation, as indicated by the \gls{RMSE} converging to, or below, standard deviation of observation error once sufficiently many modes are selected.
Increasing the reduced model dimension
\(\reducedmodeldimension > 100\) gives only a slight improvement in the asymptotic \gls{RMSE} and in reducing the time to reach the asymptotic value.

\begin{figure}[ht!]
  \centering
\begin{subfigure}[t]{0.495\linewidth}\centering
\includegraphics[width=\textwidth,
  trim={0 0 0 28pt},clip
]{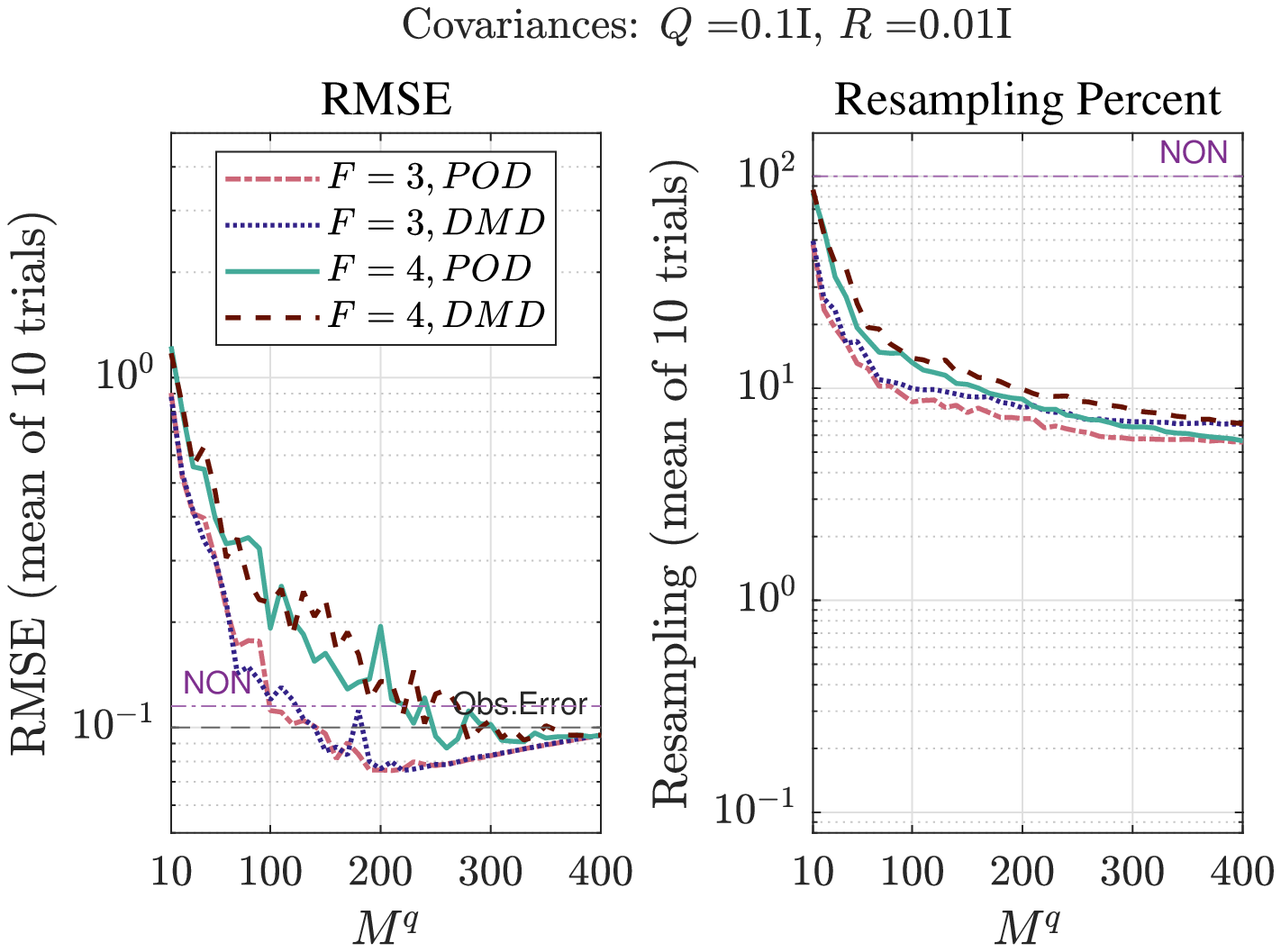}\caption{\(F = 3,4\), \(\modelerrorcovariance = 0.1\,\identity\)}
\end{subfigure}
\begin{subfigure}[t]{0.495\linewidth}\centering    \includegraphics[width=\textwidth,
  trim={0 0 0 28pt},clip
  ]{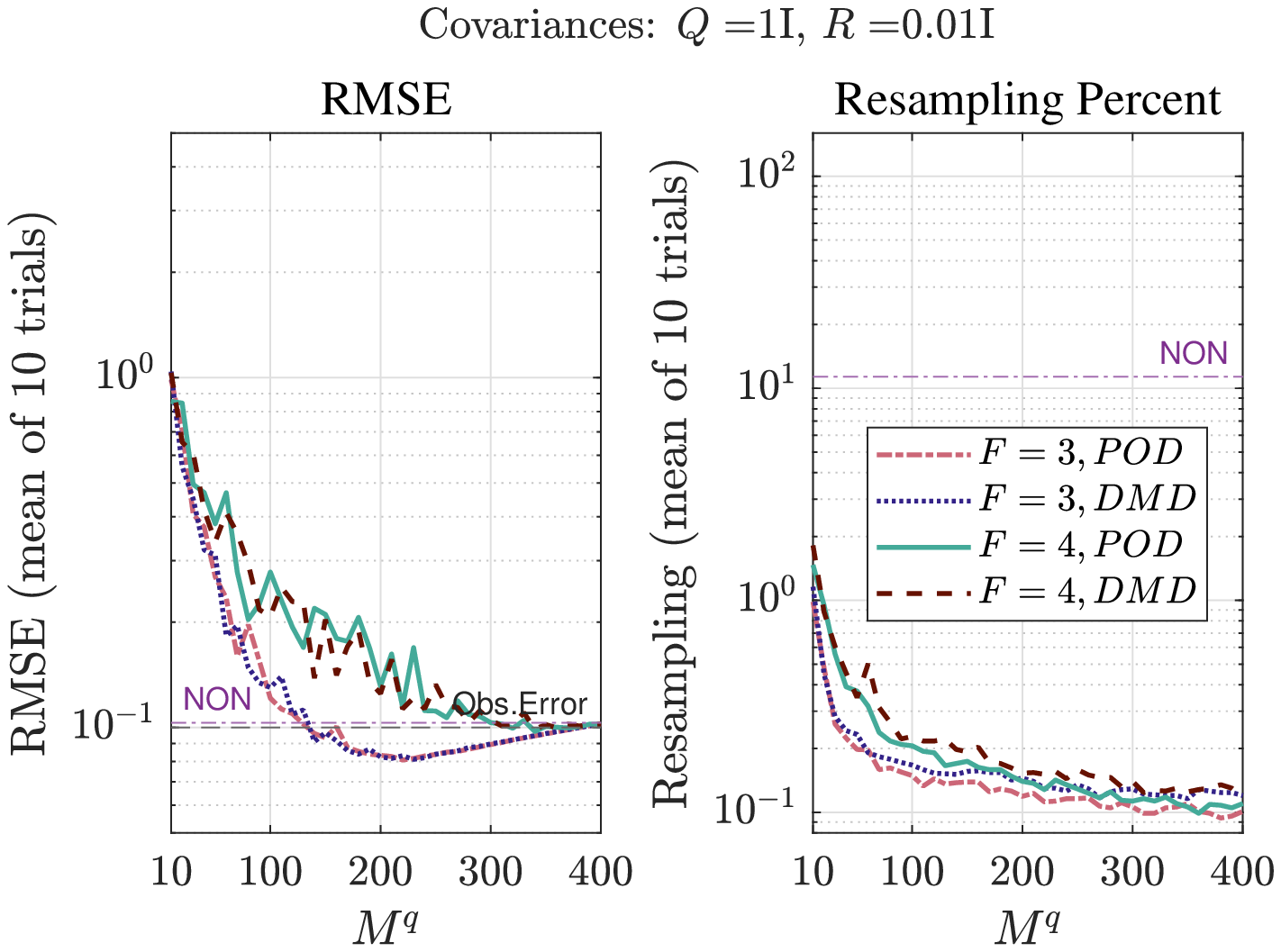}
\caption{\(F = 3,4\), \(\modelerrorcovariance = 1.0\,\identity\)}
\end{subfigure}

\begin{subfigure}[t]{0.495\linewidth}\centering    \includegraphics[width=\textwidth,
  trim={0 0 0 28pt},clip
  ]{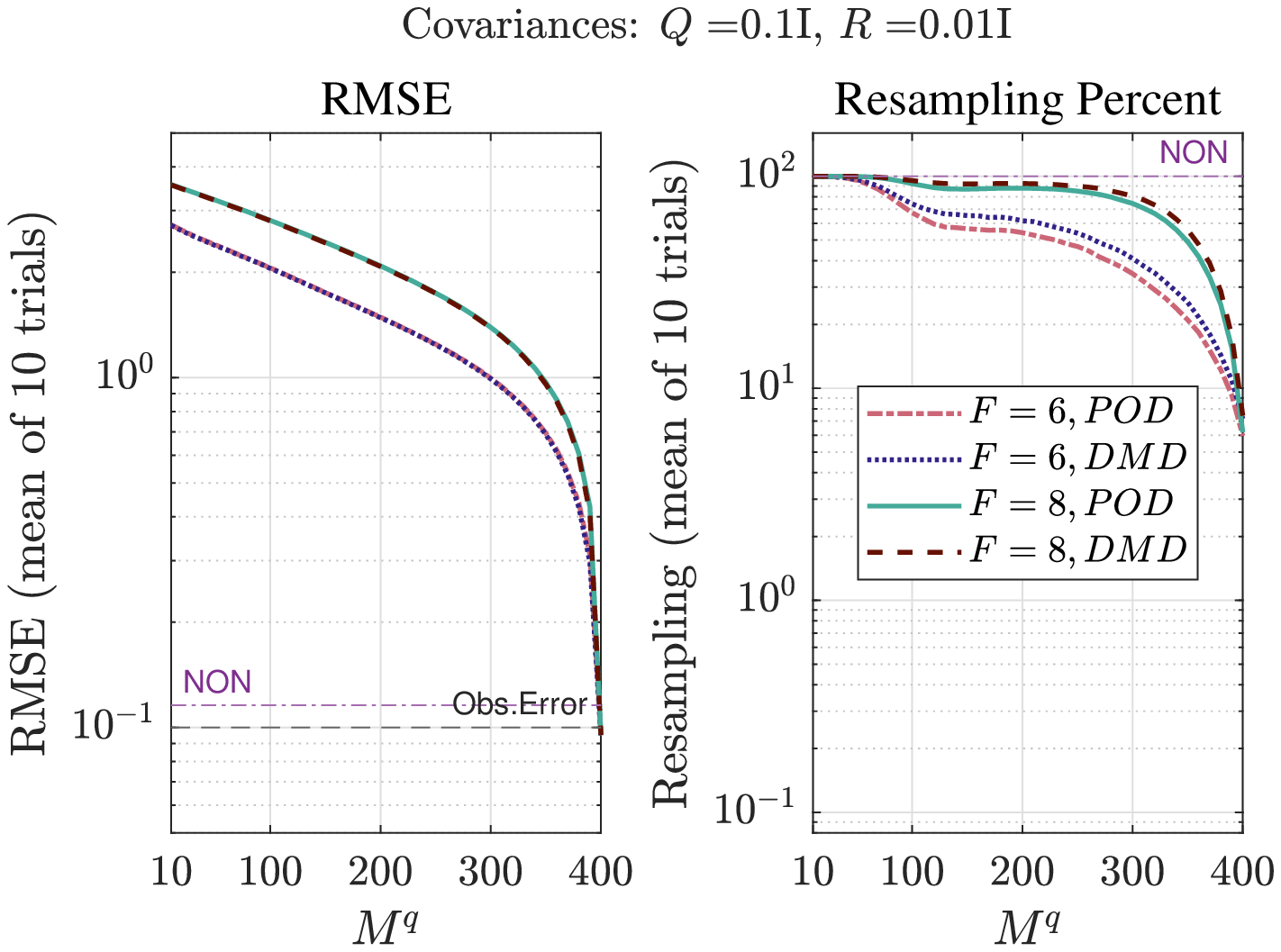}
\caption{\(F = 6,8\), \(\modelerrorcovariance = 0.1\,\identity\)}
\end{subfigure}
\begin{subfigure}[t]{0.495\linewidth}\centering  \includegraphics[width=\textwidth,
  trim={0 0 0 28pt},clip
  ]{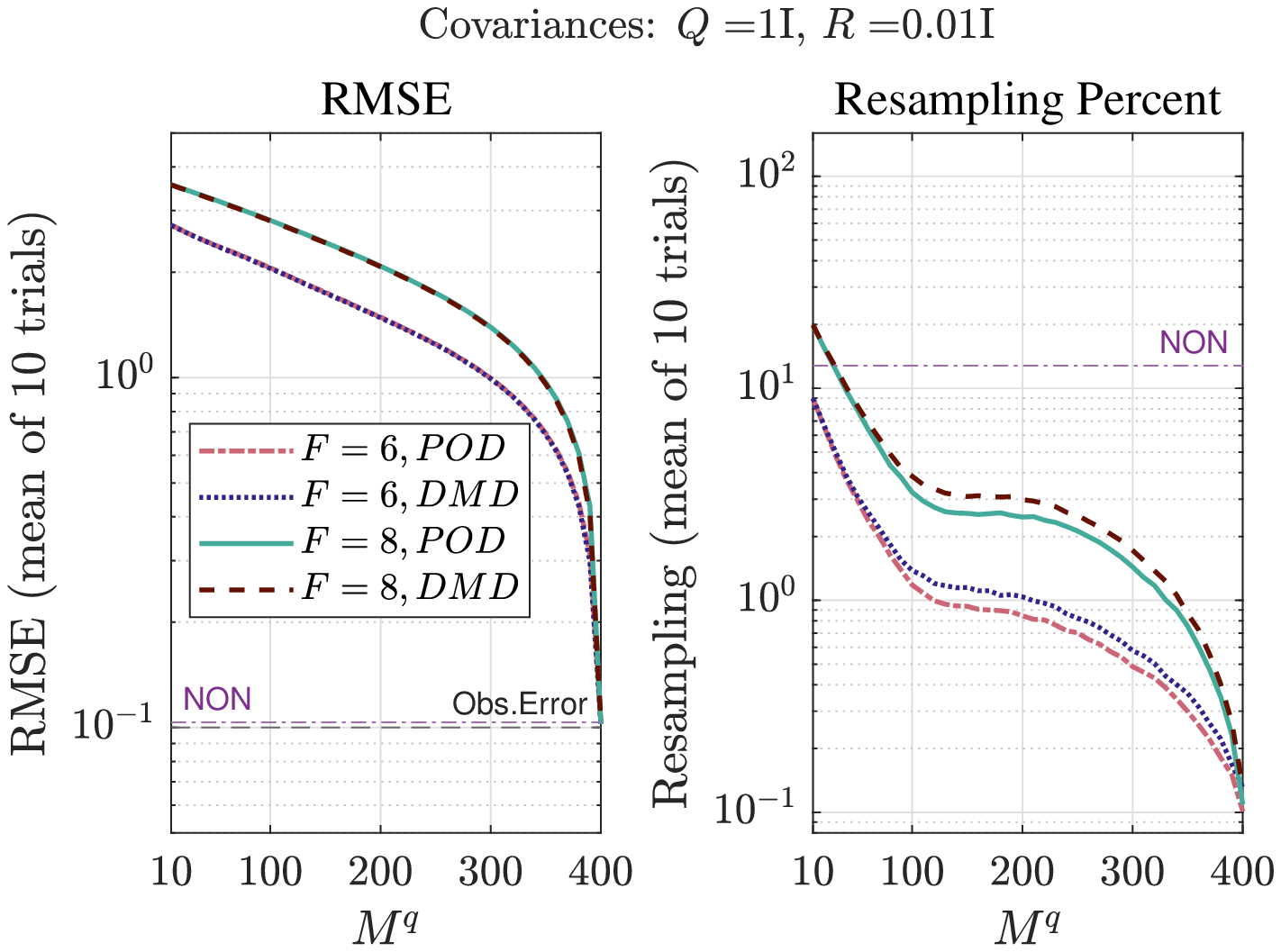}
\caption{\(F = 6,8\), \(\modelerrorcovariance = 1.0\,\identity\)}
\end{subfigure}

\caption{Experiment 3 using L96 model with $F=3,4,6,8$ and $M=400$, \gls{ProjOPPF} using  \gls{DMD}, and \gls{POD} for the physical model and the data model with fixed data projection dimension $\reduceddatadimension=5$ and varying the rank of the physical model projection $\modeldimension^q=10,20,\dots,400$.
  Columns correspond to two different model noise correlations \(\modelerrorcovariance\). NON corresponds to the optimal proposal particle filter with no reduction of the model and the data. In the case of no reduction (NON), there is little difference between $F=3$ and $F=4$ ((a) and (b)) and similarly between $F=6$ and $F=8$ ((c) and (d)). The displayed \gls{RMSE} is a time average after the initial transient (second half of the assimilation period shown in \Cref{fig:l96-RMSE-in-time}).
Standard deviation of the observation noise is fixed and given in the dash-dot line in \gls{RMSE} panels.
See \cref{sec:L96-ex3} for further details.}
  \label{l96_pod_dmd}
\end{figure}

In contrast, \(F=8\) evolution has no similar low-dimensional structure, in addition to being chaotic, in which case only small improvements are seen in \gls{RMSE} as more modes are added, but even for fairly large \(\reducedmodeldimension=350\), the \gls{RMSE} remains an order of magnitude larger than the observation error.
\Cref{l96_pod_dmd} shows trends in dependence of \gls{RMSE} and \gls{RESAMP} across more values of forcing \(F\), and additionally compares \gls{POD} and \gls{DMD} model reduction.
In all cases, \gls{RMSE} shown is a time-average of values in the second half of the assimilation period, after the transient (to the right of the vertical dashed line in \Cref{l96_pod_dmd}).

We see that for larger values of $F$ the results are of the order of the observation error only for projection ranks near the underlying model dimension $\modeldimension = 400$.
However, for $F=3$ and $F=4$ we obtain \gls{RMSE} less than one for relatively small physical model projection ranks and \gls{RMSE} of less than $0.25$ for \gls{POD} for physical model projections with ranks of approximately $50$ and higher.
Overall, as with $\modeldimension=40$ and $F=8$, we observe little gain in reduction of the physical model for $\modeldimension=400$ and $F=6,8$. On the other hand, for $F=3,4$ we observe plateauing ($F=4$) and a minimum ($F=3$) as a function of the reduced model dimension $\reducedmodeldimension$.

\glsreset{SWE}
\subsection{\gls{SWE}}

\subsubsection{Model, simulation, and parameters}

 The \gls{SWE} are frequently used in science and engineering applications to model free-surface flows where the depth is small compared to the horizontal scale(s) of the domain.
Successful applications of the SWEs include modeling dam breaks, hurricane storm surges, tsunamis and atmospheric flows \cite{Lozovskiy2017}.
Motivated by this wide utility, our second example features \gls{SWE} on a rectangular domain, configured to approximate a barotropic instability.
Detailed derivation of \gls{SWE} model can be found in standard textbooks on geophysical fluid dynamics, e.g. \cite[\S 3]{Pedlosky2013}.

The governing equations for this system are
\begin{equation}
\label{eq:swe}\begin{aligned}
  \frac{\partial u}{\partial t} &= \left(-\frac{\partial u}{\partial y}+f \right)v-\frac{\partial}{\partial x}\left(\frac{1}{2}u^2+gh \right)+\nu\nabla^2 u-c_b u, \\
  \frac{\partial v}{\partial t} &= -\left(\frac{\partial v}{\partial x}+f \right)u-\frac{\partial}{\partial y}\left(\frac{1}{2}v^2+gh \right)+\nu \nabla^2 v-c_b v, \\
  \frac{\partial h}{\partial t} &= -\frac{\partial}{\partial x}((h+\underline{h})u)-\frac{\partial}{\partial y}((h+\underline{h})v).
\end{aligned}
\end{equation}
Here  $u$ and $v$ are the \(x\)- and \(y\)- components of the   velocity field. The total height of the water column is $h+\underline{h}$, where $h$ is the height of the wave, and $\underline{h}$ is the depth of the ocean, although we employ the flat orography $\underline{h}\equiv 0$.
 The parameter $g$ is the gravitational constant, $f$ is the Coriolis parameter,
 $c_b$ the bottom friction coefficient, and $\nu$ is the viscosity coefficient.

\begin{figure}[H]
   \centering
\begin{subfigure}[t]{0.49\linewidth}\centering
    \centering
\includegraphics[width=\textwidth]{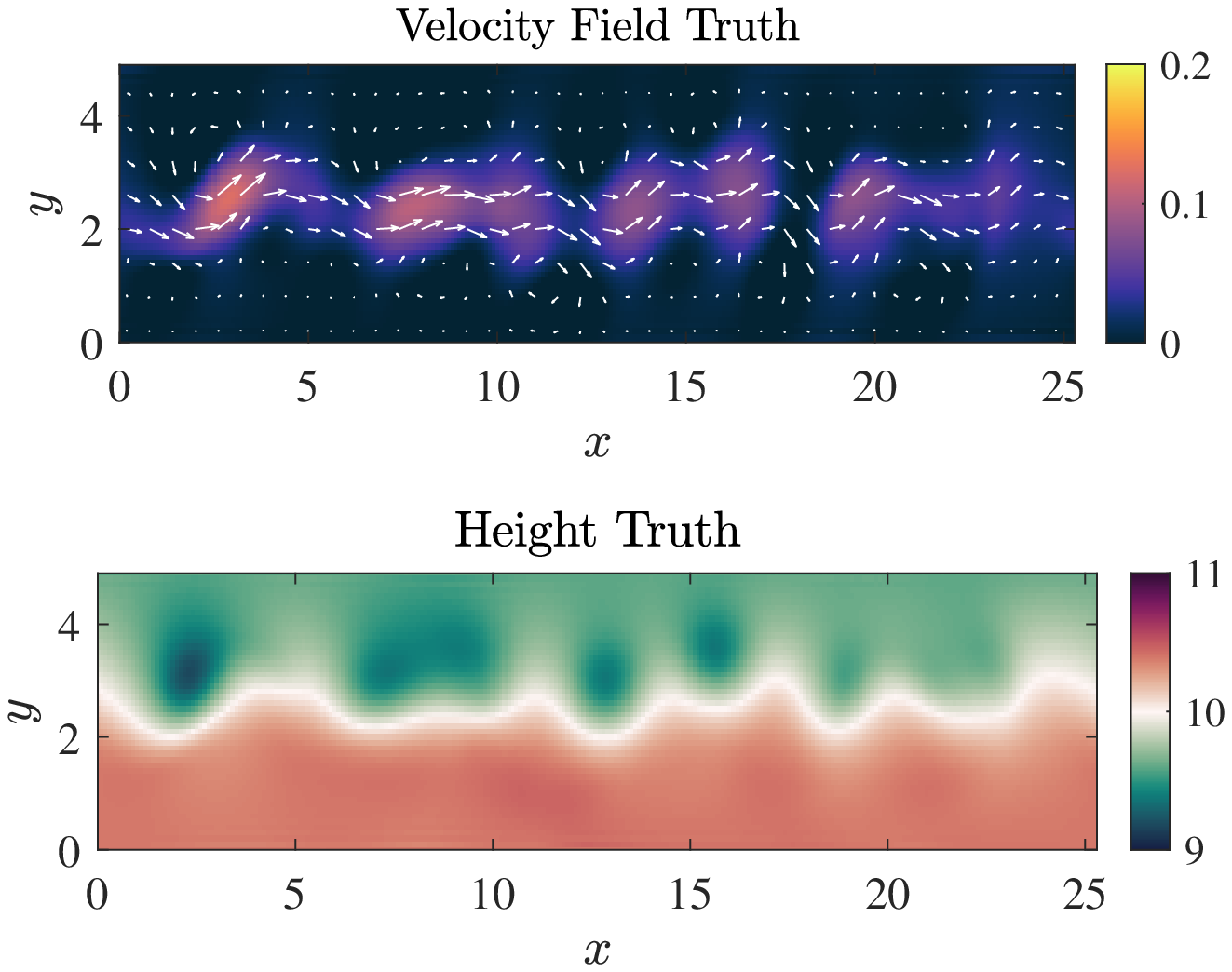}
    \caption{Truth for the velocity and height fields. Color in the top panel is local speed (magnitude of the velocity vector).}
    \label{fig:SWE_truth}
  \end{subfigure}\hfill
\begin{subfigure}[t]{0.49\linewidth}\centering
    \centering
\includegraphics[width=\textwidth]{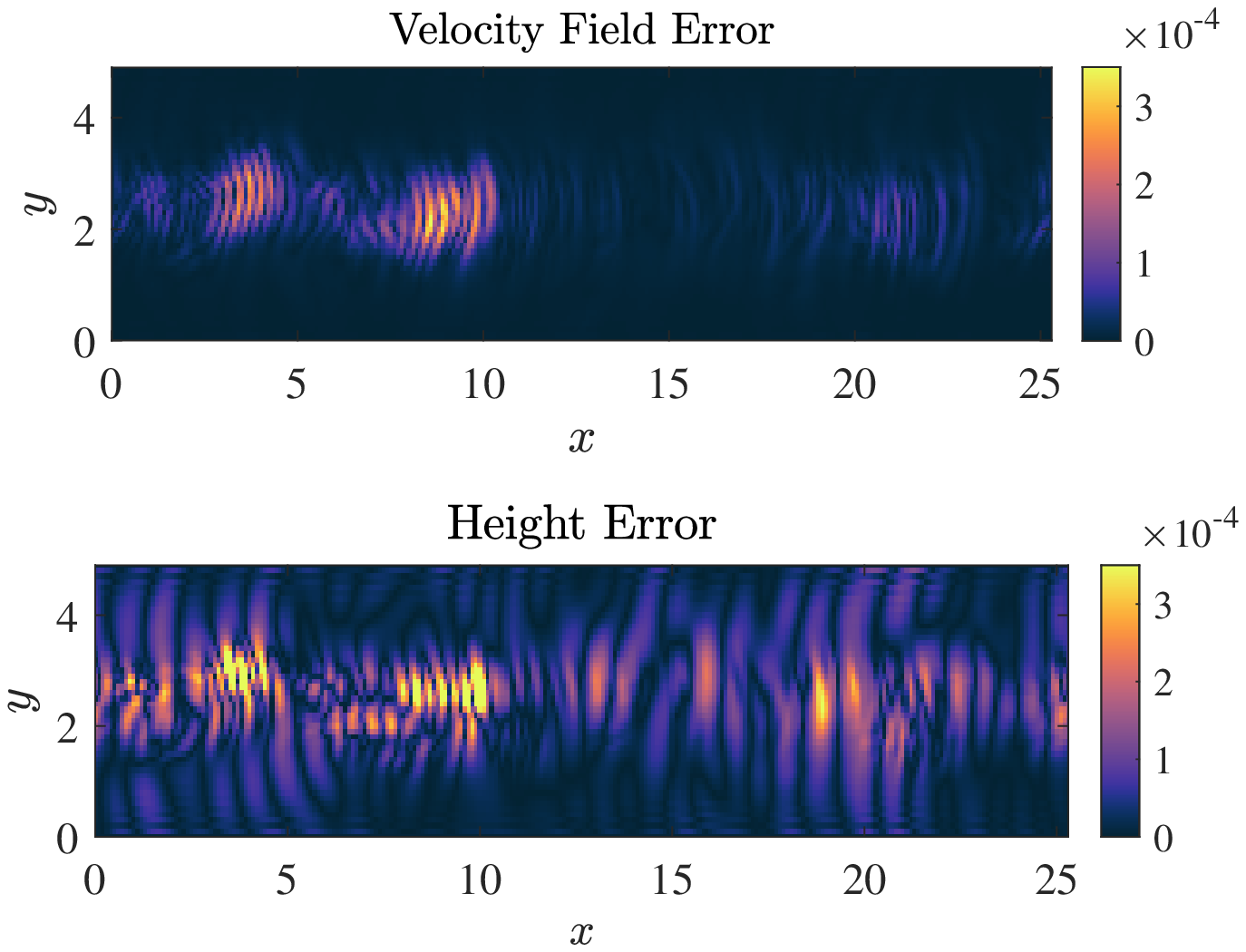}
\caption{The \(\ell^{2}\)-error between the truth and the approximation for the velocity and height fields.}
    \label{fig:SWE_difference}
\end{subfigure}
\caption{``Truth" for the state variables \(h,u,v\) of the \gls{SWE} models, and the spatial structure of the assimilation error at the end of assimilation period (\(\timeindex=72\)). Observation operator included all state variables at \(1\%\) of spatial nodes. Particle filter used \(5\) particles with error covariances \(\modelerrorcovariance=0.1\ \identity\), and \(\dataerrorcovariance = 0.01\ \identity\). Reduction performed using \gls{POD} with \(\reducedmodeldimension=20\), and \(\reduceddatadimension=10\) }
\label{fig:SWE-truth-error}
\end{figure}

The initial value problem for~\eqref{eq:swe} is solved using  R.~Hogan's finite difference code~\cite{HoganCode}.
The three fields \(u,v,h\) are evaluated at a grid of \(254 \times 50\) points in the \((x,y)\) plane, resulting in $\modeldimension = 38 100$ state variables.
Solutions are approximated using a standard finite difference scheme in space and a Lax--Wendroff finite-difference scheme in time with a fixed time step of  \(\timestepmod = 1 \si{\minute} = 1/60 \si{\hour}\) so that the Courant--Friedrichs--Lewy condition is satisfied.
We let the system evolve in time over a total of \(96 \si{\hour}\). \Cref{fig:SWE_truth} shows an example of the non-assimilated simulation output showing the barotropic instability with the flat orography.

For \gls{SWE} we consider both complete observations with all model variables observed (\(\nodespct =  100\%\)) or sparse observations with every $100$th variable observed (\(\nodespct =  1\%\)).
We additionally consider three basic observation scenarios for experiments (see \cite{Paulin2019}):
\begin{compactenum}[(i)]
\item only $u$ and $v$ variables are observed, so \(\datadimension = \frac{2}{3} \nodespct \modeldimension\) \label{scen:u-v}
\item all variables are observed, so \(\datadimension = \nodespct \modeldimension\), and \label{scen:all}
\item only $h$ variable is observed, so \(\datadimension = \frac{1}{3} \nodespct \modeldimension\). \label{scen:h}
\end{compactenum}
Since in all cases the observed variables are not transformed in any other way, the \(\dataoperator\) amounts to an identity matrix with a portion of the rows removed.

As with \gls{L96},  model and observation error noise are uncorrelated and affect all variables equally, so correlation matrices \(\modelerrorcovariance, \dataerrorcovariance\) are chosen to be scaled identity matrices.
Specifically, we employ error covariance matrices
\(\modelerrorcovariance= 1.0\ \identity, 0.1\ \identity\) and $\dataerrorcovariance = 0.01\ \identity$ and subsequently (see \ref{sec:SWE-ex-covariances}) observation error covariance matrices \(\dataerrorcovariance = 1.0\ \identity, 0.1\ \identity\).

\begin{figure}[H]
\centering
  \begin{subfigure}[t]{0.32\linewidth}    \includegraphics[width=\linewidth]{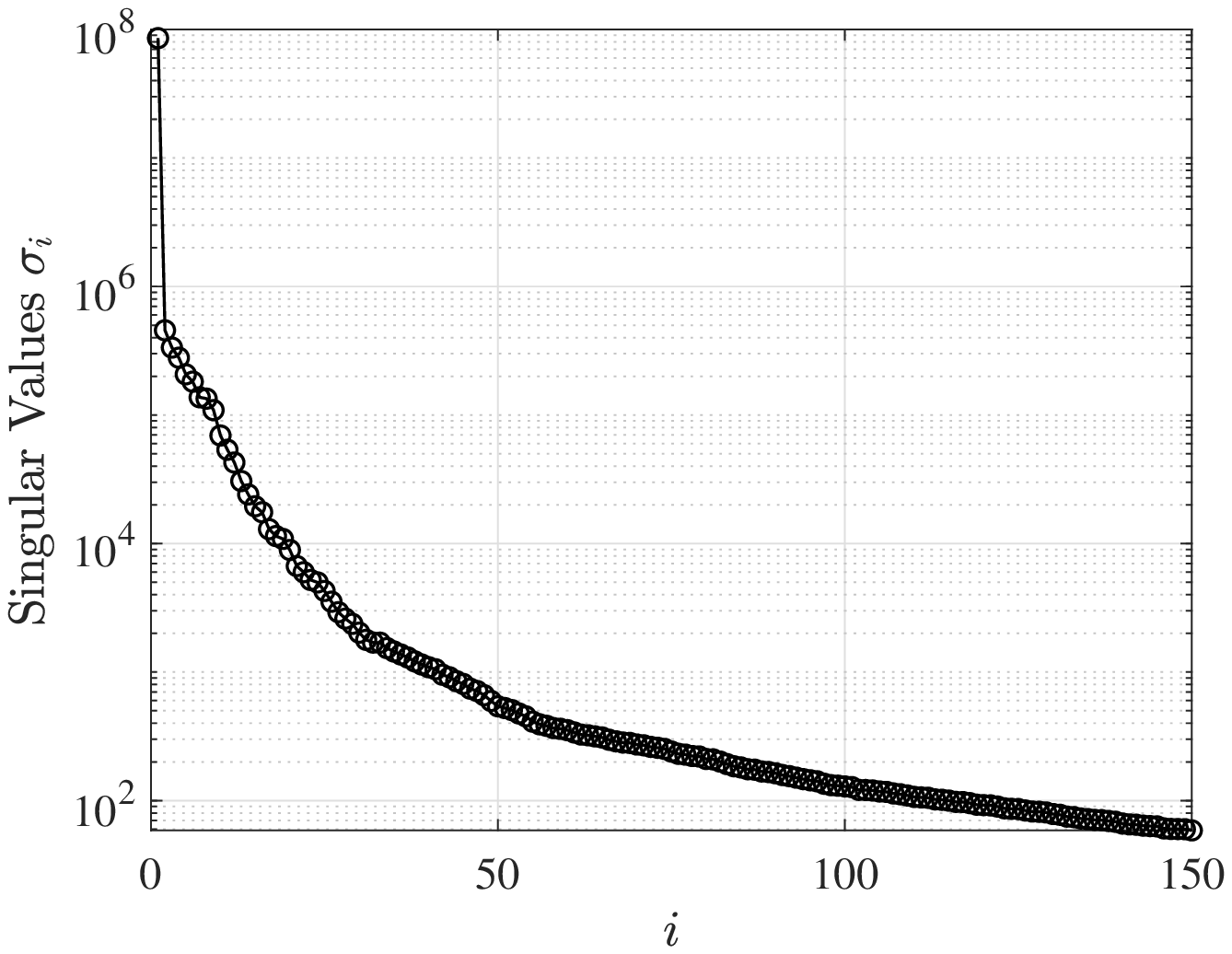}
  \caption{Singular values associated with \gls{POD} modes.}
  \end{subfigure}
  \begin{subfigure}[t]{0.32\linewidth}
\includegraphics[width=\linewidth]{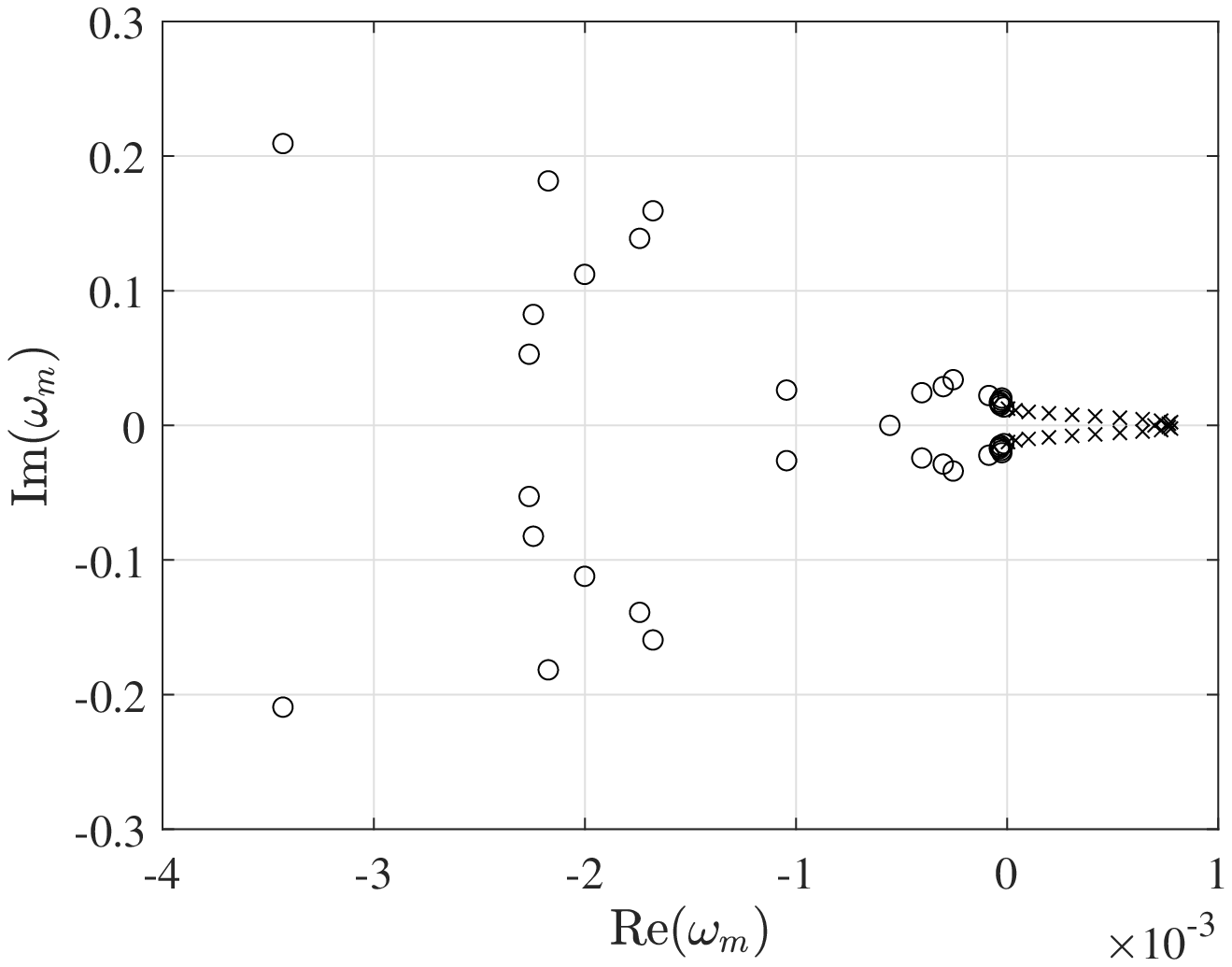}
\caption{\gls{DMD} eigenvalues in the complex plane.
\(\re (\dmdfrequency) < 0\) signify decaying modes, \(\im (\dmdfrequency)\) is proportional to the frequency of oscillations.}
  \end{subfigure}
    \begin{subfigure}[t]{0.32\linewidth}
\includegraphics[width=\linewidth]{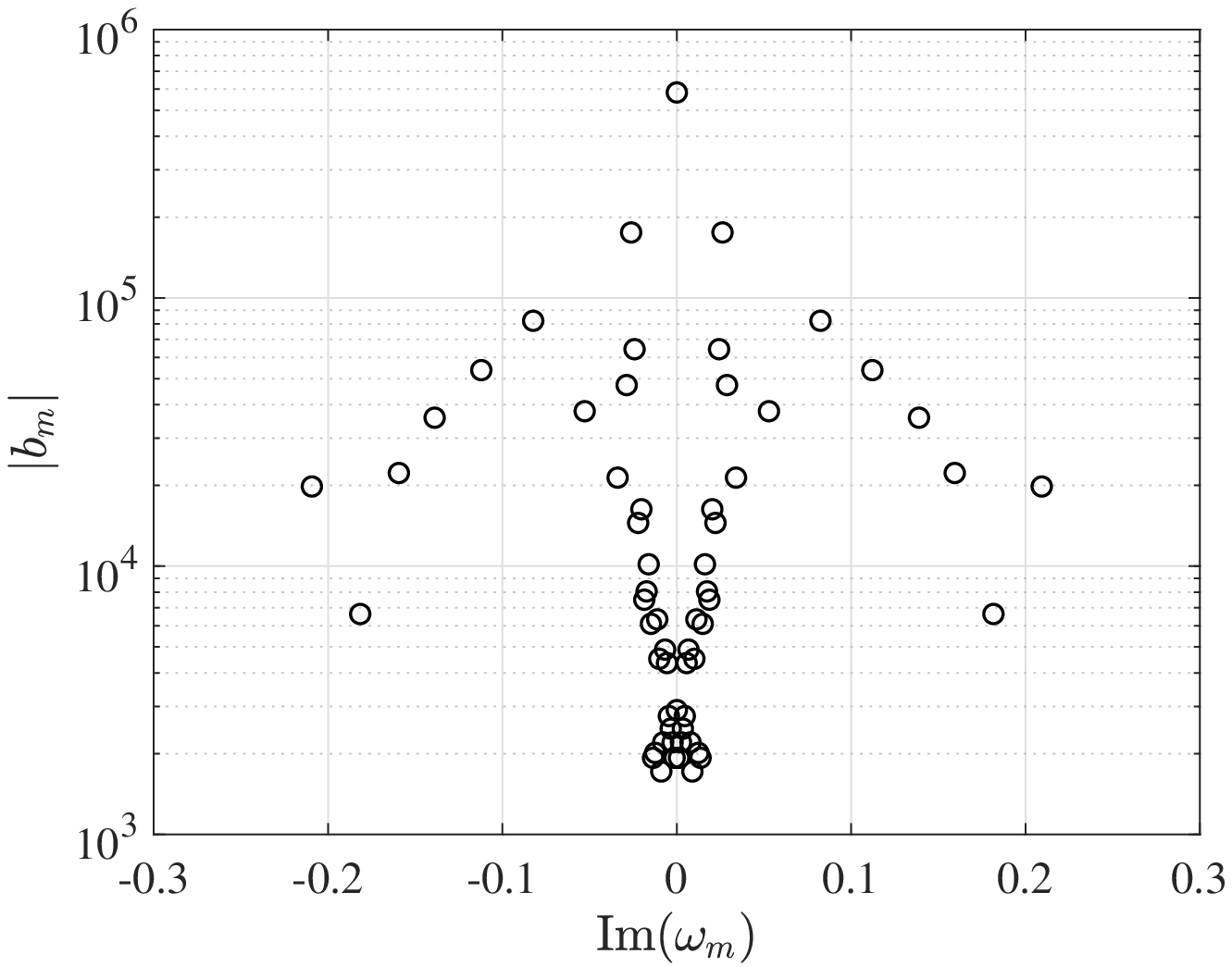}
  \caption{\gls{DMD} coefficients~\eqref{eq:DMD-coefficients}.}
  \end{subfigure}
\caption{Singular values and DMD eigenvalues spectra for \gls{SWE} model between \(\timeindex = 24 \si{\hour}\) and \(48 \si{\hour}\).}
\label{fig:swe-spectrum}
\end{figure}

Although we will not employ \gls{AUS} projections with \gls{SWE}, we calculated the approximate Lyapunov spectrum for the version of \gls{SWE} we are employing. In particular,
the largest computed approximate Lyapunov exponents over the 4 day time interval are relatively smaller than those for \gls{L96} and of the order $10^{-3} \si{\minute}^{-1}$.
We found $\approx 30$ positive exponents and a computed Kaplan--Yorke dimension of $\approx 110$. For \gls{POD} and \gls{DMD} type projections,
\Cref{fig:swe-spectrum} shows the spectrum of singular values and DMD frequencies.
Changes in the slope of the graph of singular values are sometimes used as a  indicator of an inherent dimensionality of the problem, with the assumption that different component phenomena, such as multiscale oscillations or noise sources, may have different slopes of variance associated with them.
For \gls{SWE} the \gls{POD} and \gls{DMD} projections are obtained using a fixed reference solution (the truth) for each of the 10 trials.

The \gls{DMD} coefficients \(\dmdcoefficient_{i}\) do not show many isolated peaks, suggesting that the dynamics is not prominently low-dimensional.
The same conclusion can be drawn from a lack of jumps or gaps in the spectrum of singular values.
Nevertheless, we can evaluate how well the reduced-order assimilation performs for various choices of $\reducedmodeldimension$.

With the projection matrices computed, we assimilate using \gls{ProjOPPF} starting at \(\timeindex = 48 \si{\hour}\) and continuing for the next \(24 \si{\hour}\), with observations performed and assimilated with the timestep of \(\timestepobs = 60 \si{\minute}\) (as compared to a one-minute observational time scale in \cite{Paulin2019}).
This effectively discards the first day as transient between the initial condition and the development of coherent structures.

The time scales, from the length of the simulation down to the observation timescale, were chosen to allow for an efficient proof-of-principle demonstration of the use of data-driven model reduction with data assimilation. We make no claim here that the same choices should be made in general when the \gls{SWE} is used as the model.
In general, the parametrization of the model, the broader context in which \gls{DA} is used, and the variation of the computed projection matrices with respect to the duration and start of the window would all be driving the choice of the chosen timescales.
We expect to further explore these issues in future work.

\Cref{fig:SWE-truth-error} shows the true state of the model at the end of the assimilation period (\(\timeindex=72\)) and the assimilation error at the same time.
Here we used \gls{POD} to reduce the dimension of the model space to $\modeldimension^q=20$ and the data space to $\datadimension^q=10$.
We observed \(\nodespct = 1\%\) of spatial nodes, used $L=5$ particles, and employed error covariances of $\modelerrorcovariance = 0.1\,\identity$ and $\dataerrorcovariance = 0.01\,\identity$.
The magnitude of error in both components several orders of magnitude smaller than the value of the state, indicating that the assimilation was successful.
The ideal error fields would show no structure, since the observation noise was homogeneous across all spatial nodes.
However, the error field in~\Cref{fig:SWE_difference} shows a striated pattern which is likely related to the structure of the first POD mode removed from the model.
The peaks in the pattern are concentrated along the midline, where the state variables change rapidly, which is expected as sensitivity to the spatial location above/below the midline would likely result in different assimilation particles taking different values of state variables at those spatial nodes.

\glsreset{ESS}

For a systematic evaluation of \gls{OP-PF} across ranges of parameters, we performed experiments in which the quality of assimilated state was measured by  \gls{RMSE} and \gls{RESAMP}, as explained in \cref{sec:exp-setup}. Details of the setup of each experiment are given in \Cref{tab:swe-exp}. The \gls{RMSE} and \gls{RESAMP} are calculated based on hourly observations in the $24\si{\hour}$ observation window. We also tested the dependence on the projected data dimension $(\reduceddatadimension)$ similar to what we illustrated in \Cref{L96_data} for \gls{L96}.
For example, for SWE with $\reducedmodeldimension=40$ with \(\nodespct=1\%\) spatial nodes observed in scenario (\ref{scen:all}) with $\modelerrorcovariance = 0.1\,\identity$ and $\dataerrorcovariance = 0.01\,\identity$ we found little variation in the \gls{RMSE} and \gls{RESAMP} as the projected observation dimension $\reduceddatadimension$ varied from $1$ to $40$.
In particular, with $L=5$ particles we found a mean \gls{ESS} of nearly $4$, no resampling, and mean \gls{RMSE} of $0.0465$ in model space and mean \gls{RMSE} of $0.0250$ in the projected model space.
\begin{table}[H]
  \centering
  \rowcolors{1}{}{white!90!black}\small
  \begin{tabular}{ >{\bfseries\normalsize}c | c | c | c | c | c | c | c | c | c }
    {Exp.} & {Obs. Scenario} & \(\nodespct\) & \(\modelerrorcovariance\)   & \(\dataerrorcovariance\) & Reduction & {\(\reducedmodeldimension\)} & {\(\reduceddatadimension\)} & \(L\)  \\ \hline
    1 & (\ref{scen:all}) & \(100\%\)  & \(0.1\ \identity\) & \(0.01\ \identity\) & DMD & 10 -- 100 & 10 & 5,15,20,30  \\
    2 & (\ref{scen:u-v}), (\ref{scen:all}), (\ref{scen:h}) & \(1\%\)  & \(0.1\ \identity\), \(1.0\ \identity\) & \(0.01\ \identity\) & POD & 10 -- 50 & 10 & 5  \\
    3 & (\ref{scen:all}) & \(1\% , 100\%\)  & \(0.1\ \identity\), \(1.0 \ \identity\) & \(0.01\ \identity\) & POD, DMD & 10 -- 100 & 10 & 5  \\
    4 & (\ref{scen:all}) & \(1\%\)  & \(0.1\ \identity\), \(1.0 \ \identity\) & \(0.1\ \identity\),\(1.0\ \identity\) & POD, DMD & 10 -- 100 & 10 & 5  \\
    \end{tabular}
    \caption{Parametrization of experiments used to test \gls{ProjOPPF} for \gls{SWE}.}
    \label{tab:swe-exp}
  \end{table}

\subsubsection{Experiment 1 (Number of Particles)}\label{sec:SWE-ex-particles}

\begin{figure}[H]
  \centering
\begin{subfigure}[t]{0.495\linewidth}\centering
  \includegraphics[trim={0 0 0 28pt},clip,width=\textwidth]{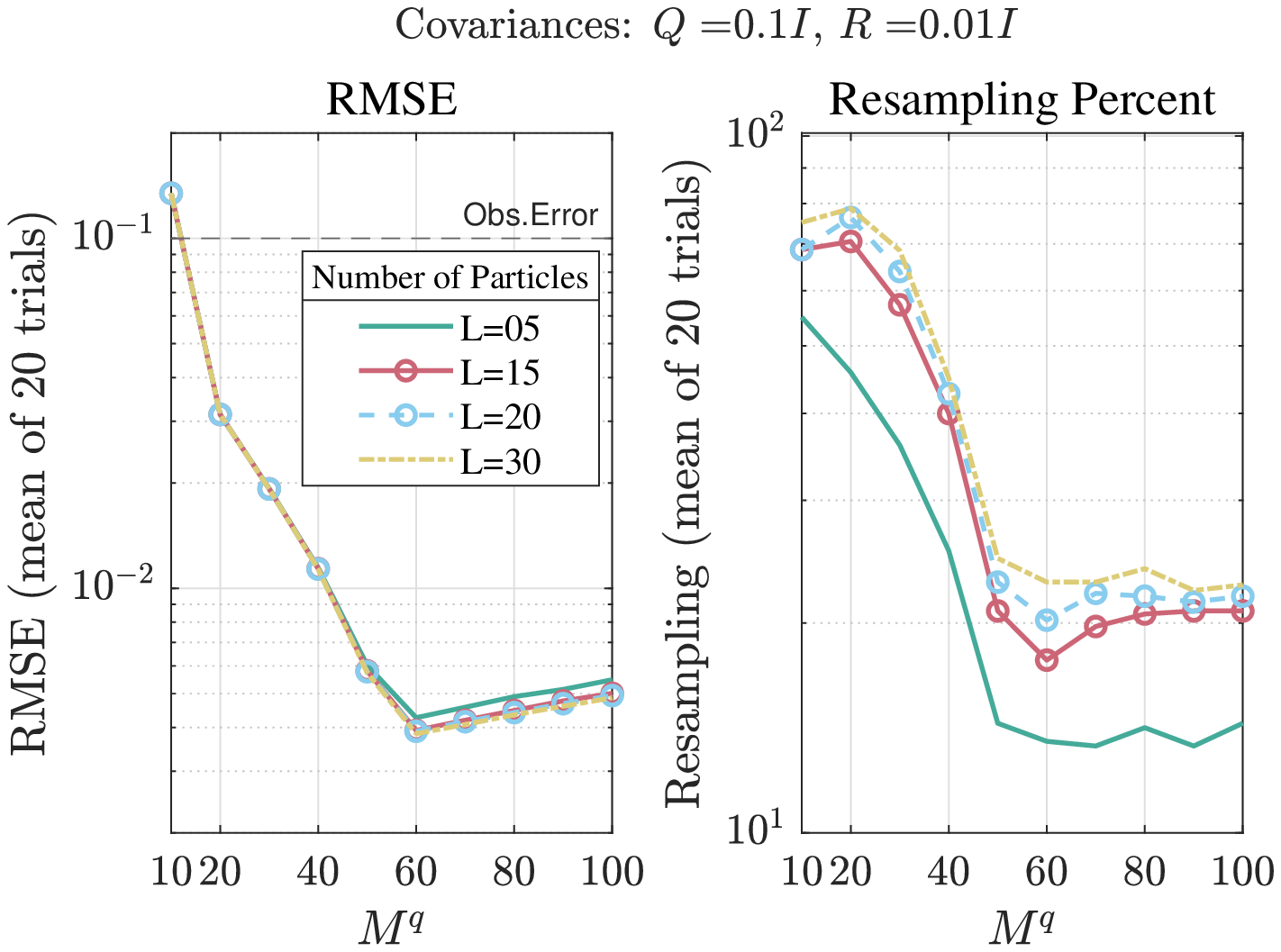}
  \caption{\(\modelerrorcovariance = 0.1\,\identity\)}\label{fig:SWEa}\end{subfigure}
\begin{subfigure}[t]{0.495\linewidth}\centering
  \includegraphics[trim={0 0 0 28pt},clip,width=\textwidth]{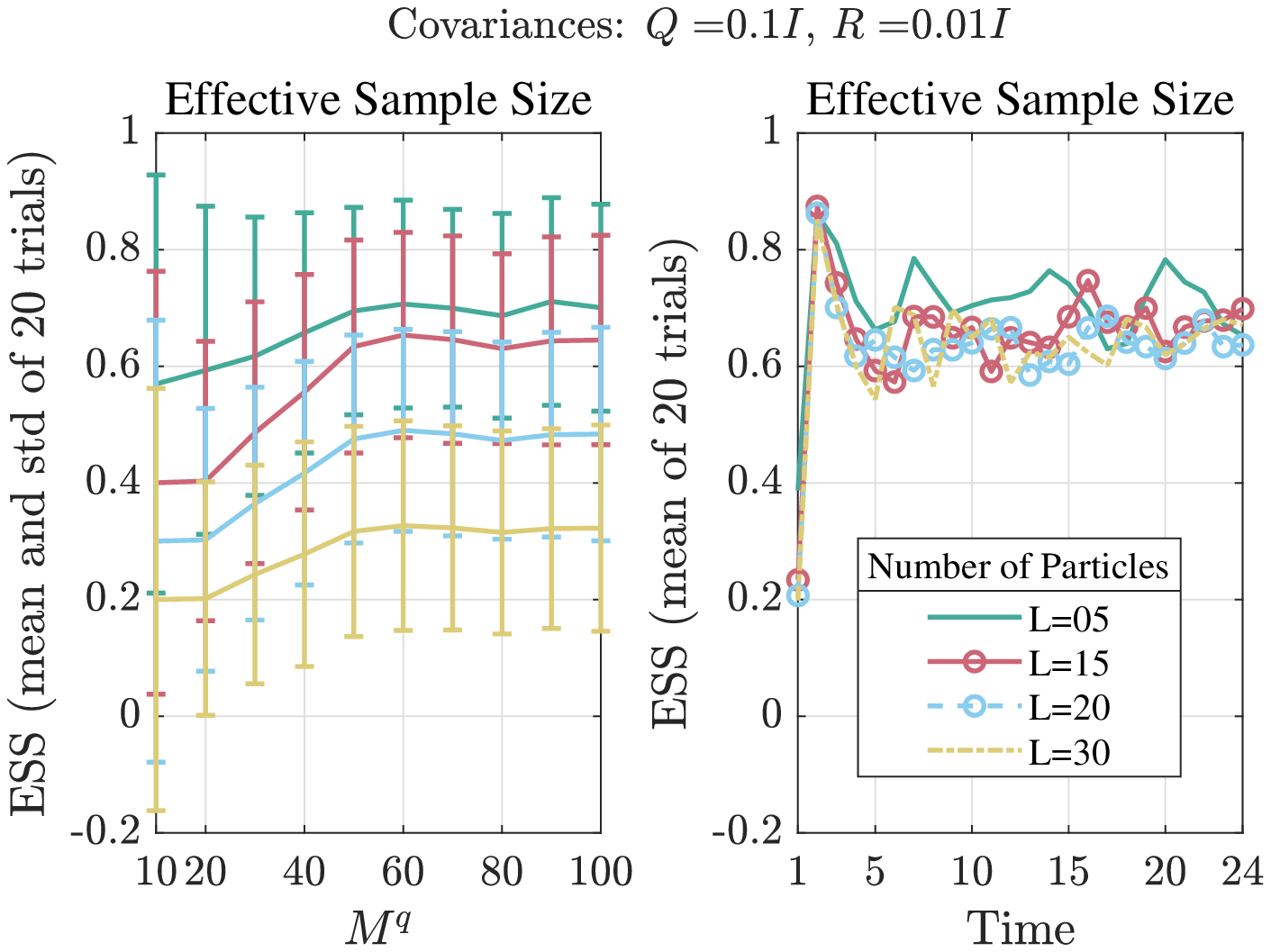}
  \caption{\(\modelerrorcovariance = 0.1\,\identity\)}\label{fig:SWEb}\end{subfigure}
\caption{\gls{SWE} model, \gls{ProjOPPF} using \gls{DMD} for the physical and data models with fixed data projection dimension $\reduceddatadimension=10$ and varying the physical model projection from $\reducedmodeldimension=10$ to $\reducedmodeldimension=100$.
Observation scenario (\ref{scen:all}) with number of particles $L=5,15,20,30$. In (b) the mean scaled \gls{ESS} (\gls{ESS}/$\totalparticles$) is plotted against the reduced model dimension $\reducedmodeldimension$ (error bars correspond to a single standard deviation) and against time.}
  \label{SWE_dmd_L}
\end{figure}

In \Cref{SWE_dmd_L}, we vary the number of particles $\totalparticles$ and find that we obtain similar results for $\totalparticles=5, 15, 20, 30$ particles. Figure (\ref{fig:SWEa}) represents the \gls{RMSE} and \gls{RESAMP} where we obtain minimum \gls{RMSE} for reduced model dimension of $\reducedmodeldimension = 60$. Figure (\ref{fig:SWEb}) shows the scaled \gls{ESS} (\gls{ESS}/$\totalparticles$) where the left graph is showing the mean and standard deviation of $20$ trials for $\reducedmodeldimension=10, 20, \dots, \modeldimension$. Note that although the scaled \gls{ESS} is in general larger for $\totalparticles=5$ particles, in an absolute sense for $L=15,20,30$ the \gls{ESS} is much larger than with $\totalparticles=5.$  The right graph is showing the mean of \gls{ESS} of $20$ trials over time for $\totalparticles=5, 15, 20, 30$.
For the remaining experiments with \gls{SWE} we employ $L=5$ particles since it is computationally more efficient.

\subsubsection{Experiment 2 (Observation Scenarios)}\label{sec:SWE-ex-scenarios}

\begin{figure}[H]
  \centering
\begin{subfigure}[t]{0.495\linewidth}\centering
  \includegraphics[trim={0 0 0 28pt},clip,width=\textwidth]{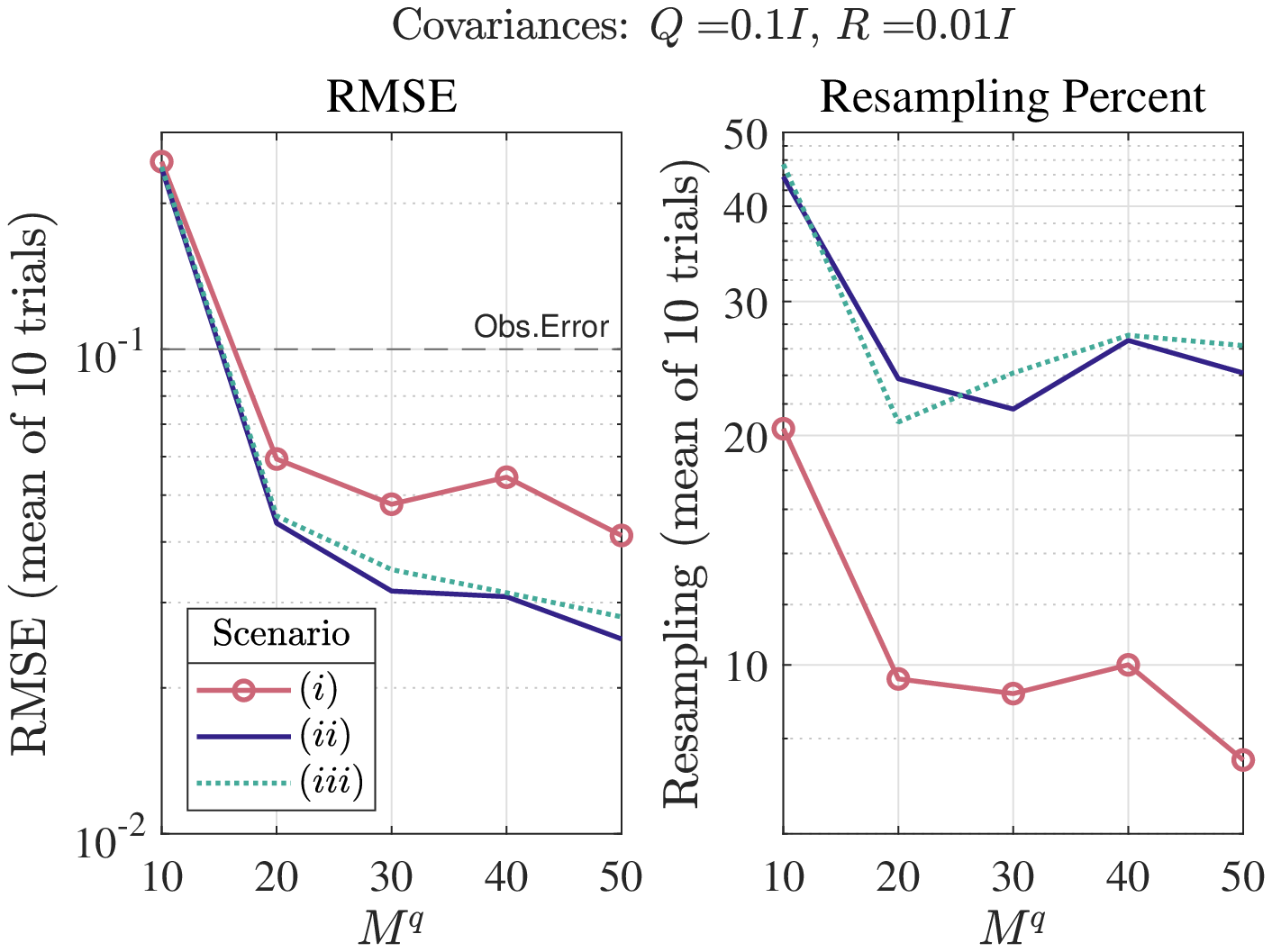}
  \caption{\(\modelerrorcovariance = 0.1\,\identity\)}
\end{subfigure}
\begin{subfigure}[t]{0.495\linewidth}\centering
  \includegraphics[trim={0 0 0 28pt},clip,width=\textwidth]{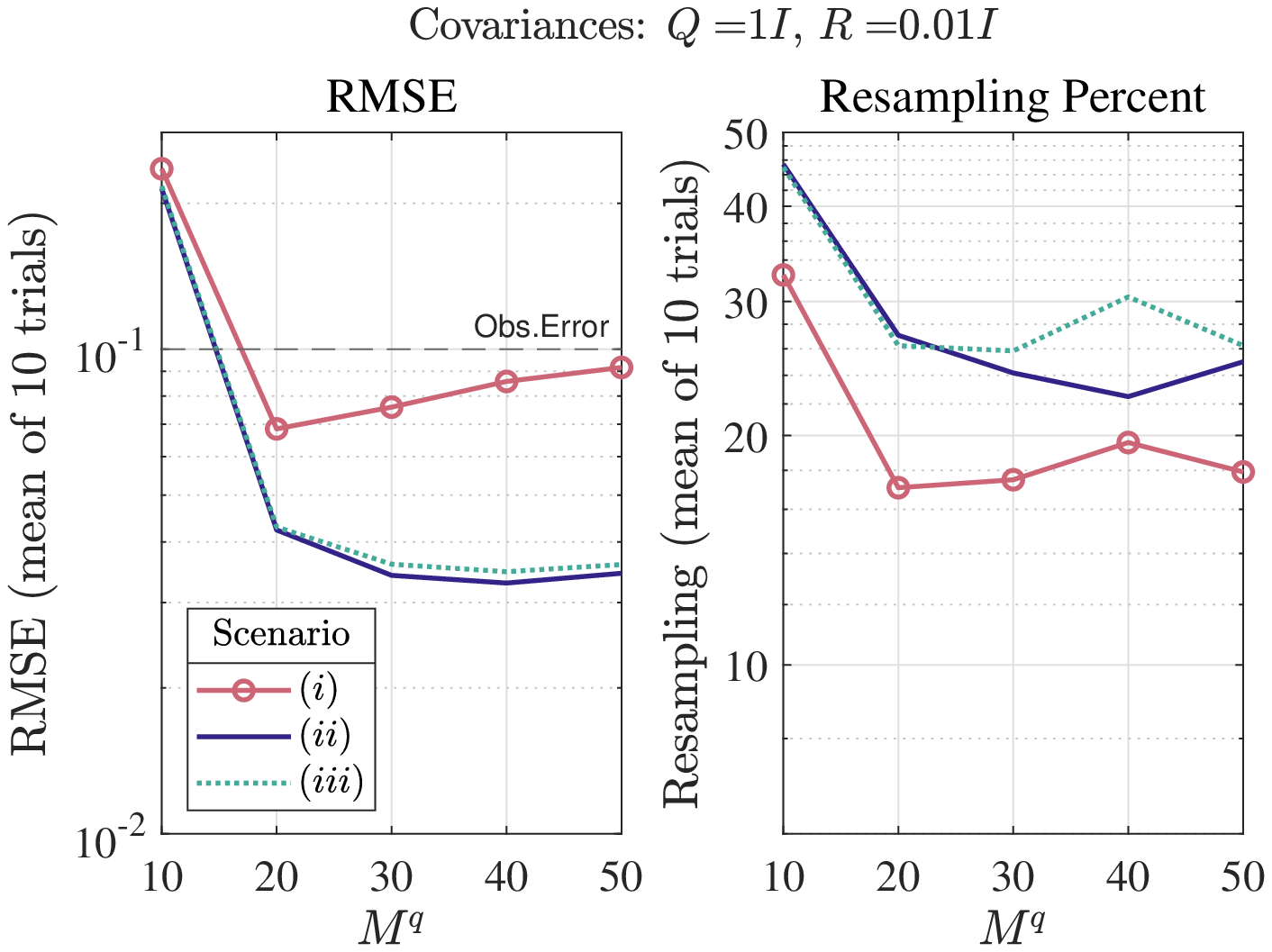}
    \caption{\(\modelerrorcovariance = 1.0\,\identity\)}
\end{subfigure}
  \caption{\gls{SWE} model, \gls{ProjOPPF} using \gls{POD} for both the physical and data models with fixed data projection rank $\reduceddatadimension=10$ and varying the rank of the physical model projection \(\reducedmodeldimension\). In all cases, assimilation uses $L=5$ particles.}
  \label{Scenarios}
  \end{figure}

\Cref{Scenarios}, illustrates the performance of \gls{OP-PF} on \gls{SWE} where \gls{POD} is used to derive the reduced order physical and data models.
It shows the mean \gls{RMSE} the resampling parentage for each of the three scenarios considered here.
Scenarios (\ref{scen:all}) and (\ref{scen:h}), in which variables $u$,$v$ and/or $h$ are observed, can be seen to have a lower \gls{RMSE} than the case of scenario (\ref{scen:u-v}), in which just $u$ and $v$ are observed.

\subsubsection{Experiment 3 (Sparse and Complete Observations)}\label{sec:SWE-ex-sparsity}

\begin{figure}[H]
  \centering
\begin{subfigure}[t]{0.495\linewidth}\centering
  \includegraphics[trim={0 0 0 28pt},clip,width=\textwidth]{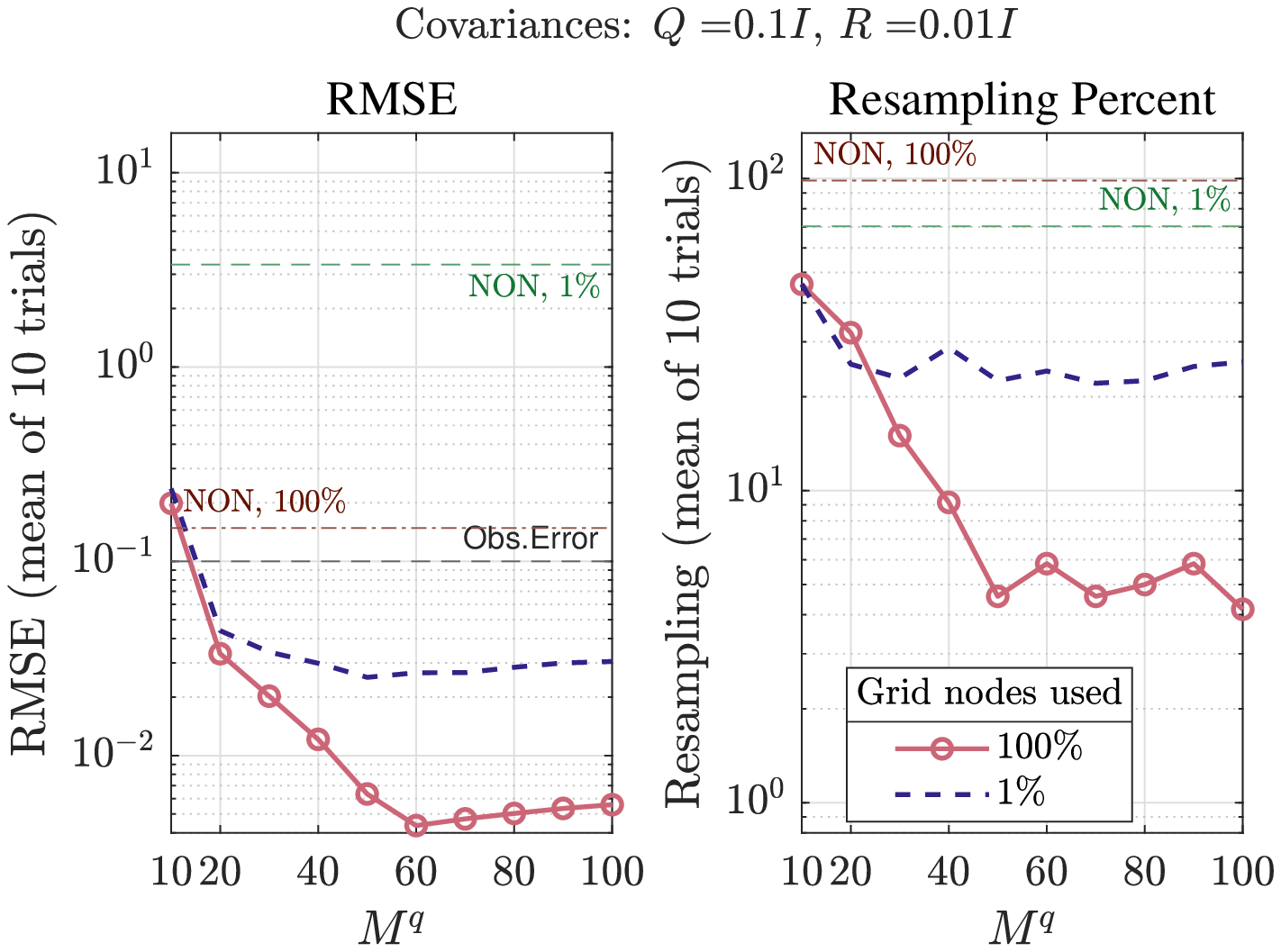}
  \caption{\(\modelerrorcovariance = 0.1\,\identity\)}\end{subfigure}
\begin{subfigure}[t]{0.495\linewidth}\centering
  \includegraphics[trim={0 0 0 28pt},clip,width=\textwidth]{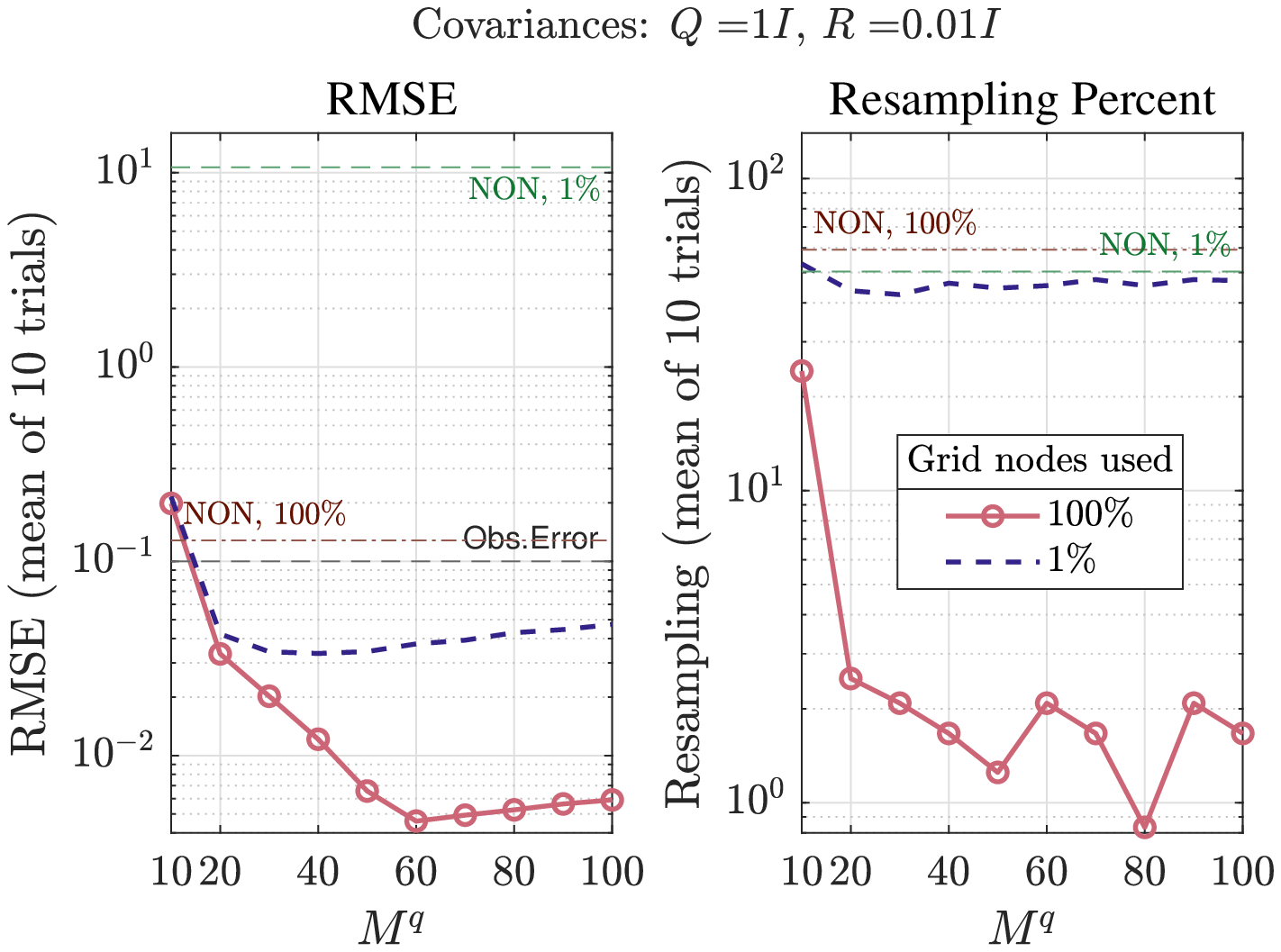}
  \caption{\(\modelerrorcovariance = 1.0\,\identity\)}\end{subfigure}
\caption{\gls{SWE} model, \gls{ProjOPPF} using \gls{POD} for the physical and data models with fixed Data projection $\reduceddatadimension=10$ and varying the physical model projection from $\reducedmodeldimension=10$ to $\reducedmodeldimension=100$. NON corresponds to the optimal proposal particle filter with no reduction of the model and the data.
   All variables $(u,v,h)$ are observed (observation scenario (\ref{scen:all})) at each observation time but at varying percentages \(\nodespct = 1\%\) and \(\nodespct= 100\%\) of spatial nodes.}
  \label{SWE_Pod_inth}
  \end{figure}

\begin{figure}[H]
      \centering
\begin{subfigure}[t]{0.495\linewidth}\centering
  \includegraphics[trim={0 0 0 28pt},clip,width=\textwidth]{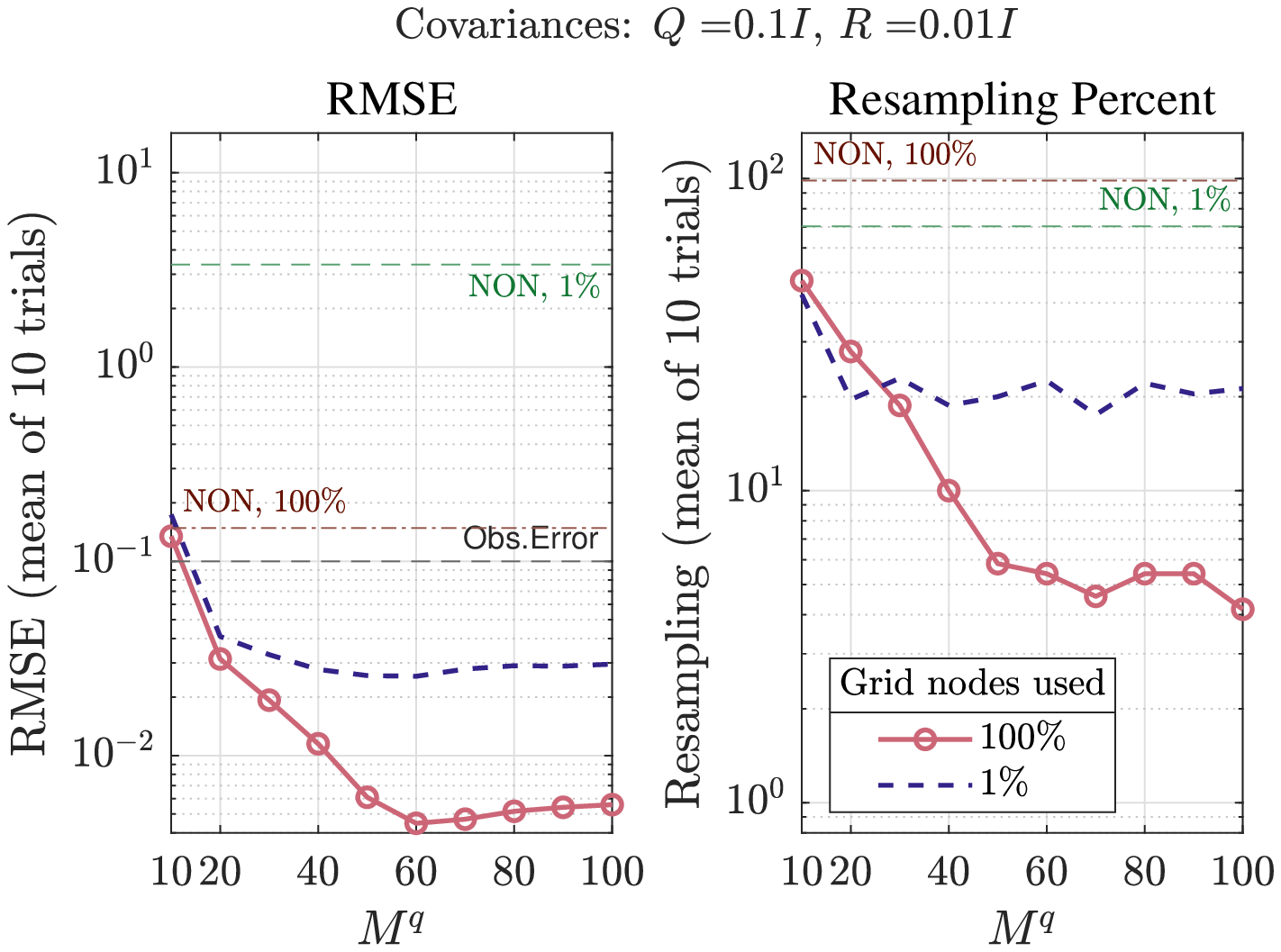}
    \caption{\(\modelerrorcovariance = 0.1\,\identity\)}
\end{subfigure}
\begin{subfigure}[t]{0.495\linewidth}\centering
  \includegraphics[trim={0 0 0 28pt},clip,width=\textwidth]{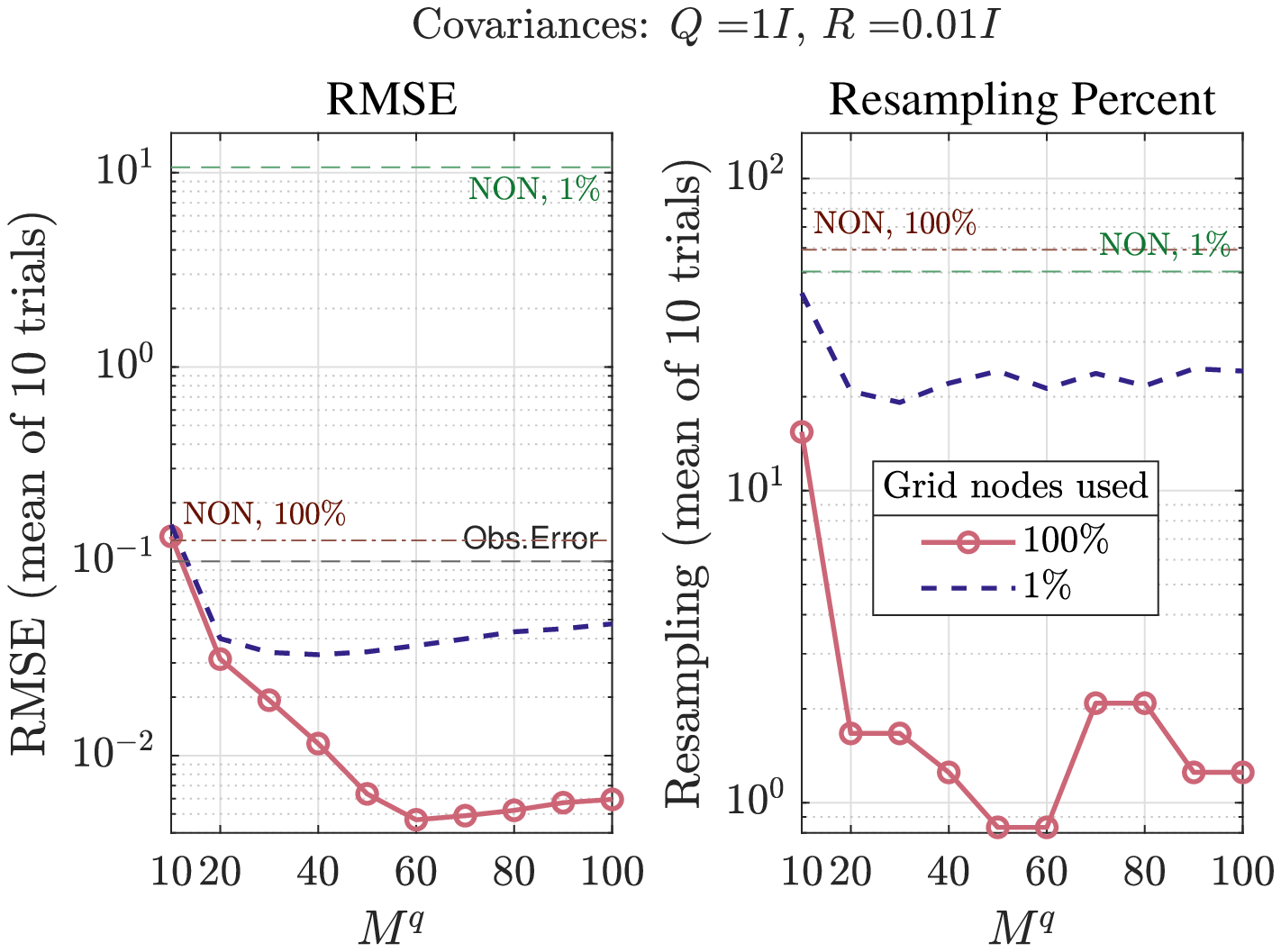}
  \caption{\(\modelerrorcovariance = 1.0\,\identity\)}
\end{subfigure}
\caption{\gls{SWE} model, \gls{ProjOPPF} using \gls{DMD} for the physical and data models with fixed Data projection $\reduceddatadimension=10$ and varying the physical model projection from dimension $\reducedmodeldimension=10$ to $\reducedmodeldimension=100$. NON corresponds to optimal proposal particle filter with no model and no data reduction. All variables $(u,v,h)$ are observed (observation scenario (\ref{scen:all})) at each observation time but at varying percentages \(\nodespct = 1\%\) and  \(\nodespct= 100\%\) of spatial nodes.}
  \label{SWE_DMD_inth}
  \end{figure}

  \Cref{SWE_Pod_inth,SWE_DMD_inth} compare the \gls{RMSE} and \gls{RESAMP} for the case when all variables $(u,v,h)$ are observed at only a fraction of spatial nodes, \(\nodespct = 1\%\) and  \(\nodespct= 100\%\) of nodes observed, and implying that the dimension of the observation space scales as $\datadimension =  \nodespct \modeldimension$. In all cases the minimum \gls{RMSE} obtained over $10$ trials occurs for reduced model dimension $\reducedmodeldimension = 60$. Both \gls{RMSE} and \gls{RESAMP} are much smaller compared to the optimal proposal particle filter with no model and data reduction (NON) in all cases. The resampling percentage \gls{RESAMP} is lower for \(\nodespct= 100\%\) as are the calculated \gls{RMSE}.

\subsubsection{Experiment 4 (Larger Observation Error Covariances $\dataerrorcovariance$) }\label{sec:SWE-ex-covariances}
\begin{figure}[H]
  \centering
\begin{subfigure}[t]{0.495\linewidth}\centering
  \includegraphics[trim={0 0 0 28pt},clip,width=\textwidth]{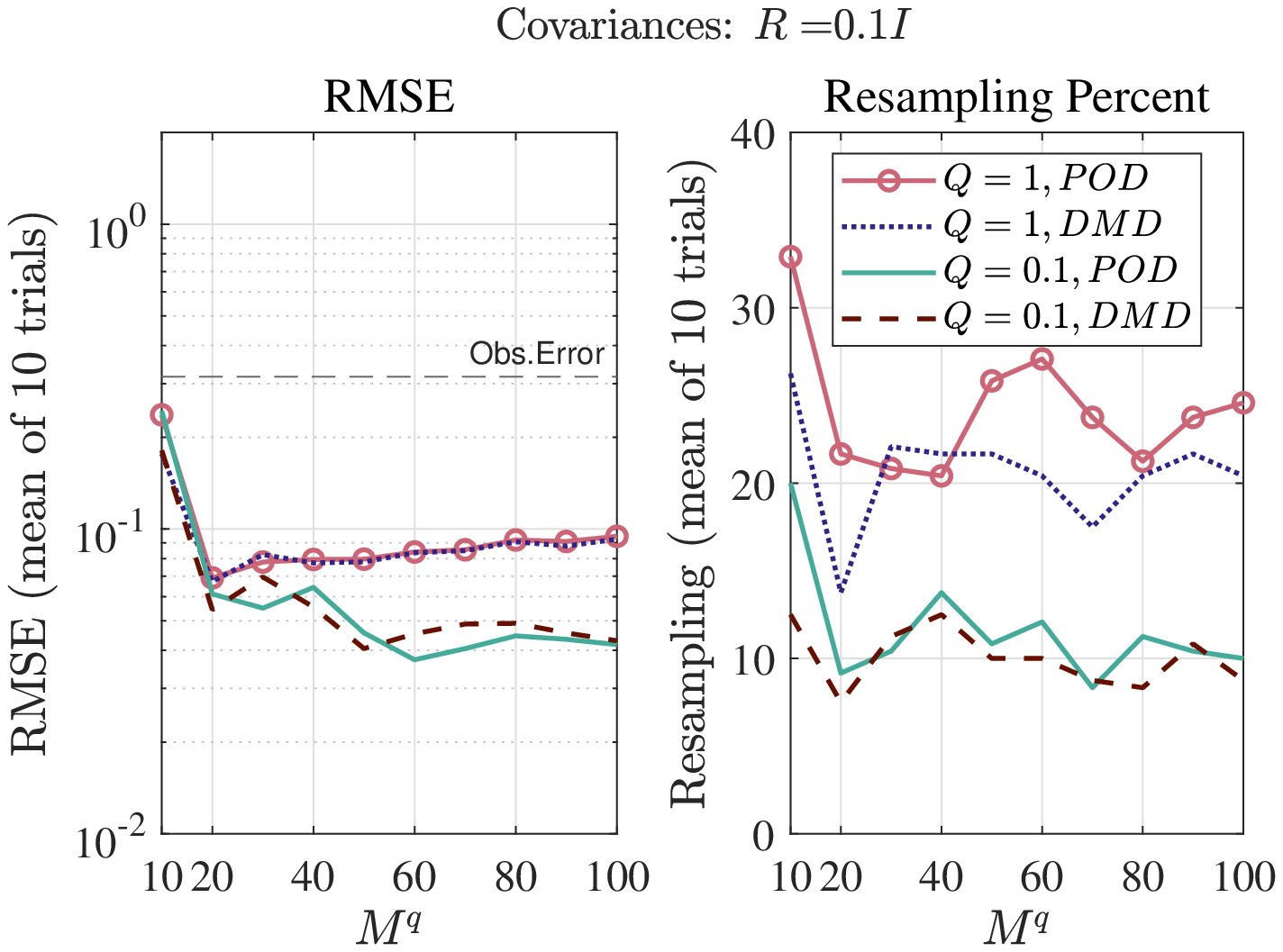}
    \caption{\(\dataerrorcovariance = 0.1\,\identity\)}
\end{subfigure}
\begin{subfigure}[t]{0.495\linewidth}\centering
  \includegraphics[trim={0 0 0 28pt},clip,width=\textwidth]{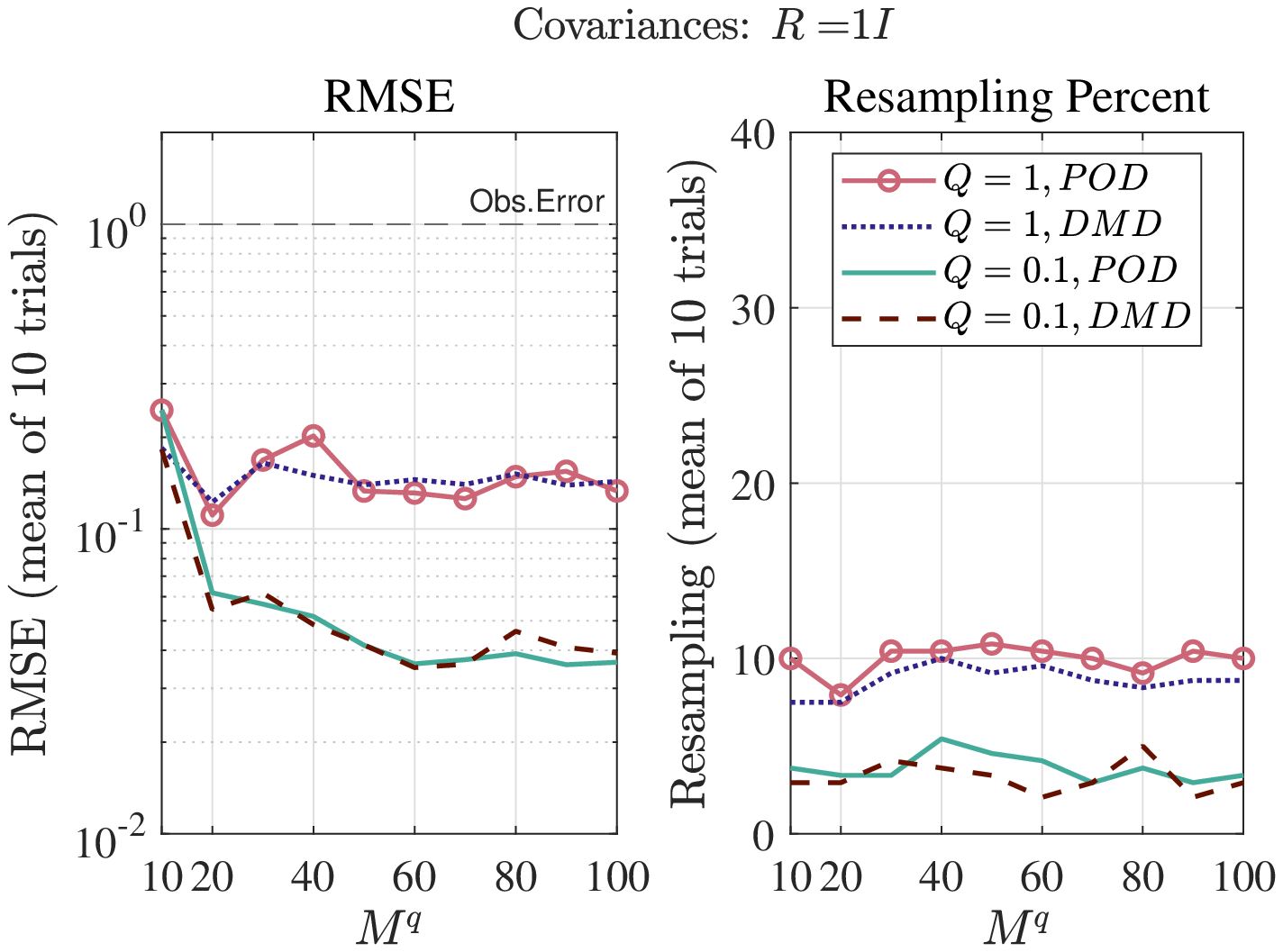}
    \caption{\(\dataerrorcovariance = 1.0\,\identity\)}
\end{subfigure}
\caption{Influence of larger observation error covariances on assimilation for the \gls{SWE} model.}
    \label{fig:SWE_R}
\end{figure}

In Figure \ref{fig:SWE_R}, \gls{POD} and \gls{DMD} are employed to reduce the dimension of the model space to $\modeldimension^q=10, 20, ..., 100$ and the data space to $\datadimension^q=10$.
We observe \(\nodespct = 1\%\) of spatial nodes and use $L=5$ particles. We compare \gls{POD} and \gls{DMD} with error covariance matrices $\modelerrorcovariance = 1.0\,\identity$ and $0.1\,\identity$ and $\dataerrorcovariance = 0.1\,\identity$ and $1.0 \,\identity$. We observe a plateauing of the \gls{RMSE} starting from approximately $\reducedmodeldimension = 60$ with model error covariance $\modelerrorcovariance = 0.1\,\identity$ and a minimum at $\reducedmodeldimension = 20$ for $\modelerrorcovariance = 1.0\,\identity$. The resampling percentages are relatively constant as a function of $\reducedmodeldimension$ with larger values for $\modelerrorcovariance = 1.0\,\identity$ and $\dataerrorcovariance = 0.1\,\identity$.

\glsresetall
\section{Discussion and Conclusions}\label{discuss_sec}

In this paper, we have derived a projected bootstrap and optimal proposal Particle Filters for physical and data models used in a standard data assimilation framework.
Our focus is on state space based projections formed using \gls{AUS}, \gls{POD}, or \gls{DMD}.
This framework provides a basis for employing Particle Filters for high dimensional nonlinear problems, and extensively tests a
projected optimal proposal particle filter algorithm, \gls{ProjOPPF}, that combines projected and unprojected models.
It is shown that stable assimilation results are obtained for the \gls{L96} model and \gls{SWE} in terms of \gls{RMSE} and resampling percentage.
The results are particularly promising for the \gls{SWE}, where \gls{ProjOPPF} with minimal tuning provides good results for severely truncated physical models and low dimensional observation operators from full to sparse observations. That is, we have successfully applied a Particle Filter to a $38,100$-dimensional multi-component physical system.
Essentially, the methods developed here perform effectively when either the physical model or the observational data have a lower effective dimension.
If these effective dimensions are sufficiently small, then the resampling percentage is low due to working with lower dimensional projected models.
When there is sufficient resolution in the reduced dimensional physical model solutions and in the reduced dimensional data, then \gls{RMSE}s are obtained on the order of the observation error.

There are several interesting avenues for further exploration.
These include the application of more sophisticated ocean, atmosphere, and coupled models.
Extension to nonlinear observation operators would require a suitable linearization of the nonlinear observation operator be obtained, to employ in the projected data models. Incorporating localization techniques for particle filters is theoretically trivial (although we have not yet tried it), as the projected approach is phrased as a Bayesian problem suitable for any data assimilation algorithm: a projected, localized Particle Filter would illustrate the value of these techniques in a more realistic context.
Another avenue we are planning to explore is to develop time dependent \gls{POD} and \gls{DMD} modes using appropriately sized moving windows of snapshots and an updating/downdating procedure.
We also plan to explore observation space projections using windowed snapshots of the observations.

\bibliographystyle{elsarticle-num}
\bibliography{references}
\appendix
\printglossaries

\end{document}